\documentclass[authoryear,noinfoline,11pt]{imsart}
\RequirePackage[OT1]{fontenc}
\RequirePackage{amsthm,amsmath,natbib}
\RequirePackage{hypernat}
\startlocaldefs
\numberwithin{equation}{section}
\theoremstyle{plain}

\endlocaldefs
\usepackage{psfrag, epsfig, graphicx}
\usepackage{amsfonts,amsmath,latexsym,amssymb}
\usepackage[active]{srcltx}
\usepackage{fullpage}

\newtheorem{lemma}{Lemma}
\newtheorem{theorem}{Theorem}
\newtheorem{proposition}{Proposition}
\newtheorem{conjecture}{Conjecture}

\newtheorem{remark}{Remark}
\newtheorem{definition}{Definition}
\newtheorem{corollary}{Corollary}

\renewcommand{\kappa}{\varkappa}

\newcommand{\rd}{{\rm d}}
\newcommand{\e}{\varepsilon}
\newcommand{\cA}{{\cal A}}

\newcommand{\cF}{{\cal F}}

\newcommand{\cH}{{\cal H}}

\newcommand{\cJ}{{\cal J}}

\newcommand{\cM}{{\cal M}}

\newcommand{\cP}{{\cal P}}

\newcommand{\cR}{{\cal R}}

\newcommand{\cX}{{\cal X}}

\newcommand{\bB}{\mathbb B}
\newcommand{\bC}{\mathbb C}

\newcommand{\bE}{\mathbb E}
\newcommand{\bF}{\mathbb F}
\newcommand{\bG}{\mathbb G}

\newcommand{\bK}{{\mathbb K}}
\newcommand{\bL}{{\mathbb L}}

\newcommand{\bN}{{\mathbb N}}

\newcommand{\bP}{{\mathbb P}}

\newcommand{\bR}{{\mathbb R}}

\newcommand{\bZ}{{\mathbb Z}}
\newcommand{\mB}{\mathfrak{B}}

\newcommand{\mm}{\mathfrak{m}}

\newcommand{\ma}{\mathfrak{a}}

\newcommand{\rk}{\mathrm{k}}

\newcommand{\epr}{\hfill\hbox{\hskip 4pt
                \vrule width 5pt height 6pt depth 1.5pt}\vspace{0.5cm}\par}

\begin{document}
\begin{frontmatter}
\title{
On adaptive minimax density estimation on~$R^d$}
\runtitle{Density estimation}
\begin{aug}
\author[t1]{\fnms{A.} \snm{Goldenshluger}
\ead[label=e1]{goldensh@stat.haifa.ac.il}}
\thankstext{t1}{Supported by the ISF grant No. 104/11.}\ \ \
\author{\fnms{O.} \snm{Lepski}
%\thanksref{t3}
\ead[label=e2]{lepski@cmi.univ-mrs.fr}}
%\and
%\author{\fnms{Third} \snm{Author}\thanksref{t1}
%\ead[label=e3]{third@somewhere.com}
%\ead[label=u1,url]{http://www.foo.com}}
%\thankstext{t1}{Some comment}
%\thankstext{t2}{Supported by the ISF grant 389/07}
%\thankstext{t3}{Supported by the ANR grant 0000}
%\thankstext{t3}{Second supporter of the project}
\runauthor{A. Goldenshluger and O. Lepski}

\affiliation{University of Haifa and  Universit\'e Aix--Marseille }

\address{Department of Statistics\\
University of Haifa
\\
Mount Carmel
\\ Haifa 31905, Israel\\
\printead{e1}\\}
%\phantom{E-mail:\ }}

\address{Laboratoire d'Analyse, Topologie, Probabilit\'es\\
 Universit\'e Aix-Marseille  \\
 39, rue F. Joliot-Curie \\
13453 Marseille, France\\
\printead{e2}\\ }
%\printead{u1}}
\end{aug}

\maketitle
\begin{abstract}
We address the problem of adaptive minimax density estimation on $\bR^d$
with $\bL_p$--loss on the anisotropic Nikol'skii classes.
We fully characterize behavior of the minimax risk for different
relationships between regularity parameters and
norm indexes in  definitions of the functional class and of the risk.
In particular, we show that there are four
different regimes  with respect
to the behavior of the minimax risk.
We develop
a single estimator which is (nearly) optimal in order
over the complete  scale of the anisotropic Nikol'skii classes.
Our estimation procedure is based
on a data-driven selection of an estimator from a fixed
family of
kernel estimators.
\end{abstract}

\end{frontmatter}

%%%%%%%%%%%%%%%%%%%%
%\setlength{\evensidemargin}{-0.5cm}
%
%\def\huh{\hbox{\vrule width 2pt height 8pt depth 2pt}}
%\def\blacksquare{\hbox{\vrule width 4pt height 4pt depth 0pt}}
%\def\square{\hbox{\vrule\vbox{\hrule\phantom{o}\hrule}\vrule}}
%\def\inter{\mathop{{\rm int}}}
%

%
\section{Introduction }
Let $X_1,\ldots, X_{n}$ be independent copies of random
vector $X\in \bR^d$ having density $f$ with respect to the Lebesgue measure.
We want to estimate $f$ using observations $X^{(n)}=(X_1,\ldots,X_n)$.
By estimator we mean any $X^{(n)}$-measurable map
$\hat{f}:\bR^n\to \bL_p\big(\bR^d\big)$. Accuracy of an estimator $\hat{f}$
is measured by the $\bL_p$--risk
\[
 \cR^{(n)}_p[\hat{f}, f]:=\Big(\bE_f \|\hat{f}-f\|_p^p\Big)^{1/p},\;\;\;p\in [1,\infty),
\]
where $\bE_f$ denotes expectation with respect to the probability measure
$\bP_f$ of the observations $X^{(n)}=(X_1,\ldots,X_n)$,
and $\|\cdot\|_p$, $p\in [1,\infty)$, is the $\bL_p$-norm on $\bR^d$. The objective is to
construct an estimator of $f$ with small $\bL_p$--risk.
\par
In the framework of the minimax approach density $f$
is assumed to belong to a functional class $\Sigma$, which is specified on the basis
of  prior information on $f$.
Given a functional class $\Sigma$,
a natural accuracy measure of an estimator $\hat{f}$ is its maximal $\bL_p$--risk
over $\Sigma$,
\[
 \cR_p^{(n)}[\hat{f};\Sigma] = \sup_{f\in \Sigma} \cR_p^{(n)}[\hat{f},f].
\]
The main question is:
\begin{itemize}
\item[(i)] how to construct
a {\em rate--optimal}, or {\em optimal in order}, estimator $\hat{f}_*$ such
that
\[
\cR_p^{(n)}[\hat{f}_*;\Sigma] \asymp
\phi_n(\Sigma):=\inf_{\hat{f}} \cR_p^{(n)}[\hat{f};\Sigma],\;\;\;\;
n\to\infty?
\]
\end{itemize}
\noindent Here the infimum is taken over all possible estimators.
We refer to the
outlined problem as the {\em problem of minimax density estimation with $\bL_p$--loss
on the class $\Sigma$}.
\par
Although the minimax approach provides a fair and convenient
criterion for comparison between different estimators, it lacks some flexibility.
Typically
$\Sigma$ is a class of  functions
that is determined by some {\em  hyper-parameter}, say, $\alpha$.
(We write $\Sigma=\Sigma_\alpha$ in order to indicate explicitly dependence of the class $\Sigma$ on
the corresponding
hyper-parameter $\alpha$.)
In general, it turns out that an estimator which
is optimal in order on the class $\Sigma_\alpha$ is not
optimal on the class $\Sigma_{\alpha^\prime}$.
This fact motivates the following question:
\begin{itemize}
 \item [(ii)] is it possible to construct an estimator $\hat{f}_*$
 that is optimal in order
on some scale of functional classes $\{\Sigma_\alpha, \alpha \in A\}$ and not only on
one class $\Sigma_\alpha$? In other words, is it possible
to construct an estimator $\hat{f}_*$ such that
for any $\alpha\in A$ one has
\[
 \cR^{(n)}[\hat{f}_*; \Sigma_\alpha] \asymp \phi_n(\Sigma_\alpha),\;\;\;\;n\to\infty?
\]
 \end{itemize}
 We refer to this question as the {\em problem of adaptive minimax density
estimation on the scale of classes
$\{\Sigma_\alpha, \alpha\in A\}$}.
\par
The minimax and adaptive minimax
density estimation with $\bL_p$--loss  is a subject of
the vast literature, see for example
\cite{bretagnolle}, \cite{Ibr-Has1},
\cite{devroye-gyorfi,dev-lug96}, \cite{Efr,Efr2}, \cite{Has-Ibr},
\cite{Donoho}, \cite{Gol}, \cite{kerk},
 \cite{rigollet},
\cite{massart}[Chapter~7],
\cite{samarov},
\cite{rigollet-tsybakov} and \cite{birge}.
It is not our aim here to provide a complete review of the literature on density
estimation with $\bL_p$-loss.
Below we will only discuss
results that are directly related to our study. First we review
papers dealing with the one--dimensional setting; then we proceed
with the multivariate case.
\par
The problem of minimax density estimation on $\bR^1$ with $\bL_p$--loss, $p\in [2,\infty)$,
was studied
by \cite{bretagnolle}.
In this paper the functional class $\Sigma$ is the class of all densities
such that
$\big[\|f^{(\beta)}\|_p \|f\|_{p/2}^{\beta}\big]^{1/(2\beta+1)}\leq L<\infty$, where $f^{(\beta)}$ is the generalized derivative of order
$\beta$.
It was shown there that
$$
\phi_n(\Sigma)\asymp n^{-\frac{1}{2+1/\beta}},\;\; \forall p\in [2,\infty).
$$
Note that the same parameter $p$ appears in the definitions of the risk and of the functional class.
\par
The problem of adaptive minimax density estimation on a compact interval of $\bR^1$
with $\bL_p$--loss was addressed in \cite{Donoho}. In this paper
class $\Sigma$ is
the Besov functional class $\bB^\beta_{r\theta}(L)$, where
parameter $\beta$ stands for the regularity index, and
$r$ is the index of the norm in which the regularity is measured.
It is shown there that there is an elbow in the rates of convergence for the minimax risk
according to whether $p\leq r(2\beta+1)$ (called in the literature {\em the dense zone}) or $p\geq r(2\beta+1)$ ({\em the sparse zone}). In particular,
\begin{equation}
\label{eq:rate-donoho}
\phi_n\big(\bB^\beta_{r\theta}(L)\big)\geq
\left\{
\begin{array}{ll}
n^{-\frac{1}{2+1/\beta}}, & p\leq r(2\beta+1),
\\
%\phi_n\big(\bB^\beta_{r\theta}(L)\big) \asymp
(\ln n/n)^{\frac{1-1/(\beta r)+1/(\beta p)}{1-1/(\beta r)+1/(2\beta)}}, &
p\geq r(2\beta+1).
\end{array}
\right.
\end{equation}
\cite{Donoho} develop a wavelet--based hard--thresholding
estimator that achieves the indicated rates
(up to a $\ln n$--factor in the dense zone) for a scale of
the Besov classes~$\bB^\beta_{r,\theta}(L)$ under additional assumption $\beta r>1$.
\par
It is quite remarkable that if the assumption that the
underlying density has compact support  is dropped, then
the minimax risk behavior becomes completely different.
Specifically,
\cite{Juditsky} studied the problem of adaptive minimax density estimation on $\bR^1$
 with $\Sigma$ being the H\"older class $\bN_{\infty,1}(\beta,L)$.
Their results are in striking contrast
with those of \cite{Donoho}: it is shown that
$$
\phi_n\big(\bN_{\infty,1}(\beta,L)\big) \geq
\left\{
\begin{array}{ll}
n^{-\frac{1}{2+1/\beta}}, & p>2+1/\beta,\\
%\phi_n\Big(\bN_{\infty,1}(\beta,L)\Big)\geq
n^{-\frac{1-1/p}{1+1/\beta}}, & 1\leq p\leq 2+1/\beta.
\end{array}\right.
$$
\cite{Juditsky} develop
a wavelet--based estimator that
achieves the indicated rates up to a logarithmic factor
on a scale of the H\"older classes.
Note that if the aforementioned
results of \cite{Donoho} for densities with compact support are
applied to the H\"older
class, $r=\infty$, then
the rate is $n^{-1/(2+1/\beta)}$ for any $p\geq 1$. Thus,
the rate corresponding to the zone $1\leq p\leq 2+1/\beta$,
does not appear in the case
of compactly supported densities.
\par
In a recent paper, \cite{patricia} consider the problem of adaptive density
estimation on $\bR^1$ with $\bL_2$--losses on the Besov classes
$\bB_{r\theta}^\beta(L)$. It is shown there
that
$$
\phi_n\big((\bB_{r\theta}^\beta(L)\big)\geq
\left\{
\begin{array}{cc}
n^{-\frac{1}{2+1/\beta}}, & 2/(2\beta+1)<r\leq 2,
\\
%\phi_n\left(\bB^\beta_{r\theta}(L)\right)\geq
n^{-\frac{1}{1-1/(\beta r)+1/\beta}}, & r>2.
\end{array}
\right.
$$
 They also proposed
a wavelet--based estimator that achieves the indicated rates up to a
logarithmic factor for a scale of
Besov classes under additional assumption $2\beta r>2-r$. It follows from \cite{Donoho}
that if $p=2$ and the density is compactly supported
then the
corresponding rates are $\phi_n(\Sigma)\asymp n^{-1/(2+1/\beta)}$
for all $r\geq 2/(2\beta+1)$. Hence the rate corresponding to the zone $r>2$, $p=2$,
does not appear in the case
of the compactly supported densities.
\par
As for  the multivariate setting, Ibragimov and Khasminskii in a series of papers
[\cite{Ibr-Has1}, and \cite{Has-Ibr}]
studied the problem of minimax density estimation
 with $\bL_p$--loss on $\bR^d$. Together with some classes of infinitely
differentiable densities, they considered the anisotropic  Nikolskii's classes
$\Sigma=\bN_{\vec{r},d}\big(\vec{\beta},\vec{L}\big)$, where $\vec{\beta}=(\beta_1,\ldots,\beta_d)$,
$\vec{r}=(r_1,\ldots,r_d)$ and $\vec{L}=(L_1,\ldots,L_d)$
%characterized by the H\"older conditions in $\bL_p$ for the differences of the derivatives
%of various orders
(for the precise definition see
Section~\ref{sec:nikolski}).
It was shown  that if $r_i=p$ for all $i=1,\ldots,d$ then
\begin{equation}\label{eq:ibr-has}
\phi_n\big(\bN_{\vec{r},d}(\vec{\beta},\vec{L})\big)
\asymp
\left\{
\begin{array}{ll}
n^{-\frac{1-1/p}{1-1/(\beta p)+1/\beta}}, & p\in [1,2),\\
%\phi_n\left(\bN_{\vec{r},d}\big(\vec{\beta},\vec{L}\big)\right)=
n^{-\frac{1}{2+1/\beta}}, & p\in [2,\infty).
\end{array}
\right.
\end{equation}
Here $\beta$ is the parameter defined by the relation $1/\beta=\sum_{j=1}^d 1/\beta_j$.
It should be stressed that in the cited papers the same norm index $p$ is used
in the definitions of the risk and of the functional class. We also refer to
the recent paper by \cite{mason},
where further discussion of these results can be found.
\par
\cite{delyon-iod} generalized  the results of
\cite{Donoho} to the minimax density estimation on a
bounded interval of $\bR^d$, $d\geq 1$ over a collection of
the isotropic
Besov classes. In particular, they showed that the minimax rates of
convergence  given by (\ref{eq:rate-donoho}) hold
with $1/(\beta r)$ and $1/\beta$ replaced by $d/(\beta r)$ and
$d/\beta$ respectively.
Comparing rates in (\ref{eq:ibr-has}) with
the asymptotics of minimax risk found in \cite{delyon-iod} with $r=p$
we conclude that the rate in (\ref{eq:ibr-has})
in the zone $p\in [1,2)$ does not appear for
compactly supported densities.
\par
Recently \cite{GL11} developed an adaptive minimax estimator
over a scale of classes $\bN_{\vec{r},d}(\vec{\beta},\vec{L})$;
in particular, if $r_i=p$ for all $i=1,\ldots, d$ then
their estimator attains the minimax rates indicated in
(\ref{eq:ibr-has}).
%$$
%\phi_n\big(\bN_{\vec{r},d}
%(\vec{\beta},\vec{L})\big)
%\asymp
%\left\{
%\begin{array}{cc} n^{-\frac{1}{2+1/\beta}}, & p\geq 2,\\
%\phi_n\left(\bN_{\vec{r},d}\big(\vec{\beta},\vec{L}\big)\right)\asymp
%n^{-\frac{1-1/p}{1-1/(\beta p)+1/\beta}}, & p\in (1,2).
%\end{array}
%\right.
%$$
Note that in the considered setting the norm
indexes in the definitions of the risk and the functional class coincide.
\par
The results discussed above
show that there is an essential difference between
the problems of density estimation on the whole space
and on a compact interval.
The literature on density estimation on the whole space is quite fragmented, and
relationships between aforementioned results are yet to be understood.
These relationships
become even more complex and interesting in the multivariate setting
where the density to be estimated belongs
to a functional class with anisotropic and inhomogeneous smoothness.
%Another important issue is the simultaneous influence of
% the anisotropy and the inhomogeneity of the the density to be estimated
% (the coordinates of $\vec{r}$ are different and not related to the index $p$)
%on the asymptotics of minimax risk.
The problem of minimax estimation
under  $\bL_p$--loss over homogeneous Sobolev $\bL_q$--balls ($q\ne p$)
was initiated in \cite{nemirovski} in the regression model on the unit cube of $\bR^d$.
For the first time,
functional classes
with anisotropic and inhomogeneous smoothness were considered
in  \cite{lepski-kerk,lepski-kerk-08} for
the Gaussian white noise model on a compact subset of $\bR^d$.
In the density estimation model
\cite{akakpo} studied the case $p=2$
and considered compactly supported  densities on $[0,1]^d$.
%\par
%A related  problem of minimax nonparametric estimation of a multivariate regression function with $\bL_p$--losses
%on isotropic Sobolev classes was studied by \cite{nemirovski}.
%He proposed an optimal in order estimator, and showed existence of two zones
%with respect to the behavior of the minimax risk.
\par
To the best of our knowledge, the problem of estimating  a multivariate density
from anisotropic and inhomogeneous functional classes on $\bR^d$
was not considered in the literature. This problem is a subject of the current paper.
Our results cover the existing ones and
generalize them in the following directions.
\par\smallskip
1.~We fully characterize behavior of the minimax risk for all possible
relationships between regularity parameters and norm indexes in the definition of
the functional classes and of the risk.
In particular, we discover that there are four different regimes
with respect to the minimax rates of convergence:
{\em tail, dense and sparse zones}, and
the last zone, in its turn, is subdivided in two regions.
Existence of these regimes is not a consequence  of the multivariate nature
of the problem or the considered functional classes; in fact,
these regimes appear already in the dimension one.
Thus our results reveal all possible zones with respect to the
rates of convergence in the problem of
density estimation on $\bR^d$
and
explain different results on
rates of convergence in the existing literature.
In particular, results in
\cite{Juditsky} and \cite{patricia} pertain
to
the rates of convergence in the
tail and dense zones, while those
in \cite{Donoho} and \cite{delyon-iod} correspond to
the dense zone and to a subregion  of the sparse zone.
\par
2.~We propose an estimator
that  is based upon a data--driven selection from a family of kernel estimators, and
establish for it a point--wise oracle inequality.
Then we use this inequality for derivation of bounds on the
$\bL_p$--risk over a collection of the Nikol'skii functional classes.
Since the
construction of our estimator
does not use any prior information
on the class parameters, it is
adaptive minimax
over a scale of these classes. Moreover,
we believe that
the method of deriving
$\bL_p$--risk bounds from point--wise oracle inequalities
employed in the proof of Theorem~\ref{th:main}
is of interest in its own right. It is
quite general and can be applied to
other nonparametric estimation problems.
\par
3.~Another issue  studied  in the present paper
is related to the existence of the tail zone.
This zone does not exist in the problem of estimating compactly supported
densities. Then a natural question
arises: what is a general condition on $f$ which
ensures the same asymptotics of the minimax risk
on $\bR^d$ as in the case of
compactly supported densities?
We propose a {\em tail dominance condition} and show that, in a sense,
it is the weakest possible condition under which the tail zone disappears.
We also show that this condition guarantees
existence of a consistent estimator under $\bL_1$-loss.
Recall that smoothness  alone
is not sufficient in order to guarantee consistency
of density estimators in $\bL_1(\bR^d)$ [see \cite{Ibr-Has2}].
\par\smallskip
The paper is structured as follows.
In Section~\ref{sec:procedure+oracle}
we define our estimation procedure and derive the corresponding point--wise
oracle inequality. Section~\ref{sec:main-results}
presents  upper and lower bounds on the minimax risk.
We also discuss the obtained results and relate  them to the existing
results in the literature.
The same estimation problem under the tail dominance condition
is studied in
Section~\ref{sec:tail-dom}.
Sections~\ref{sec:Proof-th-3}--\ref{sec:Proof-th-2}
contain
proofs of Theorems~\ref{th:oracle-inequality}--\ref{th:upper-new}; proofs
of auxiliary results are relegated to Appendices~A and~B.
\par
The following notation and conventions are used throughout the paper.
For vectors $u, v\in \bR^d$ the operations
$u/v$, $u\vee v$, $u\wedge v$ and inequalities such as $u\leq v$
are all understood in the coordinate--wise sense. For instance,
$u\vee v =(u_1\vee v_1,\ldots, u_d\vee v_d)$.
All integrals are taken over $\bR^d$
unless the domain of integration is specified explicitly.
For a Borel set
$\cA\subset\bR^d$ symbol
$|\cA|$ stands for the Lebesgue measure of $\cA$;
if $\cA$ is a finite set, $|\cA|$ denotes the cardinality
of $\cA$.
\section{Estimation procedure and point--wise oracle inequality}
\label{sec:procedure+oracle}
In this section we define  our estimation procedure and
derive an upper bound  on  its point--wise risk.
%Now we turn to the density estimation problem formulated in the introduction section.
\subsection{Estimation procedure}
\label{sec:procedure}
Our estimation procedure is based on data-driven selection
from a family of kernel estimators. The family
of estimators is defined as follows.
\subsubsection{Family of kernel estimators}
Let $K:[-1/2,1/2]^d\to\bR^1$ be a fixed   kernel such that
\mbox{$K\in\bC(\bR^d)$},
$\int K(x)\rd x=1$, and
$\|K\|_\infty<\infty$.
Let
$$
 \cH =\Big\{h=(h_1,\ldots,h_d)\in (0,1]^d: \;h_j=2^{-k_j},
 k_j=0,\ldots, \log_2n,\;j=1,\ldots,d\Big\};
$$
without loss of generality we assume that $\log_2n$ is integer.
\par
Given  a {\em bandwidth} $h\in \cH$, define
the corresponding kernel estimator of $f$ by
the formula
\begin{equation}\label{eq:kernel-est}
\hat{f}_h (x):=
\frac{1}{nV_h} \sum_{i=1}^n K\Bigg(\frac{X_i-x}{h}\Bigg)=
\frac{1}{n}\sum_{i=1}^n K_h(X_i-x),
\end{equation}
where $V_h:=\prod_{j=1}^d h_j$, $K_h(\cdot):=(1/V_h)K(\cdot/h)$.
Consider the family of kernel estimators
\begin{equation*}
%\label{eq:kernel-family}
\cF(\cH):=\{\hat{f}_h, h\in \cH\}.
\end{equation*}
The proposed estimation procedure is based on data--driven selection of an estimator
from~$\cF(\cH)$.
\subsubsection{Auxiliary estimators}
Our selection rule uses auxiliary estimators that are constructed as follows.
For any pair $h,\eta \in \cH$ define the kernel $K_h*K_\eta$ by
the formula
$
[K_h*K_\eta](t)=\int K_h(t-y) K_\eta(y)\rd y.
$
Let
 $\hat{f}_{h,\eta}(x)$ denote the estimator  associated with this kernel:
\begin{equation*}%\label{eq:f-eta-h}
 \hat{f}_{h,\eta}(x)=\frac{1}{n}\sum_{i=1}^n K_{h,\eta}(X_i-x),\;\;\;
 K_{h,\eta}=K_h*K_\eta.
\end{equation*}
The following representation of kernels $K_{h,\eta}$
will be useful: for any $h,\eta\in \cH$
\begin{equation}\label{eq:repres-1}
 [K_h*K_\eta](t) \;=\;
\frac{1}{V_{h\vee \eta}} Q_{h,\eta}\Big(\frac{t}{h\vee \eta}\Big),
\end{equation}
where  function $Q_{h,\eta}$ is given by the formula
 \begin{equation}\label{eq:Q}
 Q_{h,\eta}(t)=\int K\big(v(y,t-\nu y)\big) K\big(v(t-\nu y,y)\big) \rd y,\;\;\; \nu:=\frac{h\wedge \eta
 }{h\vee \eta}.
\end{equation}
Here function $v:\bR^d\times \bR^d \to \bR^d$ is defined by
\[
v_j(y,z)=\left\{
\begin{array}{cc}
y_j, & h_j\leq \eta_j,\\
z_j , & h_j>\eta_j,
\end{array}
\right.,\;
\;\;\;\;j=1,\ldots,d.
\]
\par
The representation (\ref{eq:repres-1})--(\ref{eq:Q})
is obtained by a straightforward change of variables in the convolution integral
[see the proof of Lemma~12 in \cite{GL11a}].
We also note that ${\rm supp}(Q_{h,\eta})\subseteq [-1,1]^d$,
and $\|Q_{h,\eta}\|_\infty \leq \|K\|_\infty^2$ for all $h,\eta$.
In the special case where  $K(t)=\prod_{i=1}^d k(t_i)$ for some univariate kernel
$k:[-1/2, 1/2]\to \bR^1$  we have
\[
 Q_{h,\eta} (t) =\prod_{i=1}^d \int k(t_i-\nu_i u_i)k(u_i)\rd u_i,\;\;\;
\nu_i=(h_i\wedge \eta_i)/(h_i\vee \eta).
\]
We also define
\begin{equation*}%\label{eq:QQ}
 Q(t) =
\sup_{h,\eta\in \cH}\Big|\int K\big(v(y,t-\nu y)\big) K\big(v(t-\nu y,y)\big) \rd y\Big|,
\end{equation*}
and
note that ${\rm supp}(Q)\subseteq [-1,1]^d$, and $\|Q\|_\infty\leq \|K\|_\infty^2$.
\subsubsection{Stochastic errors of kernel estimators and their majorants}
Uniform moment bounds on stochastic errors of kernel estimators $\hat{f}_h(x)$ and
$\hat{f}_{h,\eta}(x)$ will play an important role in the
construction of our selection rule.
Let
\begin{eqnarray}
 \xi_{h}(x) &=& \frac{1}{n}\sum_{i=1}^n K_h(X_i-x) - \int K_h(t-x)f(t) \rd t,
\label{eq:xi-h}
\\
 \xi_{h, \eta}(x) &=&  \frac{1}{n}\sum_{i=1}^n K_{h,\eta}(X_i-x) - \int K_{h,\eta}(t-x)f(t) \rd t
%\label{eq:xi-h-eta}
\nonumber
\end{eqnarray}
denote the stochastic errors of $\hat{f}_h$ and $\hat{f}_{h,\eta}$ respectively.
In order to construct our selection rule we need
to find uniform upper bounds ({\em majorants}) on $\xi_h$ and $\xi_{h,\eta}$, i.e.
we need  to find
functions
$M_h$ and $M_{h,\eta}$
such that moments of random variables
\begin{equation}\label{eq:unif-maj}
 \sup_{h\in \cH}\big[|\xi_h(x)|- M_h(x)\big]_+,\;\;\;
\sup_{h,\eta\in \cH}\big[
|\xi_{h,\eta}(x)| - M_{h,\eta}(x)\big]_+
\end{equation}
are ``small'' for each $x\in \bR^d$. We will be also interested in the integrability
properties of these moments.
\par
It turns out that
the majorants $M_h(x)$ and $M_{h,\eta}(x)$ can be
defined in the following way.
For a function \mbox{$g:\bR^d\to \bR^1$} let
\begin{equation}
 \label{eq:A-A}
A_h(g,x)= \int |g_h(t-x)| f(t)\rd t,\;\;\;g_h(\cdot)=V^{-1}_h g\big(\cdot/h\big),\;\;\;\;h\in \cH.
\end{equation}
Now define
\begin{equation}\label{eq:A-M}
M_h(g,x)=\sqrt{\frac{\kappa A_h(g,x)\ln n}{nV_h}} + \frac{\kappa\ln n}{nV_h},
\end{equation}
where $\kappa$ is a positive constant to be specified.
In Lemma~\ref{lem:1} in Section~\ref{sec:Proof-th-3} we show that under appropriate choice of parameter
$\kappa$ functions
\begin{equation}\label{eq:M-h-M-h-eta}
M_h(x):= M_h(K,x),\;\;\;\;M_{h,\eta}(x):= M_{h\vee \eta}(Q, x)
\end{equation}
uniformly majorate  $\xi_h$ and $\xi_{h,\eta}$ in the sense that the moments
of random variables in (\ref{eq:unif-maj}) are ``small''.
\par
It should be noted, however, that functions $M_h(x)$ and $M_{h,\eta}(x)$
given by (\ref{eq:M-h-M-h-eta})
cannot be directly used
in construction of the selection rule because they depend on unknown density $f$
to be estimated. We will use empirical counterparts of $M_h(x)$ and $M_{h,\eta}(x)$ instead.
\par
For $g:\bR^d\to\bR^1$ we let
$$
\hat{A}_h(g, x)=\frac{1}{n}\sum_{i=1}^n |g_h(X_i-x)|,
$$
and define
\begin{equation}\label{eq:hat-A-M}
\hat{M}_h(g,x)= 4\sqrt{\frac{\kappa\hat{A}_h(g,x)\ln n}{nV_h}} + \frac{4\kappa\ln n}{nV_h}.
\end{equation}

\subsubsection{Selection rule and final estimator}
Now we are in a position  to define our selection rule.
For every $x\in \bR^d$ let
\begin{eqnarray}
\hat{R}_h(x)\;=\;\sup_{\eta\in \cH}\Big[|\hat{f}_{h,\eta}(x)-\hat{f}_\eta(x)| -
\hat{M}_{h\vee \eta}(Q, x)-\hat{M}_\eta(K, x)\Big]_+  &&
\nonumber
\\
+\;\;
\sup_{\eta\geq h} \hat{M}_{\eta}(Q, x) + \hat{M}_h(K,x), && h\in \cH.
\label{eq:hat-R}
\end{eqnarray}
The
selected bandwidth $\hat{h}(x)$ and the corresponding estimator are defined by
\begin{equation}\label{eq:hat-h}
 \hat{h}(x)={\rm arg}\inf_{h\in \cH} \hat{R}_h(x),\;\;\;\hat{f}(x)=\hat{f}_{\hat{h}(x)}(x),\;\;\;
x\in \bR^d.
\end{equation}
\par
Note that the estimation procedure is completely determined by the
family of kernel estimators $\cF(\cH)$ and by the constant
$\kappa$ appearing in the definition of~$\hat{M}_h$.
\par
We have to ensure that the map
$x\mapsto \hat{f}_{\hat{h}(x)}(x)$ is an $X^{(n)}$-measurable Borel function.
This follows from continuity of $K$ and the fact that $\cH$ is a discrete set;
for details see Appendix~A, Section~\ref{meas}.
\par
The main idea  behind the construction of the selection
procedure (\ref{eq:hat-R})--(\ref{eq:hat-h}) is the following.
The expression $\hat{M}_{h\vee \eta}(Q, x)+\hat{M}_\eta(K, x)$
appearing in the square brackets in (\ref{eq:hat-R}) dominates with
large probability the stochastic part of the difference
$|\hat{f}_{h,\eta}(x)-\hat{f}_\eta(x)|$. Consequently, the first
term on the right hand side of (\ref{eq:hat-R}) serves as a proxy
for the deterministic part of $|\hat{f}_{h,\eta}(x)-\hat{f}_\eta(x)|$
which is the absolute value of the difference of biases
of kernel estimates $\hat{f}_{h,\eta}(x)$
and $\hat{f}_\eta(x)$. The latter, in its own turn, is closely related to
the bias of the estimator $\hat{f}_h(x)$. Thus, the first term on
the right hand side of (\ref{eq:hat-R}) is a proxy for the bias
of $\hat{f}_h(x)$, while the second term is an upper bound on the
standard deviation of $\hat{f}_h(x)$.
\subsection{Point--wise oracle inequality}
\label{sec:Pointwise oracle inequality}
%The upper bounds on the risk of our estimator $\hat{f}(x)$ are derived
%from the point--wise oracle inequality given in the next theorem. This result is of interest on its own right. \par
Let $B_h(f,t)$ be the bias of the kernel estimator $\hat{f}_h(t)$,
\begin{equation}\label{eq:B-h}
B_h(f,t)=\int K_\eta(y-t) f(y)\rd y -f(t),
\end{equation}
and define
\begin{eqnarray}\label{eq:maximal-bias}
 \bar{B}_h(f,x)= |B_h(f,x)|\,\vee\,\sup_{\eta\in \cH} \Big|\int K_\eta(t-x) B_h(f,t)\rd t\Big|.
\end{eqnarray}
\begin{theorem}\label{th:oracle-inequality}
For any $x\in \bR^d$ one has
\begin{eqnarray}
|\hat{f}(x)-f(x)| \leq \inf_{h\in \cH}
\big\{4\bar{B}_h(f,x)+
60 \sup_{\eta\geq h} M_\eta (Q, x)+ 61M_h(K,x)\big\}
%\nonumber
%\\
+7\zeta(x)\; +\; 18\chi(x),
\label{eq:pointwise-oracle}
\end{eqnarray}
where
\begin{eqnarray}
\zeta(x)&:=&\sup_{h\in \cH} [|\xi_h(x)|- M_h(K,x)]_+ \;\vee\;
\sup_{h,\eta\in \cH}[|\xi_{h,\eta}(x)| - M_{h\vee \eta}(Q,x)]_+,
\label{eq:zeta}
\\
\chi(x)&:=& \max_{g\in \{K,Q\}}\sup_{h\in \cH}\big[|\hat{A}_h(g,x)-A_h(g,x)|- M_h(g,x)\big]_+.
\label{eq:chi}
\end{eqnarray}
Furthermore, for any $q\geq 1$
if
$\kappa \geq   [\|K\|_\infty\vee 1]^2 [(4d+2)q+4(d+1)]$ then
\begin{equation}\label{eq:integrability-of-zeta}
\int \bE_f\big\{[\zeta (x)]^q +[\chi(x)]^q\big\} \rd x
\;\leq\; C n^{-q/2}, \;\;\forall n\geq 3,
\end{equation}
where $C$  is the constant depending on $d$, $q$ and $\|K\|_\infty$ only.
\end{theorem}
We remark that
Theorem~\ref{th:oracle-inequality} does not require any conditions on the
estimated density $f$.
%
%is established over the set of  all probability densities.
%
\section{Adaptive estimation over anisotropic Nikol'skii classes}
\label{sec:main-results}
In this section
we study
%adaptive
properties of the estimator defined in
(\ref{eq:hat-R})--(\ref{eq:hat-h}).
%In this context
The point--wise oracle inequality
%presented in
of Theorem~\ref{th:oracle-inequality} is the key technical tool
for %computing
bounding $\bL_p$-risk of this estimator on the anisotropic Nikol'skii classes.
\subsection{Anisotropic Nikol'skii classes}\label{sec:nikolski}
Let $(e_1,\ldots,e_d)$ denote the canonical basis of $\bR^d$.
For function $g:\bR^d\to \bR^1$ and
real number $u\in \bR$ define
{\em the first order difference operator with step size $u$ in direction of the variable
$x_j$}~by
\[
 \Delta_{u,j}g (x)=g(x+ue_j)-g(x),\;\;\;j=1,\ldots,d.
\]
By induction,
the $k$-th order difference operator with step size $u$ in direction of the variable $x_j$ is
defined~as
\begin{equation}\label{eq:Delta}
 \Delta_{u,j}^kg(x)= \Delta_{u,j} \Delta_{u,j}^{k-1} g(x) = \sum_{l=1}^k (-1)^{l+k}\binom{k}{l}\Delta_{ul,j}g(x).
%g(x+lue_j),\;\;j=1,\ldots,d.
\end{equation}
\begin{definition}
For given  real numbers $\vec{r}=(r_1,\ldots,r_d)$, $r_j\in [1,\infty]$, $\vec{\beta}=(\beta_1,\ldots,\beta_d)$,
$\beta_j>0$, and $\vec{L}=(L_1,\ldots, L_d)$, $L_j>0$, $j=1,\ldots, d$, we
say that function $g:\bR^d\to \bR^1$ belongs to the anisotropic
Nikol'skii class $\bN_{\vec{r},d}\big(\vec{\beta},\vec{L}\big)$ if
\begin{itemize}
\item[{\rm (i)}] $\|g\|_{r_j}\leq L_{j}$ for all $j=1,\ldots,d$;
\item[{\rm (ii)}]
for every $j=1,\ldots,d$ there exists natural number  $k_j>\beta_j$ such that
\begin{equation}\label{eq:Nikolski}
 \Big\|\Delta_{u,j}^{k_j} g\Big\|_{r_j} \leq L_j |u|^{\beta_j},\;\;\;\;
\forall u\in \bR^d,\;\;\;\forall j=1,\ldots, d.
\end{equation}
\end{itemize}
\end{definition}
%\begin{remark}\mbox{}
%\begin{itemize}
%\item
%The condition that for every $j=1,\ldots,d$ there exists $k_j>\beta_j$
%such that (\ref{eq:Nikolski}) holds can be replaced by the condition that
%(\ref{eq:Nikolski}) holds for every $k>\beta_j$, $j=1,\ldots,d$;
%see, \cite[Section~4.3.3]{Nikolski}.
\par
The anisotropic Nikol'skii class is a  specific case of the anisotropic Besov class,
often encountered in the nonparametric estimation literature.
%Following  \cite[Section~4.3.4]{Nikolski}
In particular,
$\bN_{\vec{r},d}\big(\vec{\beta},\cdot\big)=
\bB_{r_1,\ldots,r_d;\infty,\ldots,\infty}^{\beta_1,\ldots,\beta_d}(\cdot)$,
see \cite[Section~4.3.4]{Nikolski}.
%\item

%\end{itemize}
%\end{remark}
%\par
%
\subsection{Construction of kernel $K$}
We will use the following specific kernel $K$ in the definition of the family $\cF(\cH)$
[see, e.g., \cite{lepski-kerk} or \cite{GL11}].
\par
 Let  $\ell$ be an integer number,
and let $w:[-1/(2\ell), 1/(2\ell)]\to \bR^1$ be a function satisfying  $\int w(y)\rd y=1$,
and $w\in\bC(\bR^1)$. Put
\begin{equation}\label{eq:w-function}
 w_\ell(y)=\sum_{i=1}^\ell \binom{\ell}{i} (-1)^{i+1}\frac{1}{i}w\Big(\frac{y}{i}\Big),\qquad
 K(t)=\prod_{j=1}^d w_\ell(t_j),\;\;\;\;t=(t_1,\ldots,t_d).
\end{equation}
%and define
%\begin{equation}\label{eq:kernel-K}
% K(t)=\prod_{j=1}^d w_\ell(t_j),\;\;\;\;t=(t_1,\ldots,t_d).
%\end{equation}
The kernel $K$ constructed in this way is bounded, supported on $[-1/2,1/2]^d$, belongs to $\bC(\bR^d)$
%$K\in\bC(\bR^d)$
and satisfies
\begin{equation*}%\label{eq:kernel}
 \int K(t)\rd t=1,\;\;\;\int K(t) t^k \rd t=0,\;\;\forall |k|=1,\ldots, \ell-1,
\end{equation*}
where $k=(k_1,\ldots,k_d)$ is the multi--index, $k_i\geq 0$, $|k|=k_1+\cdots+k_d$, and
$t^k=t_1^{k_1}\cdots t_d^{k_d}$ for $t=(t_1,\ldots, t_d)$.
\par
%Now consider the family of kernel estimatos $\cF(\cH)$ defined in
%(\ref{eq:kernel-est})--(\ref{eq:kernel-family}) and associated with kernel
% (\ref{eq:kernel-K})--(\ref{eq:kernel}).

\subsection{Main results}
\label{sec:adap-nikol}
Let $\bN_{\vec{r},d}\big(\vec{\beta},\vec{L}\big)$ be the
anisotropic Nikol'skii functional class.
%
%$\bF(M)=\left\{f:\|f\|_\infty\leq M\right\}$.
Put
\[
 \frac{1}{\beta} := \sum_{j=1}^d \frac{1}{\beta_j},\;\;\;\;\;
\frac{1}{s} := \sum_{j=1}^d \frac{1}{\beta_jr_j}, \;\;\;\;\;L_\beta:= \prod_{j=1}^d L_j^{1/\beta_j},
\]
and define
\begin{eqnarray}
 \label{eq:nu}
 \nu &=&\left\{\begin{array}{lccl}
           \frac{1-1/p}{1-1/s+1/\beta }, & p <\frac{2+1/\beta}{1+1/s},\\*[2mm]
\;\;\frac{\beta}{2\beta+1}, &    \frac{2+1/\beta}{1+1/s} \leq p \leq s(2+1/\beta),\\*[2mm]
\quad\;s/p,  & p> s(2+1/\beta),\;s<1,
\\*[2mm]
\frac{1-1/s + 1/(p\beta)}{2-2/s+1/\beta}, & p> s(2+1/\beta),\; s\geq 1,
            \end{array}
\right.
\\*[4mm]
\mu_n &=& \left\{
\begin{array}{lll}
(\ln n)^{d/p}, &  p \leq\frac{2+1/\beta}{1+1/s};
\\*[2mm]
(\ln n)^{1/p}, & p=s(2+1/\beta),
\\*[2mm]
  \quad               1, & {\rm otherwise}.
                 \end{array}
\right.
%\label{eq:mu}
\nonumber
\end{eqnarray}
\par
In contrast to Theorem \ref{th:oracle-inequality} proved over the set of all probability densities,
the adaptive results presented below
require the additional assumption:
the estimated density should be uniformly bounded. For this purpose we define   for
$M>0$
\[
 \bN_{\vec{r},d}\big(\vec{\beta},\vec{L}, M\big) \;:=\; \bN_{\vec{r},d}\big(\vec{\beta},\vec{L}\big)
\;\cap \; \left\{f:\|f\|_\infty\leq M\right\}.
\]
Note, however, that  if $J:=\{j=1,\ldots,d:\;r_j=\infty\}$ then
$
\bN_{\vec{r},d}\big(\vec{\beta},\vec{L}, M\big)=\bN_{\vec{r},d}\big(\vec{\beta},\vec{L})
$
with $M=\inf_{J}L_j$.
Moreover, in view of the embedding theorem for the anisotropic Nikol'skii classes
[see Section~\ref{subsec:preliminaries} below], condition $s>1$
implies that the density to be estimated belongs to a class of
uniformly bounded and continuous functions. Thus, if $s>1$ one has
$
\bN_{\vec{r},d}\big(\vec{\beta},\vec{L}, M\big)=\bN_{\vec{r},d}\big(\vec{\beta},\vec{L})
$
with some $M$ completely determined by $\vec{L}$.

The asymptotic behavior of  the $\bL_p$-risk on
class $\bN_{\vec{r},d}(\vec{\beta},\vec{L}, M)$
is characterized in the next two theorems.
\par
Let family $\cF(\cH)$ be associated with kernel (\ref{eq:w-function}).
Let $\hat{f}$ denote the estimator given by the
selection rule (\ref{eq:hat-R})--(\ref{eq:hat-h}) with
$\kappa = (\|K\|_\infty \vee 1)^2[(4d+2)p+4(d+1)]$
that is
applied to the family $\cF(\cH)$.
\begin{theorem}\label{th:main}
For any  $M>0$, $L_0>0$,  $\ell\in\bN^*$, any
$\vec{\beta}\in (0,\ell]^d$, $\vec{r}\in (1,\infty]^d$, any $\vec{L}$ satisfying
$\min_{j=1,\ldots,d}L_j\geq L_0$, and any $p\in (1,\infty)$ one has
\begin{equation*}
%\label{eq:result}
 \limsup_{n\to\infty}\bigg\{
\mu_n\Big( \frac{L_\beta\ln n}{n}\Big)^{-\nu}\;
\cR_p^{(n)}\big[\hat{f}\,;\bN_{\vec{r},d}\big(\vec{\beta},\vec{L}, M\big)\big]\bigg\}\leq C <\infty.
\end{equation*}
Here constant $C$ does not depend on $\vec{L}$ in the cases $p\leq s(2+1/\beta)$ and $p\geq s(2+1/\beta)$,
$s<1$.
\end{theorem}
\begin{remark}
\mbox{}
{\rm
\begin{enumerate}
\item[\rm 1.]
Condition $\min_{j=1,\ldots,d}L_j\geq L_0$ ensures independence of the constant $C$ on $\vec{L}$
in the cases $p\leq s(2+1/\beta)$ and $p\geq s(2+1/\beta)$,
$s<1$.
If  $p\geq s(2+1/\beta)$, $s\geq 1$ then $C$ depends  on $\vec{L}$,
and the corresponding expressions can be easily extracted  from the proof of the theorem.
We note that in this case the map $\vec{L}\mapsto C(\vec{L})$
is bounded on every closed cube of $(0,\infty)^{d}$.
\item[\rm 2.]
  We consider the case $1<p<\infty$ only, not including $p=1$ and $p=\infty$.
It is well--known, \cite{Ibr-Has2},
 that smoothness  alone
is not sufficient in order to guarantee consistency of density estimators in $\bL_1(\bR^d)$; see also
Theorem \ref{th:lower-bound-in-L_p} for a lower bound.
The case $p=\infty$ was considered recently in
\cite{lep2012}.
%Contrary to Theorem~\ref{th:oracle-inequality} which holds forallany probability density
%proved over the
%set of all probability densities  the assertions of
% Theorem~\ref{th:main}
%requires
%additionally
%uniform boundedness of the estimated density, i.e. $\|f\|_\infty\leq M$.
% $f\in\bF(M)$.
%Note, however, that our selection rule led to the estimator
%$\hat{f}$ is independent of the knowledge of $M$.
%Hence, $\hat{f}$ is fully adaptive, i.e. is independent of all parameters
%involved in the description of $\bL_p$-risk.
\item[\rm 3.]
As it was discussed above,
Theorem~\ref{th:main} requires uniform boundedness of the estimated density, i.e.
$\|f\|_\infty\leq M<\infty$. We note however that
our estimator $\hat{f}$ is fully adaptive, i.e., its construction does not
use any information
on  the parameters  $\vec{\beta}, \vec{r}, \vec{L}$ and $M$.
\end{enumerate}
}
\end{remark}

\iffalse

\smallskip

\begin{remark}
Contrary to the Theorem \ref{th:oracle-inequality} proved over the set
of all probability densities  the assertions of Theorem \ref{th:main}
require additionally  the uniform boundedness of the estimated density, i.e. $f\in\bF(M)$. Note, however, that our selection rule led to the estimator
$\hat{f}$ is independent of the knowledge of $M$. Hence, $\hat{f}$ is fully adaptive,
i.e. is independent of all parameters involved in the description of $\bL_p$-risk.
\end{remark}

%\smallskip
\fi

\par
Now we present lower bounds on the minimax risk.
Define
$$
\alpha_n = \left\{
\begin{array}{ll}
\ln n, &  p> s(2+1/\beta),\; s\geq 1,
\\*[2mm]
  \quad               1, & {\rm otherwise}.
                 \end{array}
\right.
$$
%Let us  find now the asymptotic of minimax risk. Define
%\[
%\psi_n=
%\max
%\Bigg[\bigg(\frac{L_\beta}{n}\bigg)^{\frac{1-1/p}{1-1/s+1/\beta}},\;
%\bigg(\frac{L_\beta}{n}\bigg)^{\frac{1}{2+1/\beta}},\;\mathbf{1}_{\{s<1\}}\bigg(\frac{L_\beta}{n}\bigg)^{\frac{s}{p}}+
%\mathbf{1}_{\{s\geq1\}}L_\beta^{\frac{1/2 - 1/p}{1- 1/s + 1/(2\beta)}}
%\bigg(\frac{\ln n}{n}\bigg)^{\frac{1-1/s + 1/(p\beta)}{2- 2/s +
%1/\beta}}\Bigg].
%\]
%
\begin{theorem}
\label{th:lower-bound-in-L_p}
%The following statements hold.
%\begin{itemize}
%\item[{\rm (i)}]
Let $\vec{\beta}\in (0,\infty)^{d}$, $\vec{r}\in [1,\infty]^{d}$,  $\vec{L}\in (0,\infty)^{d}$ and $M>0$ be fixed.

{\rm (i)}\;\;
 There exists  $c>0$
 such that
$$
\liminf_{n\to\infty}\bigg\{
\Big( \frac{L_\beta\alpha_n}{n}\Big)^{-\nu} \inf_{\widetilde{f}}
\cR_p^{(n)}\big[\widetilde{f};\;\bN_{\vec{r},d}(\vec{\beta},\vec{L}, M)\big]
\bigg\} \geq c,\quad \forall p\in [1,\infty),
$$
where
the infimum is taken over all possible estimators~$\widetilde{f}$.
If $\min_{j=1,\ldots,d} L_j\geq L_0>0$ then
in the cases $p\leq s(2+1/\beta)$ or $p\geq s(2+1/\beta)$ and $s<1$
the  constant $c$ is
independent of  $\vec{L}$.
%\item[{\rm (ii)}]
\smallskip

 {\rm (ii)}\;\;Let $p=\infty$ and $s\leq 1$; then there is no
 consistent estimator, i.e., for some $c>0$
 %there exists
%a constant  such that
 $$
 \liminf_{n\to\infty}\inf_{\tilde{f}}
\sup_{f\in\bN_{\vec{r},d}(\vec{\beta},\vec{L}, M)}\bE_f\big\|\tilde{f}-f\big\|_\infty
 \;>\;c.
 $$
%\end{itemize}
\end{theorem}
\begin{remark}\mbox{}
{\rm
\begin{enumerate}
\item[{\rm 1}.] Inspection of the proof shows that
 if $\max_{j=1,\ldots,d} L_j\leq L_\infty<\infty$ then the statement~(i) is valid
with constant $c$
 depending on $\vec{\beta}, \vec{r}$,  $L_0$,  $L_\infty$, $d$ and $M$ only.
\item[{\rm 2}.] As it was mentioned above, adaptive minimax density estimation on $\bR^d$
under $\bL_\infty$--loss
was a subject of the recent paper
\cite{lep2012}.
A minimax adaptive estimator is constructed  in this paper under assumption $s>1$.
Thus, statement~(ii)
of Theorem~\ref{th:lower-bound-in-L_p} finalizes the research
on adaptive density estimation in the supremum norm.
It is interesting to note that the minimax rates
in the case $p=\infty$ coincide with  those
of Theorem~\ref{th:main} if we put formally
$p=\infty$.
\end{enumerate}
}
\end{remark}
\subsection{Discussion}
\label{sec:open-prob}
The results of Theorem~\ref{th:main} together with the matching
lower bounds of Theorem~\ref{th:lower-bound-in-L_p} provide complete classification
of minimax rates of convergence in the problem of density estimation on $\bR^d$.
In particular,
%these results show that there are
we discover four different zones with respect to the minimax rates of convergence.
\begin{itemize}
 \item {\em Tail zone} corresponds to ``small'' $p$, $1<p\leq\frac{2+1/\beta}{1+1/s}$.
This zone does not appear  if
density $f$ is assumed to be compactly supported, or
some tail dominance condition is imposed,
see Section~\ref{sec:tail-dom}.
\item {\em Dense zone} is characterized by  the ``intermediate'' range of $p$,
$\frac{2+1/\beta}{1+1/s}\leq p\leq s(2+1/\beta)$. Here the ``usual'' rate of convergence
$n^{-\beta/(2\beta+1)}$ holds.
\item {\em Sparse zone} corresponds to ``large'' $p$, $p\geq s(2+1/\beta)$.
As Theorems~\ref{th:main} and~\ref{th:lower-bound-in-L_p}  show,
this zone, in its turn, is subdivided into two regions with $s\geq 1$ and $s<1$. This
phenomenon   was  not observed in the existing literature
even for  settings with
compactly supported densities.
For other statistical models (regression, white Gaussian noise etc)
this result is also new.
\end{itemize}
It is important to emphasize that existence of these zones is not related to
the multivariate nature of the problem or to the anistropic smoothness of the estimated
density. In fact, these results hold already for the one--dimensional case, and this, to a
limited degree,
 was observed in the previous works.
In the subsequent remarks
we discuss relationships between
our results and the existing results in the literature, and comment on some open problems.
\par\smallskip
1. In \cite{Donoho}, \cite{delyon-iod} and \cite{lepski-kerk-08}
the sparse zone is defined as $p>2(1+1/\beta)$, $s>1$.
%In view of the embedding theorem for the anistropic Nikol'skii classes
%[see Section~\ref{subsec:preliminaries} below],
Recall that condition $s>1$
implies that the density to be estimated belongs to a class of
uniformly bounded and continuous functions. In the sparse zone
we consider also the case $s\leq 1$, but density $f$ is assumed to be
uniformly bounded. It turns out that in this zone
the rate corresponding to the index $\nu=s/p$ emerges.
%was not observed in the existing literature.
\par\smallskip
\par\smallskip
2.~The one--dimensional setting
was considered in
\cite{Juditsky} and \cite{patricia}. The setting of \cite{Juditsky}
corresponds to $s=\infty$, while \cite{patricia} deal with the
case of $p=2$
and $\beta>1/r-1/2$. Both settings rule out the sparse zone.
The rates of convergence in
the  dense zone
obtained in the aforementioned papers  are easily recovered from our results. However, in the tail zone our bound  contains additional $\ln(n)$-factor.
\par\smallskip
3.
In the previous papers on adaptive estimation of densities
with unbounded support
[cf. \cite{Juditsky} and \cite{patricia}] the
developed estimators are explicitly shrunk to zero.
This shrinkage is used in  bounding the minimax risk on the whole space.
We do not employ shrinkage in our estimator construction.
We derive bounds on the $\bL_p$--risk
by integration of the
point--wise oracle inequality  (\ref{eq:pointwise-oracle}).
The key elements of this derivation are
inequality (\ref{eq:integrability-of-zeta})
and  statement~(i) of Proposition~\ref{prop:new}.
The inequality (\ref{eq:integrability-of-zeta})
is based on the following fact:
the errors $\zeta(x)$ and $\chi(x)$ are integrable by the accurate choice
of the majorant. Indeed, Section~\ref{subsec:integrability}
shows that these
errors are not equal to zero
with probability which is integrable and ``negligible''
in the regions where the density is ``small''. This leads to integrability
of the remainders in (\ref{eq:pointwise-oracle}).
As for Proposition~\ref{prop:new}, it is related
to  the integrability of the main term in
(\ref{eq:pointwise-oracle}). The main problem here is that the majorant
$M_h(\cdot,x)$ itself  is not integrable. To overcome this difficulty
we use the integrability of the estimator $\hat{f}$,  approximation properties of
the density $f$, and (\ref{eq:integrability-of-zeta}).
\par\smallskip
4.~In the context of the Gaussian white noise model on a compact interval
\cite{lepski-kerk} developed
an adaptive estimator that
achieves the  rate of convergence
$(\ln{n}/n)^{\beta/(2\beta+1)}$ on the anisotropic Nikol'skii classes
under condition   $\sum_{i=1}^d [\frac{1}{\beta_i}(\frac{p}{r_i}-1)]_+<2$.
This restriction determines a part of the dense zone, and our
Theorem~\ref{th:main} improves
on this result. In fact, our estimator achieves the
rate $(\ln{n}/n)^{\beta/(2\beta+1)}$ in the zone
$\sum_{i=1}^d \frac{1}{\beta_i}(\frac{p}{r_i}-1)\leq 2$
which is equivalent to $p\leq s(2+1/\beta)$.
\par\smallskip
5. It follows from   Theorem~\ref{th:lower-bound-in-L_p}
that the upper bound of
Theorem~\ref{th:main} is sharp in the zone  $p> s(2+1/\beta)$, $s>1$,
and it is  nearly sharp up to a
logarithmic factor in all other zones.
This extra logarithmic factor is a consequence of the fact that we use
the point--wise selection procedure (\ref{eq:hat-R})--(\ref{eq:hat-h}). We also have
extra $\ln n$--term
on the boundaries $p=\frac{2+1/\beta}{1+1/s}$, $p=s(2+1/\beta)$.
\begin{conjecture} The rates found in Theorem \ref{th:lower-bound-in-L_p} are optimal.
\end{conjecture}
Thus, if our conjecture is true,
the construction of an estimator achieving the rates of
Theorem~\ref{th:lower-bound-in-L_p} in the tail and dense zones
remains an open problem.
\par\smallskip
6. Theorem~\ref{th:main} is proved under
assumption $\vec{r}\in (1,\infty]^{d}$, i.e., we do not include
the case where $r_j=1$ for some $j=1,\ldots,d$. This is related to
the construction of our selection rule, and to the necessity to bound
$\bL_{r_j}$--norm, $j=1,\ldots, d$ of the term
$\bar{B}_h(f,x)$; see (\ref{eq:maximal-bias}) and (\ref{eq:pointwise-oracle}).
In our derivations for this purpose we use properties of the strong
maximal operator [for details see  Section~\ref{subsec:preliminaries}], and
it is well--known that this operator is not
of the weak $(1,1)$--type in dimensions $d\geq 2$.
Nevertheless, using inequality (\ref{eq:weak-max})
we were able to obtain the following result.
\begin{corollary}
\label{cor:r=1}
Let $\vec{r}$ be such that
$r_j=1$ for some $j=1,\ldots,d$.
Then the result of Theorem \ref{th:main}  remains valid if the normalizing factor
$\big(n^{-1}\ln{n}\big)^{\nu}$ is replaced by $\big(n^{-1}[\ln{n}]^{d}\big)^{\nu}$.
\end{corollary}
\noindent
The proof of  Corollary~\ref{cor:r=1} coincides with the proof of
Theorem~\ref{th:main}
with the only difference that bounds in the proof of Proposition~\ref{prop:first-bound}
should use (\ref{eq:weak-max}) instead of the Chebyshev inequality.
This will result in
an extra $(\ln{n})^{d-1}$-factor. We note that
the results of Theorem~\ref{th:main}  and Corollary~\ref{cor:r=1}
coincide if $d=1$. It is not surprising
because in the dimension $d=1$
the strong maximal operator is the Hardy-Littlewood maximal function  which is
of  the weak (1,1)--type.
\section{Tail dominance condition}
\label{sec:tail-dom}
Let $g:\bR^d\to\bR^1$ be a locally integrable function.
Define  the map $g\mapsto g^*$  by the formula
\begin{equation}\label{eq:maximal-function}
 g^*(x):= \sup_{h\in (0,2]^d} \frac{1}{V_h} \int_{\Pi_h(x)}  g(t) \rd t,\;\;\;x\in \bR^d,
\end{equation}
where $\Pi_h(x)=[x_1-h_1/2,x_1+h_1/2]\times\cdots\times[x_d-h_d/2,x_d+h_d/2]$.
In fact, formula (\ref{eq:maximal-function}) defines the maximal operator
associated with the differential basis $\cup_{x\in \bR^d}\{ \Pi_h(x), h\in (0,2]\}$,
see \cite{Guzman}.
\par
Consider the following set of functions:
for any  $\theta\in (0,1]$ and $R\in (0,\infty)$ let
\begin{equation}\label{eq:G-theta}
\bG_\theta(R)=\big\{g:\bR^d\to\bR:\;\; \|g^*\|_\theta\leq R \big\}.
\end{equation}
Note that, although we keep the previous notation $\|g\|_\theta=(\int |g(x)|^\theta \rd x)^{1/\theta}$,
$\|\cdot\|_\theta$ is not longer a norm if $\theta\in (0,1)$.
\par
The assumption that
$f\in \bG_\theta(R)$ for some $\theta\in (0,1]$ and $R>0$
imposes  restrictions on the tail of the density $f$.
In particular,
the set of densities, uniformly bounded and
compactly supported on a cube of $\bR^d$,
is embedded in the set $\bG_\theta(\cdot)$
for any $\theta\in (0,1]$ (for details, see Section~\ref{subsec:(ii)-th-4}).
We will refer to the assumption $f\in \bG_\theta(R)$ as {\em the tail dominance condition}.
\par
In this section we study the problem of  adaptive density
estimation under the tail dominance condition.
We show that under this condition
the minimax rate of convergence can be essentially improved in the tail
zone.
In particular, if $\theta\leq \theta^*$ for some $\theta^*<1$ given below
then the tail zone disappears.
\par\smallskip
For any $\theta\in (0,1]$ let
\[
\nu^*(\theta)=\max\left\{\frac{1-\theta/p}{1-\theta/s+1/\beta },\;\frac{1}{2+1/\beta}\right\},
\]
and define
\begin{equation}\label{eq:nu-theta}
 \nu(\theta) =\left\{\begin{array}{cc}
           \nu^*(\theta),&  p \leq s(2+1/\beta),
           \\*[2mm]
\nu,  & p> s(2+ 1/\beta),
\end{array}
\right. \qquad
\mu_n(\theta)=\left\{
\begin{array}{ll}
(\ln n)^{1/p}, & p \in\{ \frac{2+1/\beta}{1/\theta+1/s},\; s(2+1/\beta)\},
\\*[2mm]
 \quad                  1, & {\rm otherwise},
                 \end{array}
\right.
\end{equation}
where $\nu$ is defined in (\ref{eq:nu}).

\begin{theorem}
The following statements hold.
\label{th:upper-new}
\begin{itemize}
\item[{\rm (i)}] For any $\theta\in (0,1]$ and $R>0$, Theorem~\ref{th:main}
remains valid if one replaces
$\bN_{\vec{r},d}(\vec{\beta},\vec{L}, M)$
by
$ \bG_\theta(R)\cap\bN_{\vec{r},d}(\vec{\beta},\vec{L}, M)$,
$\nu$ by
$\nu(\theta)$ and $\mu_n$ by $\mu_n(\theta)$. The constant $C$ may depend on $\theta$ and $R$.
\item[{\rm (ii)}]
For any $\theta\in (0,1]$, $\vec{\beta},\vec{L}\in (0,\infty)^d$,
$\vec{r}\in [1,\infty]^d$ and $M>0$ one can find $R>0$
such that
Theorem~\ref{th:lower-bound-in-L_p}
remains valid if one replaces
$\bN_{\vec{r},d}(\vec{\beta},\vec{L}, M)$
by
$\bG_\theta(R)\cap\bN_{\vec{r},d}(\vec{\beta},\vec{L}, M)$,
$\nu$ by
$\nu(\theta)$, and $\mu_n$ by $\mu_n(\theta)$.
\end{itemize}
\end{theorem}
\begin{remark}\mbox{}
{\rm
\begin{enumerate}
\item[{\rm 1.}]
The tail dominance condition  leads to
improvement of the rates of convergence in the whole tail zone.
In particular, if $f\in \bG_1(R)$ 
then the additional $\ln^{\frac{d}{p}}(n)$-factor disappears, cf. $\mu_n$ and $\mu_n(1)$.
Moreover, 
under the  tail dominance condition with $\theta<1$
the faster
convergence rate of the dense zone is achieved 
over a wider range of values of $p$,
$\frac{2+1/\beta}{1/\theta+1/s} \leq p \leq s(2+1/\beta)$. Additionally,
if
\[
 \theta < \theta^*:=\frac{ps}{s(2+1/\beta)-p},
\]
then the tail zone disappears. Note that $\theta^*\in (0,1)$
whenever $p\leq \frac{2+1/\beta}{1+1/s}$.
\par
As it was mentioned above,
the set of uniformly bounded and
compactly supported on a cube of $\bR^d$ densities
is embedded in the set $\bG_\theta(\cdot)$
for any $\theta\in (0,1]$. This fact explains why
the tail zone does not appear
in problems of estimating compactly supported densities.
\item[{\rm 2.}]
 We would like to emphasize that the couple $(\theta, R)$
is not used in the construction of the estimation procedure; thus,
our estimator is adaptive with respect to $(\theta,R)$ as well.
In particular, if the tail dominance condition does not hold,
our estimator achieves
the rate of Theorem~\ref{th:main}.
On the other hand, if this assumption holds, the rate of convergence is improved
automatically in the tail zone.
\item[{\rm 3.}]
The second statement
of the theorem is proved under
assumption that $R$ is large enough.
The fact that $R$ cannot be chosen  arbitrary small
is not technical; the parameters
$\vec{\beta}$, $\vec{L}$, $\vec{r}$, $M$,$\theta$ and $R$ are related to each other.
In particular, one can easily provide lower bounds on $R$ in terms of the
other parameters of the class. For instance, by the Lebesgue differentiation theorem,
$f(x)\leq f^*(x)$ almost everywhere; therefore
for any density $f\in \bG_\theta(R)$ such that $\|f\|_\infty\leq M$
one has
$$
1=\int f\leq M^{1-\theta}\|f^*\|^\theta_\theta\leq M^{1-\theta} R^{\theta} \;\;\Rightarrow\;\;
R\geq M^{1-1/\theta}.
$$
\par
Another lower bound on $R$ in terms of $\vec{L}$, $\vec{r}$ and $\theta$ can be established
using the Littlewood interpolation inequality [see, e.g., \cite[Section~5.5]{Garling}].
Let $0<q_0<q_1$ and $\alpha\in (0,1)$ be arbitrary numbers; then the
Littlewood inequality
states that $\|g\|_q\leq \|g\|_{q_0}^{1-\alpha}\|g\|_{q_1}^{\alpha}$,
where $q$ is defined by relation
 $\frac{1}{q}=\frac{1-\alpha}{q_0}+\frac{\alpha}{q_1}$.
Now, suppose that $f\in\bG_\theta(R)\cap\bN_{\vec{r},d}(\vec{\beta},\vec{L})$, and choose
$q_0=\theta$, $q_1=r_i$ and $\alpha=\frac{1-\theta}{1-\theta/r_i}$; then $q=1$ and
by the Littlewood inequality we have
$$
1=\|f\|_1\leq
\|f\|_\theta^{1-\alpha}\|f\|_{r_i}^\alpha\leq
R^{\frac{r_i\theta-\theta}{r_i-\theta}}L_{i}^{\frac{r_i-r_i\theta}{r_i-\theta}},\quad
i=1,\ldots, d\;\;\Rightarrow\;\;\;R\geq \max_{i=1,\ldots,d}
L_i^{\frac{r_i\theta-r_i}{r_i-\theta}}~.
$$
\item[{\rm 4.}] Another interesting observation is related
to the specific case $p=1$. Recall
that the condition $f\in \bN_{\vec{r},d}(\vec{\beta},\vec{L},M)$
alone is not sufficient for existence of consistent estimators.
%, see \cite{Ibr-Has2}.
However,
for any  $\theta\in (0,1)$  we can show
\[
 \inf_{\tilde{f}}
\cR_1^{(n)}\big[\tilde{f};\;
\bG_\theta(R)\cap \bN_{\vec{r},d}(\vec{\beta},\vec{L}, M)\big]
\;\leq\;C\bigg[\frac{L_\beta (\ln n)^d}{n}
\bigg]^{\frac{1-\theta}{1-\theta/s+1/\beta}}\;\to\; 0,
\;\;n\to\infty.
\]
This result follows from the proof of
Theorem~\ref{th:upper-new} and
(\ref{eq:weak-max}).
\end{enumerate}
}
\end{remark}
\par\smallskip
Now we argue that condition $f\in \bG_{\theta^*} (R)$
is, in a sense,
the weakest possible
ensuring
the ``usual'' rate of convergence,
corresponding to the index $\nu=\beta/(2\beta+1)$,
in the whole zone $p\leq s(2+1/\beta)$.
Indeed,  in view of Theorem~\ref{th:upper-new},
the minimax rate of convergence on the class
$\bG_{\theta^*}(R)\cap \bN_{\vec{r},d}(\vec{\beta},\vec{L}, M)$, say
$\overline{\psi}_n(\theta^*)$, satisfies
\begin{eqnarray*}
&&
c\big(L_\beta/n\big)^{\frac{\beta}{2\beta+1}}\;\leq\;
\overline{\psi}_n(\theta^*)\;\leq\;
C(\ln{n})^{1/p}\big(L_\beta\ln{n}/n\big)^{\frac{\beta}{2\beta+1}},
\end{eqnarray*}
where the constants $c$ and $C$ may depend on
$R$. On the other hand, if $\underline{\psi}_n(\theta^*)$ denotes
the minimax rate of convergence on the class
$\bN_{\vec{r},d}(\vec{\beta},\vec{L},M)\setminus\bG_{\theta^*}(R)$
then
\begin{eqnarray}
\label{eq:lower-bound-new}
c\big(L_\beta/n\big)^{\frac{1-1/p}{1-1/s+1/\beta}}
\;\leq\;
\underline{\psi}_n(\theta^*)
\;\leq\;
C\big(L_\beta\ln{n}/n\big)^{\frac{1-1/p}{1-1/s+1/\beta}},
\end{eqnarray}
provided that $p\leq \frac{2+1/\beta}{1+1/s}$.
The upper bound in  (\ref{eq:lower-bound-new}) is one of the statements of
Theorem \ref{th:main}, while
the lower bound is proved
in  Section~\ref{sec:proof-of-eq:lower-bound-new}.
\par
Thus we conclude that there is no tail zone
in estimation over the
class $G_{\theta^*}(R)\cap \bN_{\vec{r},d}(\vec{\beta}, \vec{L},M)$,
and this zone appears when considering
the class
$\bN_{\vec{r},d}(\vec{\beta}, \vec{L},M)\setminus \bG_{\theta^*}(R)$.
In this sense $f\in
\bG_{\theta^*}(R)\cap \bN_{\vec{r}}(\vec{\beta},\vec{L},M)$
is the necessary and sufficient condition
eliminating the tail zone.
\section{Proof of Theorem~\ref{th:oracle-inequality}}
\label{sec:Proof-th-3}
First we state two auxiliary results, Lemmas~\ref{lem:proportional} and~\ref{lem:1},
and  then turn to the proof of the theorem.
Proof of measurability of our estimator and proofs of
Lemmas~\ref{lem:proportional} and~\ref{lem:1} are given
in Appendix~A.
\subsection{Auxiliary lemmas}\label{sec:aux-th-1}
For any $g:\bR^d\to\bR^1$ denote
$$
 \check{M}_h(g,x)
% \frac{1}{4} \hat{M}_h(g,x)
 =\sqrt{\frac{\kappa\hat{A}_h(g,x)\ln n}{nV_h}} +
\frac{\kappa \ln n}{nV_h}~.
$$
\begin{lemma}\label{lem:proportional}
Let
$\chi_h(g,x)= \big[|\hat{A}_h(g,x)-A_h(g,x)|- M_h(g,x)\big]_+$, $h\in \cH$;
then
\begin{eqnarray*}%\label{eq:M-proportional}
[\check{M}_h(g,x)- 5M_h(g,x)]_+ \leq \frac{1}{2} \chi_h(g,x),\;\;\;\;
[M_h(g,x)- 4\check{M}_h(g,x)]_+ \leq 2 \chi_h(g,x).
\end{eqnarray*}
%In particular, (\ref{eq:M-proportional}) holds with $c_1=\sqrt{7}+\frac{3}{2}$ and $c_2=2$.
\end{lemma}
The next lemma establishes moment bounds on the following four random variables:
\begin{equation}\label{eq:zetas}
\begin{array}{lll}
&&
{\displaystyle \zeta_1(x)= \sup_{h\in \cH} [|\xi_h(x)|- M_h(K,x)]_+;\;\;}
\\
&&{\displaystyle \zeta_2(x)=\sup_{h,\eta\in \cH} [|\xi_{h,\eta}(x)|-M_{h\vee \eta}(Q,x)]_+;}
\\
&&{\displaystyle \zeta_3(x):=\sup_{h\in \cH} [ |A_h(K,x)-\hat{A}_h(K,x)|-M_h(K,x)]_+;\;}
\\
&&
{\displaystyle \zeta_4(x):=\sup_{h\in \cH} [ |A_h(Q,x)-\hat{A}_h(Q,x)|-M_h(Q,x)]_+.}
\end{array}
\end{equation}
\par
Denote $\rk_\infty=\|K\|_\infty\vee 1$ and
\[
%A(x)=\sup_{h\in \cH} A_h(K, x)\;\vee\;\sup_{h\in \cH} A_{h}(Q,x).
F(x)=\int \mathbf{1}_{[-1,1]^d}(t-x)f(t)\rd t.
\]
\begin{lemma}\label{lem:1}
Let  $q\geq 1$, $l\geq 1$ be arbitrary  numbers.
If
$\kappa \geq   \rk^2_\infty [(2q+4)d + 2l]$ then for all $ x\in \bR^d$
\begin{eqnarray}
\bE_f [\zeta_j(x)]^q  &\leq&  C_0
n^{-q/2}
\big\{ F(x)\vee n^{-l}\big\},\quad j=1,2,3,4,
\label{eq:ineq-1}
\end{eqnarray}
where constant $C_0$ depends on $d$, $q$, and $\rk_\infty$ only.
\end{lemma}

\subsection{Proof of oracle inequality (\ref{eq:pointwise-oracle})}
We recall the standard error decomposition of the kernel estimator:
for any $h\in \cH$ one has
\[
|\hat{f}_h(x) - f(x)| \leq |B_h(f,x)| + |\xi_h(x)|,
\]
where $B_h(f,x)$ and $\xi_h(x)$ are given in (\ref{eq:B-h}) and (\ref{eq:xi-h}) respectively.
Similar error decomposition holds for  auxiliary estimators $\hat{f}_{h,\eta}(x)$;
the corresponding  bias and
stochastic error are denoted by $B_{h,\eta}(f,x)$ and $\xi_{h,\eta}(x)$.
\par
$1^0$. The following relation for the bias
$B_{h,\eta}(f,x)$ of $\hat{f}_{h,\eta}(x)$ holds:
\begin{eqnarray}\label{eq:bias-relation}
B_{h,\eta}(f,x) - B_\eta(f,x) = \int K_{\eta}(t-x) B_h(f,t)\rd t,\;\;\;\forall h, \eta\in \cH.
\end{eqnarray}
Indeed, using the Fubini theorem and the fact that $\int K_h(x)\rd x=1$ for all $h\in \cH$
we have
\begin{eqnarray*}
 \int [K_h * K_\eta](t-x) f(t) \rd t &=& \int \biggl[\int K_h(t-y) K_\eta(y-x)\rd y\biggr]f(t)\rd t
\\
&=& \int K_\eta(y-x) f(y) \rd y
\\
&&\;\;+\;\;
\int K_\eta (y-x) \biggl[\int K_h(t-y)[f(t)-f(y)]\rd t\biggr] \rd y.
\end{eqnarray*}
It remains to note that
$ \int K_h(t-y)[f(t)-f(y)]\rd t=B_h(f,y)$ and to subtract  $f(x)$ from the both sides of the above equality.
Thus, (\ref{eq:bias-relation}) is proved.
\par\medskip
$2^0$. By the triangle inequality we have for any $h\in \cH$
\begin{equation}\label{eq:decomposition}
 |\hat{f}_{\hat{h}}(x)-f(x)| \;\leq\; |\hat{f}_{\hat{h}}(x)- \hat{f}_{\hat{h}, h}(x)| +
|\hat{f}_{\hat{h},h}(x) - \hat{f}_h(x)| + |\hat{f}_h(x) -f(x)|.
\end{equation}
We bound each term on the right hand side separately.
\par
First we note that, by  (\ref{eq:bias-relation}) and (\ref{eq:maximal-bias}), for any $h\in \cH$
\begin{eqnarray*}
\hat{R}_h(x) - \sup_{\eta\geq h} \hat{M}_\eta(Q,x) - \hat{M}_h(K,x) \;=\;
\sup_{\eta\in \cH}\Big[|\hat{f}_{h,\eta}(x)-\hat{f}_\eta(x)| -
\hat{M}_{h\vee \eta}(Q,x)-\hat{M}_\eta(K,x)\Big]_+
&&
\\
\;\;\leq  \bar{B}_h(f,x) +
\sup_{\eta\in \cH}\Big[|\xi_{h,\eta}(x)-\xi_\eta(x)| -
\hat{M}_{h\vee \eta}(Q, x)-\hat{M}_\eta(K, x)\Big]_+.
&&
\end{eqnarray*}
Thus, for any $h\in \cH$
\begin{equation}
\label{eq:bound-on-R}
 \hat{R}_h(x) \;\leq\; \bar{B}_h(f,x) + 2\hat{\zeta}(x) + \hat{M}_h(K,x) +
\sup_{\eta\geq h} \hat{M}_\eta(Q,x),
\end{equation}
where we put
$$
 \hat{\zeta} (x) := \sup_{h,\eta\in \cH}\big[|\xi_{h,\eta}(x)|-
\hat{M}_{h\vee \eta}(Q, x)\big]_+ \;\vee\;\sup_{h \in \cH}\big[|\xi_{h}(x)|-
\hat{M}_{h}(K,x)\big]_+.
$$

\smallskip

Second,
by (\ref{eq:bias-relation}) and    $\hat{f}_{h,\eta}\equiv \hat{f}_{\eta, h}$ for any $h, \eta \in \cH$ we have
\begin{eqnarray*}
&& |\hat{f}_{h,\eta}(x) - \hat{f}_h(x)| \;\leq\; |B_{\eta,h}(f,x)-B_{h}(f,x)| +
|\xi_{h,\eta}(x)-\xi_{h}(x)|
\\*[2mm]
&&\;\;\;\leq\; B_\eta(f,x) + \big[|\xi_{h,\eta}(x)-\xi_h(x)|- \hat{M}_{h\vee \eta}(Q,x)-
\hat{M}_h(K,x)\big] +
\sup_{\eta\geq h} \hat{M}_{\eta}(Q,x) + \hat{M}_h(K,x)
\\ &&\;\;\;\leq\; \bar{B}_\eta(f,x) + 2\hat{\zeta}(x) + \hat{R}_h(x),
\end{eqnarray*}
where the last inequality holds by definition of $\hat{R}_h(x)$ [see (\ref{eq:hat-R})].
There inequalities imply
the following upper bound on the first term on the right hand side of
(\ref{eq:decomposition}):
for any $h\in \cH$
\begin{eqnarray}
&& |\hat{f}_{\hat{h},h}(x) - \hat{f}_{\hat{h}}(x) | \;\leq\; \bar{B}_h(f,x) + \hat{R}_{\hat{h}}(x) + 2
\hat{\zeta}(x)
\nonumber
\\
&&\;\;\;\leq\; \bar{B}_h(f,x) + \hat{R}_h(x) + 2\hat{\zeta}(x)
\nonumber
\\
&&\;\;\;\leq\;
2\bar{B}_h(f,x) + 4\hat{\zeta}(x) + \sup_{\eta\geq h} \hat{M}_\eta (Q,x) +\hat{M}_h(K,x);
\label{eq:1}
\end{eqnarray}
where we have used the fact that $\hat{R}_{\hat{h}}(x) \leq \hat{R}_h(x)$ for all $h\in \cH$, and
inequality (\ref{eq:bound-on-R}).
\par
Now we turn to bounding the second term on the right hand side of (\ref{eq:decomposition}).
We get for any $h\in \cH$
\begin{eqnarray}
&& |\hat{f}_{\hat{h}, h}(x) - \hat{f}_h(x)| \;=\; |\hat{f}_{\hat{h}, h}(x) - \hat{f}_h(x)|
\pm \big[ \hat{M}_{\hat{h}\vee h}(Q,x) + \hat{M}_h(K,x)\big]
\nonumber
\\
&&\;\;\;\leq\; \hat{R}_{\hat{h}}(x) + \sup_{\eta\geq h} \hat{M}_\eta (Q,x)+\hat{M}_h(K,x)
\nonumber
\\
&&\;\;\;\leq\;
\bar{B}_h(f,x)+2\hat{\zeta}(x)+2\sup_{\eta\geq h} \hat{M}_\eta (Q,x)+2\hat{M}_h(K,x),
\label{eq:2}
\end{eqnarray}
where we again used (\ref{eq:bound-on-R}) and the fact that $\hat{R}_{\hat{h}}(x)\leq \hat{R}_h(x)$
for all $h\in \cH$.
\par
Finally
for any $h\in \cH$
\[
 |\hat{f}_h(x)-f(x)| \leq |B_h(f, x)| +|\xi_h(x)| \leq \bar{B}_h(f,x) + M_h(K,x) + \zeta(x).
\]
Thus,
combining (\ref{eq:1}), (\ref{eq:2}) and (\ref{eq:decomposition}) we obtain
\begin{equation}\label{eq:oracle-0}
 |\hat{f}_{\hat{h}}(x)-f(x)| \;\leq\; \inf_{h\in \cH} \Big\{4\bar{B}_h(f,x)+
3 \sup_{\eta\geq h} \hat{M}_\eta (Q, x)+ 3\hat{M}_h(K,x) + M_h(K,x)\Big\}
+ 6\hat{\zeta}(x) + \zeta(x).
\end{equation}
\par\medskip
$3^0$.
In order to complete the proof we note
that by the first inequality of  Lemma~\ref{lem:proportional} for any $g:\bR^d\to \bR^1$
\[
 \hat{M}_h(g,x) \leq 20 M_h(g,x) + 2 \chi_h(g,x).
\]
In addition, by the second inequality in Lemma~\ref{lem:proportional}
\begin{eqnarray*}
&& |\xi_h(x)| - \hat{M}_h(K,x)= |\xi_h(x)| - M_h(K,x) + M_h(K,x)- \hat{M}_h(K,x) \leq
\zeta(x) + 2\chi(x),
\\
&& |\xi_{h,\eta}(x)| - \hat{M}_{h\vee \eta}(Q,x) = |\xi_{h,\eta}(x)| - M_{h\vee \eta}(Q,x) +
M_{h\vee \eta}(Q,x) - \hat{M}_{h\vee \eta}(Q,x)\leq \zeta(x) + 2\chi(x),
\end{eqnarray*}
so that
$\hat{\zeta}(x) \leq \zeta(x)+ 2\chi(x)$.
Substituting these bounds in (\ref{eq:oracle-0}) we obtain
\[
 |\hat{f}(x)-f(x)| \leq \inf_{h\in \cH}
\Big\{4\bar{B}_h(f,x)+
60 \sup_{\eta\geq h} M_\eta (Q, x)+ 61M_h(K,x)\Big\}
+ 7\zeta(x) + 18\chi(x),
\]
as claimed.
\epr
\subsection{Proof of moment bounds (\ref{eq:integrability-of-zeta})}
\label{subsec:integrability}
%%%%%%%%%%%%%%%%%%%%%%%%%%%%%%%%%%%%%%%%%%%%%%%%%%%%%%%%%%%%%%%%%%%%%%%%%%%%%%
\iffalse
Let $f^*$ be the strong maximal function of $f$, see Section \ref{sec:max-f}.
First we note that
\[
\sup_{h\in \cH} \int |K_h(t-x)| f(t) \rd t \leq \rk_\infty f^*(x),\;\;\;
\sup_{h\in \cH} \int |Q_{h}(t-x)| f(t) \rd t \leq \rk^2_\infty f^*(x),
\]
so that $A(x)\leq \rk_\infty^2 f^*(x)$, $\forall x$.
\fi
%%%%%%%%%%%%%%%%%%%%%%%%%%%%%%%%%%%%%%%%%%%%%%%%%%%%%%%%%%%%%%%%%%%%%%%%%%%%%%
Let $\zeta_j(x)$, $j=1,\ldots, 4$ be defined by (\ref{eq:zetas}).
Then
\begin{eqnarray*}
\bE_f [\zeta_j(x)]^q
 &\leq&  C_0 n^{-q/2}
\big\{F(x)\vee n^{-l}\big\},
\end{eqnarray*}
as claimed in Lemma~\ref{lem:1}.
\par
Let $T_1=\{x\in \bR^d: F(x)\geq n^{-l}\}$ and $T_2=\bR^d\setminus T_1$.
Therefore
\begin{eqnarray}
\label{eq1:new-int}
&&\int_{T_1} \bE_f [\zeta_j(x)]^q \rd x \leq C_0 n^{-q/2}
 \int_{T_1}F(x)  \rd x\leq C_0 n^{-q/2}\int F(x)  \rd x
%\\
=
2^{d}C_0n^{-q/2}.
\end{eqnarray}
Now we analyze integrability on the set $T_2$.
We consider only the case $j=1,2$
since computations for $j=3,4$ are the same as for $j=1$.
\par
Let
$
U_{\max}(x)=[x-1, x+1]^d
%,\;\;\;\;U_{\max}:=U_{h^{\max}}.
$
and define the event
$D(x)=\big\{\sum_{i=1}^n {{\bf 1}}[X_i\in U_{\max}(x)]<2\big\}$, and let $\bar{D}(x)$
denote the complementary event.
First we argue that for $j=1,2$
\begin{equation}\label{eq:D(x)}
\zeta_j(x) {\mathbf{1}}\{D(x)\}=0,\;\;\;\forall x\in T_2.
\end{equation}
Indeed,
if $x\in T_2$ then for any $h\in\cH$
\begin{eqnarray*}
\left|\bE_f K_h(X_i-x)\right|  \leq n^{d}\rk_\infty F(x)\leq \rk_\infty n^{d-l},
\quad
\left|\bE_f Q_h(X_i-x)\right|  \leq n^{d}\rk^{2}_\infty F(x)\leq \rk^{2}_\infty n^{d-l}.
\end{eqnarray*}
Here we have used that $\cH=\big[1/n,1]^d$ and that $\text{supp}(K)=[-1/2,1/2]^d,\;\text{supp}(Q)=[-1,1]^d$.

Hence, by definition of $\xi_h(x)$, for any $h\in \cH$ one has for any $l\geq d+1$
\begin{eqnarray*}
 |\xi_h(x)| {\bf 1}\{D(x)\}\;\leq\; \Big|\frac{1}{n}\sum_{i=1}^n K_h(X_i-x)\Big| {\bf 1}\{D(x)\} +
\rk_\infty n^{d-l}
\\
 \leq\;
\frac{2\rk_\infty}{nV_h} + \rk_\infty n^{d-l}
\;\leq\; \frac{4\rk_\infty}{nV_h}\leq M_h(K,x),
\end{eqnarray*}
where we have used that $n^{d-l}\leq (nV_h)^{-1}$ for $l\geq d+1$, $\kappa\ln n\geq 4\rk_\infty$
by the condition on $\kappa$
[see also definition of $M_h(K,x)$], and $n\geq 3$.
Therefore $\zeta_1(x){\bf 1}\{D(x)\}=0$ for $x\in T_2$.
By the same reasoning for $\zeta_2(x)$ we obtain that
$\zeta_2(x){\bf 1}\{D(x)\}=0$, $\forall x\in T_2$ because $\kappa\ln n \geq 4\rk^2_\infty$. Thus
(\ref{eq:D(x)}) is proved.
Using (\ref{eq:D(x)}) we can write
\begin{eqnarray}
\int_{T_2}  \bE_f [\zeta_j(x)\big]^q {\bf 1}\{\bar{D}(x)\}
\rd x
&\leq& \int_{T_2}  \bE_f \Big(\big[\sup_{h\in \cH} |\xi_h(x)|^q
\vee \sup_{h,\eta\in \cH} |\xi_{h,\eta}(x)|^q\big]
{\bf 1}\{\bar{D}(x)\}\Big)
\rd x
\nonumber
\\
&\leq& \big(2\rk^2_\infty n^d\big)^q \int_{T_2} \bP_f\{\bar{D}(x)\} \rd x.
\label{eq:E-T2}
\end{eqnarray}
Now we bound from above the integral on the right hand side of the last display formula.
For any $z>0$ we have in view of the exponential  Markov inequality
\begin{eqnarray*}
\bP_f\{\bar{D}(x)\} &=& \bP_f\Big\{ \sum_{i=1}^n {\bf 1}[X_i\in U_{\max}(x)] \geq 2\Big\}
\leq e^{-2z} \big[e^{z} F(x) + 1-F(x)\big]^n
\\
&=&
e^{-2z}\big[(e^z-1)F(x)+1]^n \leq \exp\{-2z + n(e^{z}-1)F(x)\}.
\end{eqnarray*}
Minimizing the right hand side w.r.t. $z$ we find $z=\ln 2 -\ln{\{nF(x)\}}$ and, therefore,
\[
 \bP_f\big\{\bar{D}(x)\big\} \leq
4^{-1}n^2F^2(x) \exp\{2-nF(x)\}\leq (e^2/4) n^2 F^2(x).
\]
Since $F(x)\leq n^{-l}$ for any $x\in T_2$
we obtain
\begin{eqnarray*}
 \int_{T_2} \bP_f\{\bar{D}(x)\} \rd x \leq  (e^2/4)n^{2-l} \int F(x)\rd x
=2^d \big(e^2/4\big)  n^{2-l}.
\end{eqnarray*}
Combining this inequality with (\ref{eq:E-T2}) we obtain
\begin{eqnarray}
\label{eq2:new-int}
&&\int_{T_2}  \bE_f [\zeta_j(x)\big]^q {\bf 1}\{\bar{D}(x)\}
\rd x \leq 2^d (2\rk^2_\infty)^q(e^2/4)   n^{2+dq-l}.
\end{eqnarray}
Choosing $l=(d+1)q+2$ we come to the assertion of the theorem in view of
(\ref{eq1:new-int})
and (\ref{eq2:new-int}).
\epr

\section{Proofs of Theorem~\ref{th:main} and statement~(i) of
Theorem~\ref{th:upper-new}}
\label{sec:proofs}
The proofs of Theorem~\ref{th:main} and
of statement~(i) of Theorem~\ref{th:upper-new}
go along similar lines. That is why
we state our auxiliary results
(Propositions~\ref{prop:first-bound} and~\ref{prop:new})
in the form that is suitable
for the use in the proof of Theorem~\ref{th:upper-new}.
%For this purpose it will be convenient to extend the definition of the class$\bG_\theta(R)$ [see (\ref{eq:G-theta})]
%to the case $\theta=1$.In the sequel by $\bG_1(R)$ we mean the set of all probability densities on $\bR^d$,
%no matter what the value of $R$ is.
\par
This section is organized as follows.
First, in Subsection~\ref{subsec:preliminaries}
we present and discuss some facts from
functional analysis.
Then in Lemma~\ref{lem:bias-norm-bound} of Subsection~\ref{sec:approximation-K}
we state an auxiliary result on approximation properties of the kernel $K$
defined in (\ref{eq:w-function}).
Proof outline and notation are discussed
in Subsection~\ref{subsec:outline+notation}.
Subsection~\ref{subsec:aux-proposition}
presents two auxiliary propositions,
and the proofs of Theorem~\ref{th:main}
and statement~(i) of Theorem~\ref{th:upper-new}
are completed in Subsections~\ref{sec:2} and~\ref{sec:upper-4}.
Proofs of the auxiliary results, Lemma~\ref{lem:bias-norm-bound}
and Propositions~\ref{prop:first-bound} and~\ref{prop:new} are given in Appendix~B.
\par
In the subsequent proof
$c_i,C_i,\bar{c}_i, \bar{C}_i, \hat{c}_i,\hat{C}_i, \tilde{c}_i, \tilde{C}_i, \ldots$,
stand for constants that can depend on $L_0$, $M$
$\vec{\beta}$, $\vec{r}$, $d$ and $p$, but
are independent of $\vec{L}$ and $n$. These constants can be different on
different appearances. In the case when the assumption
$f\in\bG_\theta(R)$ with $\theta\in (0,1]$ 
is imposed, they may also depend on $\theta$ and $R$.

\subsection{Preliminaries}
\label{subsec:preliminaries}
We present
an embedding theorem for the anisotropic Nikol'skii classes and discuss
some properties of the strong
maximal operator.
\subsubsection{Embedding theorem}
The statement given below in (\ref{eq:embedd-nik})
is a particular case of the embedding theorem for anisotropic Nikol'skii classes
$\bN_{\vec{r},d}(\vec{\beta},\vec{L})$; see
\cite[Section~6.9.1.]{Nikolski}.
\par
For the fixed class parameters $\vec{\beta}$ and $\vec{r}$ define
\[
\tau(p)=1-\sum_{j=1}^d\frac{1}{\beta_j}\left(\frac{1}{r_j}-\frac{1}{p}\right),\quad
\tau_i=1-\sum_{j=1}^d\frac{1}{\beta_j}\left(\frac{1}{r_j}-\frac{1}{r_i}\right),\;\;i=1,\ldots,d,
\]
and put
\begin{equation}
\label{eq:gamma-and-q}
q_i=r_i\vee p,\;\;\;
\qquad\;\;
\gamma_i=\left\{
\begin{array}{ll}
\frac{\beta_i\tau(p)}{\tau_i},\quad & r_i<p,
\\
\beta_i,\quad & r_i\geq p.
\end{array}
\right.
\end{equation}
Let $\tau(p)>0$ and $\tau_i>0$ for all $i=1,\ldots, d$; then for any $p\geq 1$ one has
\begin{equation}
\label{eq:embedd-nik}
 \bN_{\vec{r},d}\big(\vec{\beta},\vec{L}\big) \subseteq
\bN_{\vec{q},d}\big(\vec{\gamma},c\vec{L}\big),
\end{equation}
where constant $c>0$ is independent of $\vec{L}$ and $p$.
%
\iffalse
\paragraph{Strong maximal function}
\label{sec:max-f}
Let $g:\bR^d\to\bR$  be a locally integrable function, and let $g^*$ denote
the strong maximal function of $g$, see
(\ref{eq:maximal-function}).
It is worth noting that the {\em Hardy--Littlewood maximal function} is defined by
(\ref{eq:maximal-function}) with the supremum taken over all cubes with sides parallel to the
coordinate axes,
centered at $x$.
\par
It is well known that the strong maximal operator $g\mapsto g^*$
is of the strong $(p,p)$--type for all $1<p\leq \infty$, i.e.,
if $g\in \bL_p(\bR^d)$ then $g^*\in \bL_p(\bR^d)$ and there exists a constant $\bar{C}$
depending on $p$ only such that
\begin{equation}\label{eq:strong-max}
 \|g^*\|_p \leq \bar{C} \|g\|_p,\;\;\;p\in (1,\infty].
\end{equation}
In distinction to the Hardy--Littlewood maximal function,
the strong maximal operator is not of
the weak (1,1)--type.
In fact, the following statement holds: there exists constant $C$ depending on $d$
only such that
\begin{equation}\label{eq:weak-max}
 \big|\{x: g^*(x)\geq \alpha\}\big| \leq C \int \frac{|g(x)|}{\alpha}
\bigg\{ 1 +
\bigg(\ln_+\frac{|g(x)|}{\alpha}\bigg)^{d-1}\bigg\}\rd x,\quad\forall\alpha>0.
\end{equation}
We refer to \cite{Guzman} for more details.
%
\fi
\subsubsection{Strong maximal function}
\label{sec:max-f}
Let $g:\bR^d\to\bR$  be a locally integrable function. We define  the strong maximal function $g^\star$ of $g$ by formula
\begin{equation}\label{eq2:maximal-function}
 g^\star(x):= \sup_{H} \frac{1}{|H|} \int_H g(t) \rd t,\;\;\;x\in \bR^d,
\end{equation}
where the supremum is taken over all
possible rectangles $H$  in $\bR^d$ with sides parallel to the coordinate
axes, containing point $x$.
It is worth noting that the {\em Hardy--Littlewood maximal function} is defined by
(\ref{eq2:maximal-function}) with the supremum taken over all cubes with sides parallel to the
coordinate axes,
centered at $x$.
\par
It is well known that the strong maximal operator $g\mapsto g^\star$
is of the strong $(p,p)$--type for all $1<p\leq \infty$, i.e.,
if $g\in \bL_p(\bR^d)$ then $g^\star\in \bL_p(\bR^d)$ and there exists a constant $\bar{C}$
depending on $p$ only such that
\begin{equation*}
 \|g^\star\|_p \leq \bar{C} \|g\|_p,\;\;\;p\in (1,\infty].
\end{equation*}
Let $g^*$ be defined  defined in
(\ref{eq:maximal-function}). Since obviously $g^*(x)\leq g^\star(x)$ for all $x\in\bR^d$  we have
\begin{equation}\label{eq:strong-max}
 \|g^*\|_p \leq \bar{C} \|g\|_p,\;\;\;p\in (1,\infty].
\end{equation}

In distinction to the Hardy--Littlewood maximal function, the strong maximal operator is not
of
the weak (1,1)--type.
In fact, the following statement holds: there exists constant $C$ depending on $d$
only such that
\begin{equation}\label{eq:weak-max}
 \big|\{x: g^\star(x)\geq \alpha\}\big| \leq C \int \frac{|g(x)|}{\alpha}
\bigg\{ 1 +
\bigg(\ln_+\frac{|g(x)|}{\alpha}\bigg)^{d-1}\bigg\}\rd x,\quad\forall\alpha>0.
\end{equation}
We refer to \cite{Guzman} for more details.
\subsection{Approximation properties of kernel $K$}
\label{sec:approximation-K}
%\subsubsection{Bias representation and its properties}
%Let $B_h(f,x)$ stand for the bias of kernel estimator $\hat{f}_h(x)$
%associated with kernel $K$ given in (\ref{eq:w-function}); see (\ref{eq:B-h}).
%\begin{equation*}%\label{eq:bias-kernel-est}
% B_h (f,x)= \int K_h (t-x) \big[f(t)-f(x)\big] \rd t,\;\;\;x\in \bR^d.
%\end{equation*}
%\par
The
next lemma establishes an upper bound on norm of
the bias $B_h(f,\cdot)$ of kernel estimator $\hat{f}_h$
when $f$ belongs to the anisotropic Nikol'skii class.
\begin{lemma}\label{lem:bias-norm-bound}
 Let $f\in \bN_{\vec{r},d}\big(\vec{\beta},\vec{L}\big)$. Let $\hat{f}_h$ be the
estimator  (\ref{eq:kernel-est})
associated with kernel (\ref{eq:w-function})
%(\ref{eq:kernel})
with $\ell> \max_{j=1,\ldots,d}\beta_j$.
Then $B_h(f,x)$ can represented as the sum
$B_h(f, x)=\sum_{j=1}^d B_{h,j}(f,x)$
with  functions $B_{h,j}(f,x)$ satisfying the following inequalities:
\begin{equation}\label{eq:bias-norms-1}
 \big\|B_{h,j}(f, \cdot)\big\|_{r_j} \leq C_1L_j h_j^{\beta_j},\;\;\;\forall j=1,\ldots,d.
\end{equation}
Moreover,  if  $s \geq 1$, then for any $p\geq 1$
%\begin{equation}\label{eq:beta-prime}
%\gamma_j :=\beta_j \frac{\tau(p)}{\tau_j},\;\;\;j=1,\ldots,d,
%\end{equation}
%where
%\begin{equation}
%\label{eq:tau}
%$\tau(p)=1-1/s+1/\beta p,\;\tau_j=1-1/s+1/\beta r_j,\;j=\overline{1,d}.$
%\end{equation}
%If $\tau(p)>0$ and $\tau_j>0$ for all $j=1,\ldots,d$
%Then, if $s>1$
\begin{equation}\label{eq:bias-norms-2}
\big\|B_{h,j}(f,\cdot)\big\|_{q_j} \leq C_2 L_j h_j^{\gamma_j} ,\;\;\;\forall j=1,\ldots, d,
\end{equation}
where $\vec{\gamma}=\vec{\gamma}(p)$ and $\vec{q}=\vec{q}(p)$ are defined in (\ref{eq:gamma-and-q}).
Here $C_1$ and  $C_2$ are constants
independent of $\vec{L}$ and $p$.
\end{lemma}
\subsection{Proof outline and notation}\label{subsec:outline+notation}
%\subsubsection{Preliminary remarks}
The starting point of our proof is
the pointwise oracle inequality
(\ref{eq:pointwise-oracle}) together with the moment bound (\ref{eq:integrability-of-zeta}).
Denote
\begin{equation}\label{eq:U-bar}
\bar{U}_f(x)= \inf_{h\in\cH}\Big\{\bar{B}_h(f,x)+\sup_{\eta\geq h}
M_\eta (K\vee Q, x)\Big\};
\end{equation}
then, taking into account that
$M_\eta(K\vee Q, x)$ is greater than $M_\eta(K,x)$ and $M_\eta(Q,x)$ for any $x$ and $\eta$
[see (\ref{eq:A-A}) and  (\ref{eq:hat-A-M})], and using
(\ref{eq:pointwise-oracle}), we have
\begin{equation*}%\label{eq:oracle-inequality}
 |\hat{f}(x)-f(x)| \leq c_0 [\bar{U}_f(x)
 + \omega(x)],
\end{equation*}
where $c_0$ is an absolute constant, and
$\omega(x):=\zeta(x)+\chi(x)$ with
$\zeta(x)$ and $\chi(x)$ defined in  (\ref{eq:zeta}) and (\ref{eq:chi}).
Therefore, by (\ref{eq:integrability-of-zeta}) applied with $q=p$ and by the Fubini theorem,
there exists
constant $\bar{c}_0>0$ such that for any probability density
$f$ and any Borel set $\cA\subseteq\bR^d$ one has
\begin{equation}
\label{eq00-new:oracle-inequality}
\bE_f\int_{\cA}|\hat{f}(x)-f(x)|^{p}\rd x \;\leq\; \bar{c}_0\bigg[\,
\int_{\cA}\bar{U}_f^p(x)\rd x +n^{-p/2}\,\bigg].
\end{equation}
Recall that $\rk_\infty=\|K\|_\infty\vee 1$; by definition of $\bar{B}_h(f,x)$ [see (\ref{eq:maximal-bias})]
and by Lemma~\ref{lem:bias-norm-bound} one has
\begin{equation*}%\label{eq:B-bar}
 \bar{B}_h(f,x)\leq \rk_\infty\sum_{j=1}^d B_{h,j}^*(f,x),
\end{equation*}
where $B_{h,j}^*(f,x)$ is the strong maximal function of $|B_{h,j}(f,x)|$, $j=1,\ldots,d$.
Therefore if we let
\begin{equation}\label{eq:U}
 U_f(x):= \inf_{h\in \cH}\Big\{\max_{j=1,\ldots,d} B^*_{h,j}(f,x) + \sup_{\eta\geq h} M_\eta(K\vee Q, x)
\Big\},
\end{equation}
then
\begin{equation}
\label{eq01:oracle-inequality}
\bar{U}_f(x)\leq \rk_\infty U_f(x),\quad\forall x\in\bR^d.
\end{equation}
\par
The key element of the proof
is derivation of upper bounds on the integral
\[
J:=\int_{\bR^d} U_f^p(x) \rd x.
\]
These bounds
will be established by division of
$\bR^d$ in ``slices'', and appropriate choice of bandwidth
$h\in \cH$   on every ``slice''. For this purpose
the following bounds on norms of $B^*_{h,j}(f, \cdot)$
will be used.
Inequality  (\ref{eq:strong-max}) and the first assertion of Lemma~\ref{lem:bias-norm-bound}
imply that  for any $p>1$, $\vec{r}\in (1,\infty]^{d}$ and
 any $f\in \bN_{\vec{r},d}\big(\vec{\beta},\vec{L}\big)$
one has
\begin{eqnarray}
&& \big\|B^*_{h,j}(f, \cdot)\big\|_{r_j} \leq \bar{c}_1L_j h_j^{\beta_j},\;\;\;
\forall j=1,\ldots, d,
\label{eq:B*-1}
\end{eqnarray}
Moreover, if
 $s\geq 1$ then, by
the second assertion of Lemma~\ref{lem:bias-norm-bound},
for any $p>1$, $\vec{r}\in (1,\infty]^{d}$ and
 $f\in \bN_{\vec{r},d}\big(\vec{\beta},\vec{L}\big)$
\begin{eqnarray}
&& \big\|B^*_{h,j}(f,\cdot)\big\|_{q_j} \leq \bar{c}_2 L_j h_j^{\gamma_j} ,\;\;\;\forall j=1,\ldots,d.
\label{eq:B*-2}
\end{eqnarray}
\par\medskip
Let
$\delta:=\ln n/n$, $\varphi:= (L_\beta \delta)^{\beta/(2\beta+1)}$.
Let $m_0(\theta),\;\theta\in (0,1],$ be an integer number to be specified later; see (\ref{eq:m-0}) below.
For $m\in \bZ$, $m\geq m_0(\theta)$ define ``slices''
\[
 \cX_m := \big\{x\in \bR^d: 2^m \varphi < U_f(x) \leq  2^{m+1} \varphi \big\},\;\;
\cX^-_{m_0(\theta)} := \big\{x\in \bR^d: U_f(x) \leq 2^{m_0(\theta)}\varphi\},
\]
and consider the corresponding integrals
\begin{equation*}%\label{eq:J-m-def}
 J_m:=\int_{\cX_m} U_f^p(x) \rd x,\;\;\;J_{m_0}^-:=\int_{\cX_{m_0(\theta)}^-} U_f^p(x)\rd x.
\end{equation*}
With this notation,
 using (\ref{eq00-new:oracle-inequality}) and (\ref{eq01:oracle-inequality})
we can write
\begin{eqnarray}
 \bE_f\|\hat{f}-f\|_p^p  &\leq&  \bE_f\int_{\cX_{m_0(\theta)}^{-}} |\hat{f}(x)-f(x)|^p \rd x \;+\;
\tilde{c}_1
\sum_{m=m_0(\theta)}^\infty \int_{\cX_m} U^p_f(x) \rd x \;+ \;\tilde{c}_2n^{-p/2}
\nonumber
\\
&=:&
J_{m_0(\theta)}^-
+\;\tilde{c}_1 \sum_{m=m_0(\theta)}^\infty J_m \;+\; \tilde{c}_2n^{-p/2}.
\label{eq:start}
\end{eqnarray}
The rest of the proof consists of bounding the integrals
$J_{m_0(\theta)}^-$ and $J_m$
on the  right hand side of (\ref{eq:start}) and combining these bounds in different zones.
\par
\par
The following notation will be used in the subsequent proof.
For the sake of brevity we will write
\[
M_\eta(x):=M_\eta(K\vee Q, x),\;\;\;A_\eta(x):=A_\eta(K\vee Q,x),\;\;\;\forall \eta\in \cH.
\]
We let $I:=\{1,\ldots,d\}$, and
\[
I_+:=\{j\in I: p\leq r_j<\infty\},\;\;I_-:=\{j\in I:1<r_j<p\},\;\;I_\infty:=\{j\in I:r_j=\infty\}.
\]
With
$\vec{\gamma}=(\gamma_1,\ldots,\gamma_d)$ and $\vec{q}=(q_1,\ldots,q_d)$ given by
(\ref{eq:gamma-and-q}) we  define  quantities $\gamma$, $\upsilon$ and $L_\gamma$
by the formulas
\begin{equation}\label{eq:gamma-upsilon}
\frac{1}{\gamma}:=\sum_{j=1}^d \frac{1}{\gamma_j},\qquad
\frac{1}{\upsilon}:=\sum_{j=1}^d \frac{1}{\gamma_j q_j},\qquad L_\gamma:=\prod_{j=1}^dL_j^{1/\gamma_j}.
\end{equation}
Note some useful inequalities between the quantities defined above.
First,
$\gamma_j < \beta_j$ for all $j\in I_-$ which is a consequence of the fact that $\tau(p)<\tau_j$
for $j\in I_-$. This implies
\begin{equation}
\label{eq:gamma<beta}
 \frac{1}{\gamma} - \frac{1}{\beta} = \sum_{j\in I_-} \Big(\frac{1}{\gamma_j}-\frac{1}{\beta_j}\Big) >0.
\end{equation}
Next, if $s\geq 1$ then
\begin{equation}
\label{eq:s<upsilon}
\frac{1}{s} > \frac{1}{\upsilon}.
\end{equation}
We have
\begin{equation*}
 \frac{1}{s} - \frac{1}{\upsilon} = \sum_{j\in I_-} \Big(\frac{1}{\beta_jr_j}-\frac{1}{\gamma_j p}\Big)
=\sum_{j\in I_-}\frac{1}{\beta_j}\Big(\frac{1}{r_j}-\frac{\tau_j}{\tau(p)p}\Big).
\end{equation*}
Hence  (\ref{eq:s<upsilon}) will be proved if we show that $r_j^{-1}\tau(p)p\geq \tau_j$ for all $j\in I_-$. Indeed,
$$
\frac{\tau(p)p}{r_j}=\frac{p(1-1/s)+1/\beta}{r_j}\geq 1-\frac{1}{s}+\frac{1}{\beta r_j}:=\tau_j,
$$
where to get the second inequality we have used that $r_j\leq p$ for any $j\in I_-$ and that $s\geq 1$.

\par
Finally, %At last
remark also that
\begin{equation}\label{eq003:proof-th}
p-\upsilon(2+1/\gamma)<0.
\end{equation}
Indeed,  since $r_j\geq p$ for any $j\in I_+\cup I_\infty$,
$$
\frac{p}{\upsilon}=\sum_{j\in I_+}\frac{p}{\beta_jr_j}+\sum_{j\in I_-}
\frac{1}{\gamma_j}\leq \sum_{j\in I_+}\frac{1}{\beta_j}+\sum_{j\in I_-}\frac{1}{\gamma_j} =\frac{1}{\gamma}.
$$
This  yields $p\leq\upsilon/\gamma$, and (\ref{eq003:proof-th}) follows.
\par
\subsection{Auxiliary results}
%Bounds on $J_{m_0}^-$ and $J_m$}
\label{subsec:aux-proposition}
%In the proof of the upper bound of Theorem~\ref{th:upper-new}
%parameter $m_0$ will be chosen dependent on $\theta$, and we will write
%$m_0=m_0(\theta)$. In the proof of Theorem~\ref{th:main}, with slight abuse of notation,
%$m_0$ will stand for $m_0(1)$.
%\par
For $\theta\in (0,1]$ and for some constant $\hat{c}_1>0$ define
\begin{equation}\label{eq:m-0}
 m_0(\theta):=\min\bigg\{ m\in \bZ:
2^{m_0(\theta)\left(
\frac{1-\theta/s+1/\beta}{1+\theta/s}
\right)}  >\hat{c}_1\kappa\varphi
\bigg\}.
\end{equation}
Note  that $1-\theta/s+1/\beta>0$
for any $\theta\in (0,1]$, since $s\geq\beta$ by $r_j\geq 1$, $j=1,\ldots, d$. Therefore
$m_0(\theta)<0$ for large  enough $n$.
\par
It will be convenient to introduce the following notation
\begin{equation}
 \label{eq:m-1-def}
 m_1:= \min\Big\{m\in \bZ:
2^{m\left[\upsilon(2+1/\gamma)-s(2+1/\beta)\right]}\geq
\big(L_\gamma/L_\beta\big)^{\upsilon} \varphi^{\upsilon(1/\beta-1/\gamma)}
\Big\}.
\end{equation}
It follows from this definition that
\begin{equation}\label{eq:m-1-bounds}
\Big[
 \big(L_\gamma/L_\beta\big) \varphi^{1/\beta-1/\gamma}
\Big]^{\frac{\upsilon}{\upsilon(2+1/\gamma)-s(2+1/\beta)}}\leq
2^{m_1}\leq
2\Big[
 \big(L_\gamma/L_\beta\big) \varphi^{1/\beta-1/\gamma}
\Big]^{\frac{\upsilon}{\upsilon(2+1/\gamma)-s(2+1/\beta)}}.
\end{equation}
In view of (\ref{eq:gamma<beta}) and (\ref{eq:s<upsilon})
\begin{equation}\label{eq:m-1-deff}
 \upsilon\Big(2+\frac{1}{\gamma}\Big)- s\Big(2+\frac{1}{\beta}\Big)
=
s\upsilon \bigg[\Big(2+\frac{1}{\beta}\Big)\Big(\frac{1}{s}-\frac{1}{\upsilon}\Big)+ \frac{1}{s}\Big(
\frac{1}{\gamma}-\frac{1}{\beta}\Big)\bigg]
%> s\Big(\frac{1}{\gamma}-\frac{1}{\beta}\Big)
>0;
\end{equation}
hence $m_1>1$ for large $n$.
\par
The bounds on $J_{m_0(\theta)}^-$ and $J_m$ are given in the next two propositions.
\par

\begin{proposition}\label{prop:first-bound}
There exist constants
$\hat{c}_1,\hat{c}_2>0$ and $\hat{C}_1,\hat{C}_2>0$
such that any $n$ large enough the following statements hold.
\begin{itemize}
\item[\rm (i)] For any probability density $f$ and any  $m_0(1)\leq m\leq 0$
\begin{eqnarray}\label{eq001:m-large}
&&
J_m \leq \hat{C}_1\,
2^{m\left(p-\frac{2+1/\beta}{1+1/s}\right)} \varphi^p.
\end{eqnarray}
\item[\rm (ii)]
Let $f\in \bG_\theta(R)$, $\theta\in (0,1]$; then for   any  $m_0(\theta)\leq m\leq 0$
one has
\begin{eqnarray}\label{eq01:m-large}
&&
J_m \leq \hat{C}_1\,
2^{m\left(p-\frac{2+1/\beta}{1/\theta+1/s}\right)} \varphi^p.
\end{eqnarray}
\item[{\rm (iii)}] For any $m\in \bZ$ satisfying
$1\leq 2^{m} \leq \hat{c}_2\varphi^{-1}$ and any probability density $f$ one has
\begin{eqnarray}\label{eq02:m}
&&
J_m\;\leq\;\hat{C}_2 2^{m [p- s(2+1/\beta)]} \varphi^p.
\end{eqnarray}
\item[{\rm (iv)}] Let $s\geq 1$; then  for any $m\in\bZ$ such that
%\begin{equation}\label{eq:m-range}
%2^{m\left[\upsilon(2+1/\gamma)-s(2+1/\beta)\right]}\geq
%\big(L_\gamma/L_\beta\big)^{\upsilon} \varphi^{\upsilon(1/\beta-1/\gamma)},\;\;
$
m\geq m_1$,
$
\;\;2^{m} \leq \hat{c}_2\varphi^{-1}
$
%\end{equation}
and any probability density $f$
one has
\begin{eqnarray}\label{eq03:m}
&&
J_m \;\leq\;\hat{C}_2\varphi^{p}
\bigg[\frac{L_\gamma\varphi^{1/\beta}}{L_\beta\varphi^{1/\gamma}}\bigg]^{\upsilon}
2^{m\left[p-\upsilon(2+1/\gamma)\right]}~.
\end{eqnarray}
% if, additionally, $s\geq 1$.
% $\vec{\beta}$ and $\vec{r}$ satisfy (\ref{eq:class-restriction}).
\end{itemize}
\end{proposition}
%
% For any $\theta\in (0,1]$ set
% $\cX^{(\theta)}_-:=\{x\in \bR^d: \bar{U}_f(x) \leq \rk_\infty 2^{m_0(\theta)}\varphi\}$.
\begin{proposition}
\label{prop:new}
There exist constants $\hat{C}_3, \hat{C}_4>0$ such that the following statements hold.
\begin{itemize}
\item[{\rm (i)}] Let $\nu$ is defined in (\ref{eq:nu}). Then for
all large enough $n$ and for any density $f$ one has
\begin{equation}
\label{eq:prop-new_1}
J^-_{m_0(1)}\;=\;\bE_f\int_{\cX^-_{m_0(1)}}\big|\hat{f}(x)-f(x)\big|^{p}\rd x\;\leq\;
 \hat{C}_3\ln^d(n)(L_\beta\delta)^{p\nu}.
\end{equation}
\item[{\rm (ii)}] Let  $\nu(\theta)$ is defined in (\ref{eq:nu-theta}). Then for  any $\theta\in (0,1]$ and for all $n$ large enough
\begin{equation}
\label{eq:prop-new_theta}
 \sup_{f\in\bG_\theta(R)}\bE_f\int_{\cX_{m_0(\theta)}^-}\big|\hat{f}(x)-f(x)\big|^{p}\rd x\leq
 \hat{C}_4(L_\beta\delta)^{p\nu(\theta)}.
\end{equation}
\end{itemize}
\end{proposition}

\subsection{Proof of Theorem~\ref{th:main}}\label{sec:2}
%As mentioned above, here $m_0=m_0(1)$.
Using (\ref{eq:start}) and inequality
(\ref{eq:prop-new_1}) of Proposition \ref{prop:new} we obtain
\begin{equation}\label{eq:start-1}
 \bE_f\|\hat{f}-f\|_p^p \;\leq\;  c_1\pi_n\big(L_\beta \delta\big)^{p\nu} +
c_2\sum_{m=m_0(1)}^\infty J_m,
\end{equation}
where $\pi_n=\ln^d(n)$ if $p\leq\frac{2+1/\beta}{1+1/s}$ and $\pi_n=1$ otherwise.

\par
We proceed with bounding
the second term on the right hand side of the last display formula. First,
because  $\|f\|_\infty \leq M$,
$$
\max_{j=1,\ldots,d}
\|B^*_{h,j}(f,\cdot)\|_\infty\leq 2^{d}M\rk^2_{\infty},\qquad
\sup_{\eta>0}\|A_\eta\|_\infty\leq 2^{d}M\rk^2_{\infty}.
$$
This  implies that there exists constant $c_3>0$ with the following property:
\[
 m_2:=\min \{m\in \bZ: 2^m \geq c_3\varphi^{-1}\} \;\;\Rightarrow\;\;J_m=0,\;\forall m\geq m_2.
\]
Thus the sum on right hand side of (\ref{eq:start-1})
extends from $m_0(1)$ to $m_2$.
\par\medskip
1$^0$. {\em Tail zone: $p<\frac{2+1/\beta}{1+1/s}$.}
Using bounds (\ref{eq001:m-large})  and (\ref{eq02:m}) of Proposition~\ref{prop:first-bound},
we obtain
\begin{eqnarray*}
 \sum_{m=m_0(1)}^\infty J_m \;\leq\; c_4\varphi^p\Big[
 \sum_{m=m_0(1)}^{0} 2^{m(p-\frac{2+1/\beta}{1+1/s})} + \sum_{m=1}^{m_2} 2^{m[p-s(2+1/\beta)]}\Big]
\;\leq\; c_5\,\varphi^p 2^{m_0(1)(p-\frac{2+1/\beta}{1+1/s})},
\end{eqnarray*}
where the last inequality follows from the fact that $m_0(1)<0$ and $p<\frac{2+1/\beta}{1+1/s}<s(2+1/\beta)$.
Using (\ref{eq:m-0}), after straightforward algebra we obtain  that
\[
 \sum_{m=m_0(1)}^\infty J_m \;\leq\;
c_6 \big(L_\beta \delta\big)^{\frac{p-1}{1+1/\beta -1/s}}\leq c_6~\big(L_\beta \delta\big)^{p\nu}.
\]
\par
2$^0$. {\em Dense zone: $\frac{2+1/\beta}{1+1/s} < p< s(2+\frac{1}{\beta})$.}
Because $p>\frac{2+1/\beta}{1+1/s}$,
by Proposition~\ref{prop:first-bound}, inequality (\ref{eq001:m-large}) with $\theta=1$,
\begin{equation}\label{eq:J-m-small}
 \sum_{m=m_0(1)}^0 J_m \;\leq\; c_{7} \varphi^p \sum_{m=m_0(1)}^0 2^{m(p-\frac{2+1/\beta}{1+1/s})}
\;\leq\; c_{8}\,\varphi^p
= c_{8} (L_\beta\delta)^{\frac{p\beta}{2\beta+1}}.
\end{equation}
Furthermore, because $p<s(2+\frac{1}{\beta})$ we have by Proposition~\ref{prop:first-bound}, inequality (\ref{eq02:m}),
that
\[
 \sum_{m=1}^{m_2} J_m \;\leq\; c_{9}
\varphi^p \sum_{m=1}^{m_2} 2^{m(p-s(2+\frac{1}{\beta}))} = c_{10}
(L_\beta\delta)^{\frac{p\beta}{2\beta+1}}.
\]
Thus, in the dense zone
\[
 \sum_{m=m_0(1)}^{m_2} J_m \leq c_{11} (L_\beta\delta)^{\frac{p\beta}{2\beta+1}}\leq c_{11}(L_\beta\delta)^{p\nu} .
\]
\par
3$^0$. {\em Sparse zone: $p>s(2+1/\beta)$, $s<1$.}
First we note that the bound in (\ref{eq:J-m-small}) remains true since
$p>s(2+1/\beta)$. By the same reason in view of Proposition~\ref{prop:first-bound}, inequality (\ref{eq02:m}),
\begin{equation}\label{eq:new-proof}
\sum_{m=1}^{m_2} J_m\leq c_{12} \varphi^p  2^{m_2(p-s(2+\frac{1}{\beta}))}
\leq c_{13} \varphi^{s(2+\frac{1}{\beta})}
= c_{13} \big(L_\beta \delta\big)^{s}\leq c_{13}(L_\beta\delta)^{p\nu}.
\end{equation}
Here we have used the definition of $m_2$. It remains to note that conditions $p>s(2+1/\beta)$, $s<1$
imply that $\varphi^{p}\delta^{-s}\to 0$ as $n\to 0$. Therefore the statement of the theorem
follows from (\ref{eq:J-m-small}) and (\ref{eq:new-proof}).
\par
4$^0$. {\em Sparse zone: $p>s(2+1/\beta)$,  $s\geq 1$.}
We need to bound only $\sum_{m=1}^{m_2} J_m$,
because (\ref{eq:J-m-small}) remains true.
By inequality (\ref{eq02:m}) of
Proposition~\ref{prop:first-bound} and because $p>s(2+1/\beta)$
\begin{eqnarray*}
%\label{eq0001:proof-th}
&&\sum_{m=1}^{m_1} J_m\leq c_{14}\varphi^p2^{m_1(p-s(2+\frac{1}{\beta}))}.
\end{eqnarray*}
Next, we have in view of the inequality (\ref{eq03:m}) of Proposition ~\ref{prop:first-bound}
$$
\sum_{m=m_1+1}^{m_2} J_m\leq c_{15}
\varphi^{p}\bigg[\frac{L_\gamma\varphi^{1/\beta}}{L_\beta\varphi^{1/\gamma}}\bigg]^{\upsilon}
\sum_{m=m_1+1}^{m_2}2^{m\left[p-\upsilon(2+1/\gamma)\right]}.
$$
Since $p-\upsilon(2+1/\gamma)<0$ [see (\ref{eq003:proof-th})],
$$
\sum_{m=m_1+1}^{m_2} J_m\leq c_{16}
\varphi^{p}\bigg[\frac{L_\gamma\varphi^{1/\beta}}{L_\beta\varphi^{1/\gamma}}\bigg]^{\upsilon}
2^{m_1\left(p-\upsilon[2+1/\gamma]\right)}
\leq c_{16}\varphi^{p}2^{m_1\left(p-s[2+1/\beta]\right)}.
$$
In order to obtain the second inequality we have used (\ref{eq:m-1-bounds}). Thus,
$$
\sum_{m=1}^{m_2} J_m
\leq c_{17}\varphi^{p}2^{m_1\left[p-s(2+1/\beta)\right]}.
$$
Using equality (\ref{eq:m-1-deff}) and  (\ref{eq:m-1-bounds}) we obtain
\begin{eqnarray*}
%\label{eq000100:proof-th}
\sum_{m=1}^{m_1} J_m\leq c_{20}
\Big(L_\gamma/L_\beta\Big)^{\frac{p-s(2+1/\beta)}{s(2+1/\beta)(1/s-1/\upsilon)+(1/\gamma-1/\beta)}}
(L_\beta\delta)^{\frac{p(1/s-1/\upsilon)+1/\gamma-1/\beta}{(2+1/\beta)(1/s-1/\upsilon)+(1/\gamma-1/\beta)s^{-1}}}.
\end{eqnarray*}
\par
The statement of the theorem is now obtained  by the following routine computations.
Denote
$$
A=\frac{1}{s_-}-\frac{1}{p\beta_-},\qquad \frac{1}{s_-}=\sum_{j\in I_-}\frac{1}{\beta_jr_j},\qquad  \frac{1}{\beta_-}=\sum_{j\in I_-}\frac{1}{\beta_j},\qquad \frac{1}{\gamma_-}=\sum_{j\in I_-}\frac{1}{\gamma_j}.
$$
First, we remark that
\begin{eqnarray}
\label{eq000101:proof-th}
&&p\Big(\frac{1}{s}-\frac{1}{\upsilon}\Big)+
\frac{1}{\gamma}-\frac{1}{\beta}=\frac{p}{s_-} - \frac{1}{\gamma_-} +
\frac{1}{\gamma_-} - \frac{1}{\beta_-}= \frac{p}{s_-} - \frac{1}{\beta_-}=Ap.
\end{eqnarray}
Next,
\begin{eqnarray*}
\frac{1}{\gamma_-}=\sum_{j\in I_-}\frac{\tau_j}{\tau(p)\beta_j}
&=& \frac{1}{\tau(p)}
\sum_{j\in I_-}\frac{1}{\beta_j}[1-1/s+1/(r_j\beta)]=
\frac{1-1/s}{\tau(p)\beta_-}+\frac{1}{\tau(p)\beta s_-}
\\
&=&\frac{1-1/s}{\tau(p)\beta_-}+\frac{1}{\tau(p)\beta}\bigg(\frac{1}{s_-}-\frac{1}{p\beta_-}\bigg)+\frac{1}{\tau(p)\beta p\beta_-}
\\
&=& \frac{1}{\tau(p)\beta_-}\bigg(1-\frac{1}{s}+\frac{1}{p\beta}\bigg)+\frac{A}{\tau(p)\beta}
=
\frac{1}{\beta_-}+\frac{A}{\tau(p)\beta}.
\end{eqnarray*}
Hence,
$1/\gamma-1/\beta=1/\gamma_--1/\beta_-=A/(\tau(p)\beta)$, which implies that
$$
\frac{1}{s}-\frac{1}{\upsilon}=\frac{1}{s_-}-\frac{1}{p\gamma_-}
=A+\frac{1}{p}\bigg(\frac{1}{\beta_-}-\frac{1}{\gamma_-}\bigg)=
A\bigg(1-\frac{1}{p\tau(p)\beta}\bigg).
$$
Two last equalities yield
\begin{eqnarray*}
\Big(2+\frac{1}{\beta}\Big)
\Big(\frac{1}{s}-\frac{1}{\upsilon}\Big)+
\Big(\frac{1}{\gamma}-\frac{1}{\beta}\Big)\frac{1}{s}
=\frac{A}{\tau(p)}\Big[\Big(2+\frac{1}{\beta}\Big)\Big(\tau(p)-\frac{1}{p\beta}\Big)+\frac{1}{s\beta}\Big]
=\frac{A}{\tau(p)}\Big(2+\frac{1}{\beta}-\frac{2}{s\beta}\Big),
\end{eqnarray*}
where the last equality follows from the fact that
$\tau(p)-1/(p\beta)=1-1/s$.
% we have
%that the last expression equals
%$
%(2+1/\beta)(1/s-1/\upsilon)+(1/\gamma-1/\beta)s^{-1}=(A/\tau)\big[2+1/\beta-2/s\beta\big]
%$
 This together with (\ref{eq000101:proof-th}).
leads to  the statement of the theorem
in the sparse zone.
\par
%\smallskip
%
5$^0$. {\em Boundary zones: $p=s(2+\frac{1}{\beta})$, $p=\frac{2+1/\beta}{1+1/s}$.}
Here the proof coincides with the proof for the dense zone with the only difference that
the corresponding sums equal $|m_1|$ and $m_2$ respectively.
% This results in extra
%$\ln (1/\delta)$ factor in the final bounds.
\epr
\subsection{Proof of statement (i) of Theorem~\ref{th:upper-new}}
\label{sec:upper-4}
In view of (\ref{eq:start}) and
by bound (\ref{eq:prop-new_theta}) of Proposition~\ref{prop:new},
\begin{eqnarray*}
 \bE_f\|\hat{f}-f\|_p^p
&\leq&  c_1\big(L_\beta \delta\big)^{p\nu(\theta)} +
c_2\sum_{m=m_0(\theta)}^\infty J_m.
\end{eqnarray*}
If $p<\frac{2+1/\beta}{1/\theta+1/s}$ then,
using bounds  (\ref{eq01:m-large})  and (\ref{eq02:m}) of
Proposition~\ref{prop:first-bound}, we have
 \begin{eqnarray*}
 \sum_{m=m_0(\theta)}^\infty J_m \;\leq\; c_3\varphi^p
 \sum_{m=m_0(\theta)}^{m_2} 2^{m\big(p-\frac{2+1/\beta}{1/\theta+1/s}\big)}
\;\leq\; c_4\,\varphi^p 2^{m_0(\theta)\big(p-\frac{2+1/\beta}{1/\theta+1/s}\big)}=
c_5(L_\beta\delta)^{\frac{p-\theta}{1-\theta/s+1/\beta}},
\end{eqnarray*}
and the assertion of the theorem
follows.
%in the case $p<\frac{2+1/\beta}{1/\theta+1/s}$.
If $s(2+1/\beta)\geq p\geq\frac{2+1/\beta}{1/\theta+1/s}$ then
$$
\sum_{m=m_0(\theta)}^\infty J_m \;\leq\; c_6\mu_n^{p}(\theta)\varphi^p\leq c_7\mu_n^{p}(\theta)
\big(L_\beta \delta\big)^{p\nu(\theta)}.
$$
\section{Proofs of Theorem~\ref{th:lower-bound-in-L_p},
statement~(ii) of Theorem~\ref{th:upper-new} and  the lower bound in (\ref{eq:lower-bound-new})}
\label{sec:Proof-th-2}
The proof is organized as follows. First, we formulate
two auxiliary statements,
Lemmas~\ref{lem:tsyb_book-result} and~\ref{lem:V-G-lemma}.
Second, we present a general construction
of a finite set of functions employed in the proof of lower
bounds. Then we specialize the constructed set of functions
in different regimes and derive the announced lower bounds.
\subsection{Auxiliary lemmas}
The first statement given
in Lemma~\ref{lem:tsyb_book-result} is a
simple consequence of Theorem~2.4 from \cite{Tsybakov}. Let $\bF$ be a given
set of probability densities.
\begin{lemma}
\label{lem:tsyb_book-result}
Assume that for any sufficiently large
integer $n$ one can find a positive real number $\rho_n$
and  a finite subset of functions
$\big\{f^{(0)}, f^{(j)},\;j\in\cJ_n\big\}\subset \bF
$
such that
\begin{eqnarray}
\label{eq:ass1-klp-lemma}
  &&\big\|f^{(i)}- f^{(j)}\big\|_p \geq 2\rho_n,\qquad\;
\forall  i, j\in \cJ_n\cup\{0\}:\;i\neq j;
\\*[2mm]
\label{eq:ass2-klp-lemma}
&& \limsup_{n\to\infty}\frac{1}{|\cJ_n|^{2}}
\sum_{j\in\cJ_n}\bE_{f^{(0)}}\Bigg\{
\frac{\rd \bP_{f^{(j)}}}{\rd \bP_{f^{(0)}}}(X^{(n)})\Bigg\}^{2}=:C <\infty.
\end{eqnarray}
Then for any $q\geq 1$
$$
\liminf_{n\to\infty}
\inf_{\tilde f}\;
\sup_{f \in \bF}
\rho^{-1}_n\left(\bE_f \big\|\tilde{f} - f\big\|^{q}_p\right)^{1/q}
\geq \left(\sqrt{C} +\sqrt{C+1} \right)^{-2/q},
$$
where infimum on the left hand side is taken over all possible estimators.
\end{lemma}
We will apply Lemma \ref{lem:tsyb_book-result} with $\bF= \bN_{\vec{r},d}(\vec{\beta},\vec{L}, M)$ in the proof of Theorem \ref{th:lower-bound-in-L_p}
and with $\bF=\bG_\theta(R)\cap\bN_{\vec{r},d}(\vec{\beta},\vec{L}, M)$ in the proof of statement~(ii) of Theorem~\ref{th:upper-new}.
\par
Next we quote
the Varshamov--Gilbert lemma
[see, e.g., Lemma~2.9 in \cite{Tsybakov}].
\begin{lemma}[Varshamov--Gilbert]
\label{lem:V-G-lemma}
Let $\varrho_m$ be the  Hamming distance on $\{0,1\}^m$, $m\in\bN^*$, i.e.
$$
\varrho_m(a,b)=\sum_{j=1}^m  {\bf 1}\left\{a_j\neq b_j\right\},\quad a,b\in\{0,1\}^m.
$$
For any $m\geq 8$ there exists a subset $\cP_m$ of $\{0,1\}^m$ such that
$\big|\cP_m\big|\geq 2^{m/8}$, and
$$
\varrho_m\big(a,a^\prime\big)\geq \frac{m}{8},\;\;\;\;\;\forall a,a^\prime\in\cP_m.
$$
\end{lemma}

\subsection{Proof of Theorem \ref{th:lower-bound-in-L_p}. General construction of a finite set of functions}
%Proof of Theorem~\ref{th:lower-bound-in-L_p}}
\par
1$^0$. For any $t\in\bR$ set
$$
\Lambda(t)=
\bigg(\int_{-1}^1 e^{-1/(1-u^2)}\rd u\bigg)^{-1}
e^{-1/(1-t^2)}\; \mathbf{1}_{[-1,1]}(t).
$$
Note that
$\Lambda$ is a probability density
compactly supported on $[-1,1]$ and infinitely differentiable on the real line,
$\Lambda\in \bC^{\infty}(\bR^1)$.
Obviously, for any
$\alpha>0$ and $r\geq 1$ there exists constant $c_1=c_1(\alpha,r)<\infty$ such that
\begin{equation}
\label{eq1:proof-th:lower-bound-in-L_p}
\Lambda\in\bN_{r,1}(\alpha,c_1).
\end{equation}
\par
 Define
%Set for any $x\in\bR^d$
$$
\bar{f}^{(0)}(x)=
\prod_{l=1}^d\bigg[ \frac{1}{N}
\int_{\bR^1}\Lambda(y-x_l)\mathbf{1}_{[-\frac{N}{2},\frac{N}{2}]}(y)\rd y\bigg],\;\;\;\;\;
x=(x_1,\ldots,x_d)\in \bR^d,
$$
where  parameter $N=N(n)>8$ will be chosen later.
By construction,
$\bar{f}^{(0)}$ is a probability density for any choice of
$N$, ${\rm supp}(\bar{f}^{(0)})=[-N/2-1, N/2+1]^d$, and
\begin{equation}\label{eq:bar-f-0}
 \bar{f}^{(0)}(x)=N^{-d},\;\;\;\forall x \in \big[-N/2+1, N/2-1\big]^d.
\end{equation}
Moreover,
in view of (\ref{eq1:proof-th:lower-bound-in-L_p}) and by the Young inequality,
there exist constants $\vec{C}=\big(\tilde{C}_1,\ldots,\tilde{C}_d\big)$
depending on $\vec{\beta}$ and $\vec{r}$ only
such that
\begin{equation}
\label{eq2:proof-th:lower-bound-in-L_p}
\bar{f}^{(0)}\in\bN_{\vec{r},d}\big(\vec{\beta},\vec{C}\big).
\end{equation}
Note that
$\vec{C}$ do not depend on $N$.
\par
Let $L_0>0$ be fixed, and let $f^{(0)}(x)=\kappa^{d}\bar{f}^{(0)}\big(x\kappa\big)$,
where $\kappa>0$
is chosen in such a way that $f^{(0)}$ belongs to the class
$\bN_{\vec{r},d}(\vec{\beta},2^{-1}\vec{L}_0)$, where $\vec{L}_0=(L_0,\ldots,L_0)$.
The existence of such
$\kappa$ independent of $N$ and determined by
$\vec{\beta}$, $\vec{r}$  and $L_0$ is
guaranteed by (\ref{eq2:proof-th:lower-bound-in-L_p}).
Note also that $f^{(0)}$ is a probability density. Moreover, we remark that
$\|\bar{f}^{(0)}\|_\infty\leq N^{-d}$ since $\int |\Lambda|=1$.
%in view of (\ref{eq:bar-f-0}).
Thus,
\begin{equation}
\label{eq2-new:proof-th:lower-bound-in-L_p}
f^{(0)}\in \bN_{\vec{r},d}(\vec{\beta}, \vec{L}_0/2, M/2),
\end{equation}
provided that $N>(2M^{-1})^{1/d}\kappa$. This condition is assumed to be fulfilled.
\par
\medskip
2$^0$. Put for any $t\in\bR^1$
$$
g(t)= \int_{\bR^1}\Lambda(y-t)\Big[\mathbf{1}_{[0,1]}(y)-\mathbf{1}_{[-1,0]}(y)\Big]\rd y.
$$
We obviously have $g\in\bC^\infty(\bR^1)$, and
\begin{eqnarray}
\label{eq3:proof-th:lower-bound-in-L_p}
&& ({\rm i})\;\;\int_{\bR^1} g(y)\rd y =0,\qquad ({\rm ii})\;\;\text{supp}(g)\subseteq [-2,2],\qquad
({\rm iii})\;\;\|g\|_\infty\leq 1.
\end{eqnarray}
\par
For any
$l=1,\ldots,d$ let
$(20\kappa)^{-1}>\sigma_l=\sigma_l(n)\to 0$, $n\to\infty$,  be the sequences
to be specified later. Let
$M_l=(20\kappa\sigma_l)^{-1}N$,
and without loss of generality assume that $M_l$, $l=1,\ldots, d$ are integer numbers.
Define  also
$$
x_{j,l}=-\frac{N-4}{4\kappa}+8j\sigma_l,\;\;\; j=1,\ldots, M_l,\;\;l=1,\ldots,d,
$$
and let $\cM=\{1,\ldots, M_1\}\times\cdots\times\{1,\ldots, M_d\}$.
For any $m=(m_1,\ldots,m_d)\in \cM$ define
\begin{eqnarray*}
&&G_m(x)=\prod_{l=1}^dg\left(\frac{x_l-x_{m_l,l}}{\sigma_l}\right),\quad x\in\bR^d,
\\
&&\Pi_m=\big[x_{m_1,1}-3\sigma_1,x_{m_1,1}+3\sigma_1\big]
\times\cdots\times\big[x_{m_d,d}-3\sigma_d,x_{m_d,d}-3\sigma_d\big]\subset\bR^d.
\end{eqnarray*}
Several remarks on these definitions are in order.
First, in view of (\ref{eq3:proof-th:lower-bound-in-L_p})(ii)
\begin{eqnarray}
\label{eq5:proof-th:lower-bound-in-L_p}
&&\text{supp}\big(G_m\big)\subset \Pi_m,\quad \forall m\in\cM,
\\*[2mm]
\label{eq6:proof-th:lower-bound-in-L_p}
&& \Pi_m\cap \Pi_j=\emptyset,\quad \forall m, j\in\cM:\;\;m\neq j.
%\\*[2mm]
\end{eqnarray}
Second, since $g\in\bC^{\infty}(\bR^1)$, we have that
$G_m\in\bC^{\infty}(\bR^d)$ for any $m\in\cM$.
Moreover, for any $l=1,\ldots,d$, any $|h|\leq\sigma_l$
and any integer $k$
\begin{eqnarray}
\label{eq7:proof-th:lower-bound-in-L_p}
&&\text{supp}\Big\{\Delta_{h,l}\, (D^{k}_lG_m) \Big\}\subseteq \Pi_m,
\quad \forall m\in\cM,
\end{eqnarray}
where $D^k_l G$ stands for the $k$th order derivative of a function $G$ with respect
to the variable $x_l$, and $\Delta_{h,l}$ is the first order difference operator with step size
$h$ in direction of the variable $x_l$.
\par
For $m\in \cM$ define
\begin{eqnarray*}
\pi(m)=\sum_{j=1}^{d-1}(m_j-1)\bigg(\prod_{l=j+1}^d M_l\bigg)+m_d.
\end{eqnarray*}
It is easily checked that $\pi$ defines enumeration of the set $\cM$, and
$\pi:\cM\to \{1,2\ldots,|\cM|\} $ is a bijection.
Let $W$ be a subset of $\{0,1\}^{|\cM|}$.
Define a family of functions $\{F_w, w\in W\}$ by
$$
%\cG_{w}
F_w(x)=A\sum_{m\in\cM}w_{\pi(m)}G_m(x),\;\;\;\;x\in \bR^d,
$$
where $w_j$, $j=1,\ldots, |\cM|$ are
the coordinates of $w$,  and $A$ is a parameter to be specified.
It follows from
(\ref{eq3:proof-th:lower-bound-in-L_p})(iii), (\ref{eq5:proof-th:lower-bound-in-L_p}) and
(\ref{eq6:proof-th:lower-bound-in-L_p}) that
\begin{equation}
\label{eq8-new:proof-th:lower-bound-in-L_p}
\big\|F_{w}\big\|_\infty\leq A,\quad  \forall w\in W,
\end{equation}
and (\ref{eq3:proof-th:lower-bound-in-L_p})(i) implies that
\begin{equation}
\label{eq9:proof-th:lower-bound-in-L_p}
\int_{\bR^d}F_{w}(x)\rd x=0,\quad  \forall w\in W.
\end{equation}
\par\medskip
3$^0$.
Now we find conditions
which guarantee that
$F_{w}\in \bN_{\vec{r},d}(\vec{\beta},2^{-1}\vec{L})$
for any
$w\in W$.
\par
Fix $l=1,\ldots,d$, and
let
$k_l=\lfloor\beta_l\rfloor+1$ if $\beta_l\notin\bN^*$, and
$k_l=\lfloor\beta_l\rfloor+2$ if $\beta_l\in\bN^*$ (here $\lfloor x \rfloor$ stands for the
maximal integer number strictly less than $x$).
\par
First, for any $w\in W$  and $h\in\bR$
\begin{equation}
\label{eq10:proof-th:lower-bound-in-L_p}
\Big\|\Delta^{k_l}_{h,l} F_w\Big\|_{r_l}=
\Big\|\Delta^{k_l-1}_{h,l} (\Delta_{h,l} F_w)\Big\|_{r_l}
\leq |h|^{k_l-1}\Big\| \Delta_{h,l} (D^{k_l-1}_lF_w)\Big\|_{r_l},
\end{equation}
where the last inequality is found in \cite[Section~4.4.4]{Nikolski}.
Next, in view of (\ref{eq6:proof-th:lower-bound-in-L_p}) and
(\ref{eq7:proof-th:lower-bound-in-L_p}) we obtain for any $w\in W$
  and any $r_l\neq\infty$
\begin{eqnarray*}
&& \Big\| \Delta_{h,l} (D^{k_l-1}_l F_w)\Big\|^{r_l}_{r_l}
\;=\;
\sum_{j\in\cM}\int_{\Pi_j}\Big| \Delta_{h,l} (D^{k_l-1}_lF_w)(x)\Big|^{r_l}\rd x
\nonumber
\\
&&\;\;\;
=\; A^{r_l} \sum_{j\in\cM}w_{\pi(j)}\int_{\Pi_j}\Big| \Delta_{h,l} (D^{k_l-1}_l G_j)(x)\Big|^{r_l}\rd x
\nonumber\\*[2mm]
&&\;\;\;\leq\; A^{r_l}S_W\big\|g\big\|^{(d-1)r_l}_{r_l}\sigma_l^{-(k_l-1)r_l}
\bigg(\prod_{j=1}^d\sigma_j\bigg)\Big\|g^{(k_l-1)}\Big(\cdot-\frac{h}{\sigma_l}\Big)
-g^{(k_l-1)}(\cdot)\Big\|^{r_l}_{r_l},
\end{eqnarray*}
where we have put
$S_W:=\sup_{w\in W}|\{j:\;w_j\neq 0\}|$. Thus, for any $r_l\neq\infty$ we have
\begin{eqnarray}
\label{eq11:proof-th:lower-bound-in-L_p}
\Big\| \Delta_{h,l} (D^{k_l-1}_l F_w)\Big\|_{r_l}\leq\; A \big\|g\big\|^{(d-1)r_l}_{r_l}\sigma_l^{-(k_l-1)}
\bigg(S_W\prod_{j=1}^d\sigma_j\bigg)^{\frac{1}{r_l}}\Big\|g^{(k_l-1)}\Big(\cdot-\frac{h}{\sigma_l}\Big)
-g^{(k_l-1)}(\cdot)\Big\|_{r_l}.
\end{eqnarray}
Similarly, we get for any $w\in W$
\begin{eqnarray}
\label{eq11-new:proof-th:lower-bound-in-L_p}
&& \Big\| \Delta_{h,l} (D^{k_l-1}_l F_w)\Big\|_{\infty}
\;=\;
\sup_{j\in\cM}\sup_{x\in\Pi_j}\Big| \Delta_{h,l} (D^{k_l-1}_lF_w)(x)\Big|
\nonumber
\\
&&\;\;\;
=\; A \sup_{j\in\cM}w_{\pi(j)}\sup_{x\in\Pi_j}\Big| \Delta_{h,l} (D^{k_l-1}_l G_j)(x)\Big|
\nonumber\\*[2mm]
&&\;\;\;\leq\; A \big\|g\big\|^{(d-1)}_{\infty}\sigma_l^{-(k_l-1)}
\Big\|g^{(k_l-1)}\Big(\cdot-\frac{h}{\sigma_l}\Big)
-g^{(k_l-1)}(\cdot)\Big\|_{\infty}.
\end{eqnarray}
\par
In view of (\ref{eq3:proof-th:lower-bound-in-L_p})(ii)
and $|h|\leq\sigma_l$,
function $g^{(k_l-1)}\big(\cdot-[h/\sigma_l]\big)
-g^{(k_l-1)}(\cdot)$ is supported on $[-3,3]$. Therefore
the fact that $g\in\bC^{\infty}(\bR^1)$ implies for any $r_l\in [1,\infty]$
$$
\Big\|g^{(k_l-1)}\big(\cdot-h/\sigma_l\big)
-g^{(k_l-1)}(\cdot)\Big\|_{r_l}\leq 6^{1/r_l}\big\|g^{(k_l)}\big\|_{\infty}(h/\sigma_l)\leq
6^{1/r_l}\big\|g^{(k_l)}\big\|_{\infty}|h/\sigma_l|^{\beta_l-k_l+1}.
$$
In the last inequality
we have used that $0\leq \beta_l-k_l+1\leq 1$ by definition of $k_l$.
Combining this with
(\ref{eq10:proof-th:lower-bound-in-L_p}),  (\ref{eq11:proof-th:lower-bound-in-L_p}) and (\ref{eq11-new:proof-th:lower-bound-in-L_p})
we have for any $|h|\leq\sigma_l$ and any $r_l\in [1,\infty]$
\begin{eqnarray}
\label{eq12:proof-th:lower-bound-in-L_p}
&&\Big\|\Delta^{k_l}_{h,l} F_w\Big\|_{r_l}
\leq A|h|^{\beta_l}6^{1/r_l}\big\|g\big\|_{r_l}^{d-1}\big\|g^{(k_l)}\big\|_{\infty}
\sigma_l^{-\beta_l}\bigg(S_W\prod_{j=1}^d\sigma_j\bigg)^{1/r_l}.
\end{eqnarray}
\par
If $|h|\geq\sigma_l$ then we note that
$\Delta_{h,l}(D^{k_l-1}_l F_w)(\cdot)=(D^{k_l-1}_lF_w)(\cdot-h e_l)-(D^{k_l-1}_lF_w)(\cdot)$,
and
by the triangle inequality
$$
\Big\| \Delta_{h,l}(D^{k_l-1}_l F_w)\Big\|_{r_l}\leq
2 \Big\|D^{k_l-1}_l F_w\Big\|_{r_l}\leq 2 \Big\|D^{k_l-1}_lF_w
\Big\|_{r_l}|h/\sigma_l|^{\beta_l-k_l+1}.
$$
In view of (\ref{eq5:proof-th:lower-bound-in-L_p}) and
(\ref{eq6:proof-th:lower-bound-in-L_p}) we get for any $w\in W$ and any $r_l\neq\infty$
\begin{eqnarray*}
\Big\|D^{k_l-1}_l F_w\Big\|^{r_l}_{r_l}&=&
\sum_{j\in\cM}\int_{\Pi_j}\Big|D^{k_l-1}_l F_w(x)\Big|^{r_l}\rd x
=A^{r_l}\sum_{j\in\cM}w_{\pi(j)}\int_{\Pi_j}\Big|D^{k_l-1}_lG_j(x)\Big|^{r_l}\rd x
\nonumber\\*[2mm]
&\leq& A^{r_l}S_W\big\|g\big\|^{(d-1)r_l}_{r_l}\Big\|g^{(k_l-1)}\Big\|^{r_l}_{r_l}\sigma_l^{(1-k_l)r_l}
\bigg(\prod_{j=1}^d\sigma_j\bigg).
\end{eqnarray*}
Moreover,
\begin{eqnarray*}
\Big\|D^{k_l-1}_l F_w\Big\|_{\infty}&=&
\sup_{j\in\cM}\sup_{x\in\Pi_j}\Big|D^{k_l-1}_l F_w(x)\Big|
=A\sup_{j\in\cM}w_{\pi(j)}\sup_{x\in\Pi_j}\Big|D^{k_l-1}_lG_j(x)\Big|
\nonumber\\*[2mm]
&\leq& A\big\|g\big\|^{(d-1)}_{\infty}\Big\|g^{(k_l-1)}\Big\|_{\infty}\sigma_l^{(1-k_l)}.
\end{eqnarray*}
We obtain finally from (\ref{eq10:proof-th:lower-bound-in-L_p}) that
for any  $|h|\geq\sigma_l$ and any $r_l\in [1,\infty]$
\begin{eqnarray}
\label{eq13:proof-th:lower-bound-in-L_p}
&&\Big\|\Delta^{k_l}_{h,l} F_w\Big\|_{r_l}
\leq A|h|^{\beta_l}2\big\|g\big\|^{d-1}_{r_l}\big\|g^{(k_l-1)}\big\|_{r_l}
\sigma_l^{-\beta_l}\bigg(S_W\prod_{j=1}^d\sigma_j\bigg)^{1/r_l}.
\end{eqnarray}
Combining
(\ref{eq12:proof-th:lower-bound-in-L_p}) and (\ref{eq13:proof-th:lower-bound-in-L_p})
we conclude that for any $w\in W$ and $r_l\in [1,\infty]$
$$
\Big\|\Delta^{k_l}_{h,l} F_w\Big\|_{r_l}\leq C_1 A|h|^{\beta_l}
\sigma_l^{-\beta_l}\bigg(S_W\prod_{j=1}^d\sigma_j\bigg)^{1/r_l},\quad\forall h\in\bR^1,
$$
where
$C_1=\max_{l}(\|g\big\|_{r_l}^{d-1}
\max\{6^{1/r_l}\|g^{(k_l)}\|_{\infty},2\|g^{(k_l-1)}\|_{r_l}\})$.
Thus, if
\begin{eqnarray}
\label{eq14:proof-th:lower-bound-in-L_p}
&&A\sigma_l^{-\beta_l}
\bigg(S_W\prod_{j=1}^d\sigma_j\bigg)^{1/r_l}\leq (2C_1)^{-1} L_l,\quad\forall l=1,\ldots, d
\end{eqnarray}
then $F_{w}\in\bN_{\vec{r},d}(\vec{\beta},2^{-1}\vec{L})$ for any $w\in W$.
\par\medskip
4$^0.$ Define for any $w\in W$
$$
f_w(x)=f^{(0)}(x)+F_{w}(x),\;\;\;\;\; x\in\bR^d.
$$
Remind that $f^{(0)}$
is the  probability density belonging to $\bN_{\vec{r},d}(\vec{\beta},\vec{L}_0/2, M/2)$.
Therefore, in view of
(\ref{eq9:proof-th:lower-bound-in-L_p}) and under condition (\ref{eq14:proof-th:lower-bound-in-L_p}),
for any $w\in W$
\begin{eqnarray}
\label{eq15:proof-th:lower-bound-in-L_p}
&& \int_{\bR^d} f_w(x) \rd x=1,\;\;\;\; f_w\in \bN_{\vec{r},d}(\vec{\beta},\vec{L}),
\end{eqnarray}
where the latter inclusion holds
because  $\min_{j=1,\ldots,d}L_j\geq L_0$.
\par
By construction of $F_{w}$, for any $w\in W$
\begin{eqnarray}
\label{eq16:proof-th:lower-bound-in-L_p}
&& F_{w}(x)=0, \quad \;\;\;\;\forall x\notin
 \Big[-\frac{1}{4\kappa}(N-4), \;\frac{1}{4\kappa}(N+4)\Big]^d.
\end{eqnarray}
This yields
\begin{eqnarray}
\label{eq1600:proof-th:lower-bound-in-L_p}
&&f_w(x)=f^{(0)}(x)\geq 0, \quad \forall x\notin \Big[-\frac{1}{4\kappa}(N-4),
\;\frac{1}{4\kappa}(N+4)\Big]^d.
\end{eqnarray}
On the other hand, by (\ref{eq:bar-f-0})
\begin{eqnarray}
\label{eq160:proof-th:lower-bound-in-L_p}
f^{(0)}(x)=\kappa^{d} N^{-d},\quad\;\;\; \forall x\in
\Big[-\frac{1}{4\kappa}(N-4), \;\frac{1}{4\kappa}(N+4)\Big]^d.
\end{eqnarray}
Therefore, if we require
\begin{eqnarray}
\label{eq17:proof-th:lower-bound-in-L_p}
&& A\leq \kappa^d N^{-d},
\end{eqnarray}
this
together with (\ref{eq10:proof-th:lower-bound-in-L_p}) implies
$$
f_w(x)\geq 0,\;\;\;\quad \forall x\in\Big[-\frac{1}{4\kappa}(N-4), \;\frac{1}{4\kappa}(N+4)\Big]^d.
$$
We conclude that $f_w\geq 0$ for any $w\in W$.
Moreover, we get from  (\ref{eq2-new:proof-th:lower-bound-in-L_p}),
(\ref{eq8-new:proof-th:lower-bound-in-L_p}) and
(\ref{eq17:proof-th:lower-bound-in-L_p}) that $\|f_w\|_\infty \leq M$ for any $w\in W$.
\par
All this, together
with (\ref{eq15:proof-th:lower-bound-in-L_p}), shows that
$\{f^{(0)}, f_w, w\in W\}$ is a finite set of probability densities
from
$\bN_{\vec{r},d}(\vec{\beta},\vec{L}, M)$.
Thus Lemma~\ref{lem:tsyb_book-result} is applicable with $\cJ_n=W$ and $\bF= \bN_{\vec{r},d}(\vec{\beta},\vec{L}, M)$.
\par
\medskip
5$^0$. Suppose now that the set $W$ is chosen so that
\begin{eqnarray}
\label{eq18:proof-th:lower-bound-in-L_p}
&& \varrho_{|\cM|}\big(w,w^\prime\big)\geq B,\quad\forall w,w^\prime\in W,
\end{eqnarray}
where, we remind, $\varrho_{|\cM|}$ is the Hamming distance on $\{0,1\}^{|\cM|}$.
Here $B=B(n)\geq 1$ is a parameter to be specified.
Then
we deduce from (\ref{eq15:proof-th:lower-bound-in-L_p}),
(\ref{eq5:proof-th:lower-bound-in-L_p}) and (\ref{eq6:proof-th:lower-bound-in-L_p}), that for all $w,w^\prime\in W$
\begin{eqnarray}
\label{eq19:proof-th:lower-bound-in-L_p}
\left\|f_w-f_{w^\prime}\right\|^{p}_p&=& \left\|F_w-F_{w^\prime}\right\|^{p}_p
=A^p\sum_{j\in\cM}\left|w_{\pi(j)}-w^\prime_{\pi(j)}\right|\int_{\Pi_j}\Big|G_j(x)\Big|^{p}\rd x
\nonumber\\
&=&A^{p}\big\|g\big\|^{dp}_{p}
\bigg(\prod_{j=1}^d\sigma_j\bigg)\sum_{j\in\cM}\left|w_{\pi(j)}-w^\prime_{\pi(j)}\right|=
A^{p}\big\|g\big\|^{dp}_{p}
\bigg(\prod_{j=1}^d\sigma_j\bigg)\varrho_{|\cM|}\big(w,w^\prime\big)
\nonumber\\
&\geq& \big\|g\big\|^{dp}_{p}A^pB\bigg(\prod_{j=1}^d\sigma_j\bigg).
%\nonumber
\end{eqnarray}
Here we have used that the map $\pi$ is a bijection. Putting
$C_2=\frac{1}{2}\big\|g\big\|^{d}_{p}$,
 we conclude that condition (\ref{eq:ass1-klp-lemma})
of Lemma~\ref{lem:tsyb_book-result} is fulfilled with
\begin{eqnarray}
\label{eq20:proof-th:lower-bound-in-L_p}
&& \rho_n=C_2A\bigg(B\prod_{j=1}^d\sigma_j\bigg)^{1/p}.
\end{eqnarray}
Let us remark that (\ref{eq20:proof-th:lower-bound-in-L_p}) remains true if we formally put $p=\infty$. Indeed, similarly to
(\ref{eq19:proof-th:lower-bound-in-L_p}),
\begin{eqnarray}
\label{eq190:proof-th:lower-bound-in-L_p}
&&\left\|f_w-f_{w^\prime}\right\|_\infty= \left\|F_w-F_{w^\prime}\right\|_\infty
=A\sup_{j\in\cM}\big|w_{j}-w^\prime_{j}\big|\big\|g\big\|^{d}_{\infty}\geq A\big\|g\big\|^{d}_{\infty}.
\end{eqnarray}
Here we have used (\ref{eq6:proof-th:lower-bound-in-L_p}),
the fact that the map $\pi$ is a bijection and, that
$w\neq w^\prime$ for all $w,w^\prime\in W$ in view of (\ref{eq18:proof-th:lower-bound-in-L_p}).
\par
Now we verify condition (\ref{eq:ass2-klp-lemma}) of Lemma~\ref{lem:tsyb_book-result}.
First observe that
$$
\frac{\rd \bP_{f_w}}{\rd \bP_{f^{(0)}}}\big(X^{(n)}\big)=
\prod_{k=1}^n\frac{f_w(X_k)}{f^{(0)}(X_k)}.
$$
Since  $X_k$, $k=1,\ldots,n$ are i.i.d. random vectors, we have for any $w\in W$
\begin{eqnarray*}
%\label{eq210:proof-th:lower-bound-in-L_p}
\bE_{f^{(0)}}\Bigg\{
\prod_{k=1}^n\frac{f_w(X_k)}{f^{(0)}(X_k)}
\Bigg\}^{2}&=&\Bigg\{\int_{\bR^d}
\frac{ f^{(0)}(x)+2F_w(x)+F^{2}_w(x)}{f^{(0)}(x)}\rd x
\Bigg\}^{n}
\nonumber\\
&=&\Bigg\{1+\int_{\bR^d}\frac{F^{2}_w(x)}{f^{(0)}(x)}\rd x
\Bigg\}^{n}.
\end{eqnarray*}
The last equality follows from (\ref{eq9:proof-th:lower-bound-in-L_p}).
By
(\ref{eq16:proof-th:lower-bound-in-L_p}) and (\ref{eq160:proof-th:lower-bound-in-L_p}),
$$
\int_{\bR^d}\frac{F^{2}_w(x)}{f^{(0)}(x)}\rd x=
\kappa^{-d}N^d \|F_w\|_2^2;
$$
hence for any $w\in W$
\begin{equation*}
%\label{eq21:proof-th:lower-bound-in-L_p}
\bE_{f^{(0)}}\Bigg\{
\frac{\rd \bP_{f_w}}{\rd \bP_{f^{(0)}}}
\big(X^{(n)}\big)
\Bigg\}^{2}=\Big\{1+\kappa^{-d}N^d \|F_w\big\|_2^2
\Big\}^{n}\leq\exp{\Big\{n\kappa^{-d}N^d \|F_w\big\|_2^2
\Big\}}.
\end{equation*}
Repeating computations that led to (\ref{eq13:proof-th:lower-bound-in-L_p}) we have
$$
\big\|F_w\big\|_2^2\leq A^{2}\|g\|^{2d}_{2} S_W\prod_{j=1}^d\sigma_j.
$$
The right hand side of the latter inequality does not depend on $w$; hence we
$$
\frac{1}{|W|^{2}}\sum_{w\in W}\bE_{f^{(0)}}
\Bigg\{
\frac{\rd \bP_{f_w}}{\rd \bP_{f^{(0)}}}\big(X^{(n)}\big)
\Bigg\}^{2}\leq
\exp\Big\{C_3n A^2S_W N^d\bigg(\prod_{j=1}^d\sigma_j\bigg)-\ln{\big(|W|\big)}
\Big\},
$$
where we have put $C_3=\kappa^{-d}\|g\big\|^{2d}_{2}$.
Therefore, if
\begin{eqnarray}
\label{eq22:proof-th:lower-bound-in-L_p}
C_4nA^2S_WN^d \prod_{j=1}^d \sigma_j\leq\ln{\big(|W|\big)}
\end{eqnarray}
then condition
(\ref{eq:ass2-klp-lemma}) of Lemma \ref{lem:tsyb_book-result} is fulfilled with $C=1$.
\par
In order to apply Lemma~\ref{lem:tsyb_book-result}
it remains to specify the set $W$ and the parameters $A$,
$N$, $\sigma_j$, $j=1,\ldots,d$ so that
the relationships
(\ref{eq14:proof-th:lower-bound-in-L_p}), (\ref{eq17:proof-th:lower-bound-in-L_p}),
(\ref{eq18:proof-th:lower-bound-in-L_p}),
 and (\ref{eq22:proof-th:lower-bound-in-L_p}) are simultaneously fulfilled.
According to (\ref{lem:tsyb_book-result}),
 under these conditions the lower bound is given by $\rho_n$ in
(\ref{eq20:proof-th:lower-bound-in-L_p}).

\subsection{Proof of Theorem \ref{th:lower-bound-in-L_p}. Derivation of lower bounds in different zones}

We begin with the construction of the set $W$.
Let $m\geq 8$ be an integer number whose choice will be made later,
and, without loss of generality, assume that $|\cM|/m$ is integer.
Let $\cP_m$ be a subset of $\{0,1\}^{m}$
such that
\begin{eqnarray}
\label{eq24-new:proof-th:lower-bound-in-L_p}
 |\cP_m|\geq 2^{m/8},\;\;\;\varrho_{m}(z,z^\prime)\geq m/8,\;\;\forall z,z^\prime \in \cP_m.
\end{eqnarray}
Existence of such set $\cP_m$ is guaranteed by
Lemma~\ref{lem:V-G-lemma}.
Let $\cJ:=\{1+\frac{j}{m}|\cM|, \;j=0,\ldots, m-1\}$,
and note that $\cJ\subseteq\{1,\ldots,|\cM|\}$
with the equality in the case $m=|\cM|$.
Define the map $\Upsilon: \cP_m\to \{0,1\}^{|\cM|}$ by
$$
\Upsilon_j[a]=
\left\{
\begin{array}{ll}
a_j,\quad & j\in\cJ,
\\
0,\quad & j\in\{1,\ldots,|\cM|\}\setminus\cJ,
\end{array}
\right.
$$
and let  $W=\Upsilon(\cP_m)$.
Obviously,
$\varrho_{|\cM|}(w,w^\prime)=\varrho_{|\cM|}(\Upsilon[a],\Upsilon[a^\prime])=
\varrho_m(a,a^\prime)$
for all $w,w^\prime\in W$; therefore
(\ref{eq24-new:proof-th:lower-bound-in-L_p}) implies
that
\begin{eqnarray}
\label{eq25-new:proof-th:lower-bound-in-L_p}
 |W|\geq 2^{m/8},\;\;\;\varrho_{|\cM|}(w,w^\prime)\geq 8^{-1}m,\;\;\forall w,w^\prime \in W.
\end{eqnarray}
With such a set  $W$, $S_W\leq m$; moreover, since $\ln(|W|)\geq m\ln 2/8$,
condition (\ref{eq22:proof-th:lower-bound-in-L_p}) holds true if
\begin{eqnarray}
\label{eq24:proof-th:lower-bound-in-L_p}
&&A^{2}n N^d \prod_{j=1}^d \sigma_j \leq (8C_4)^{-1}\ln 2.
\end{eqnarray}
  We also note that condition (\ref{eq14:proof-th:lower-bound-in-L_p}) is fulfilled if we require
\begin{eqnarray}
\label{eq26-new:proof-th:lower-bound-in-L_p}
A\sigma_l^{-\beta_l} \bigg(m\prod_{j=1}^d \sigma_j\bigg)^{1/r_l}\leq (2C_1)^{-1} L_l,\quad\forall
l=1,\ldots, d.
\end{eqnarray}
In addition,
(\ref{eq18:proof-th:lower-bound-in-L_p})
holds with $B=m/8$.
\subsubsection{Tail zone: $p\leq \frac{2+1/\beta}{1+1/s}$}\label{subsec:tail}
Let $m=|\cM|$. By construction,
 $|\cM|=\prod_{l=1}^d M_l= (20\kappa)^{-d}N^d\prod_{l=1}^d\sigma_l^{-1}$ and, therefore
(\ref{eq26-new:proof-th:lower-bound-in-L_p}) is reduced to
\begin{eqnarray}
\label{eq2700-new:proof-th:lower-bound-in-L_p}
A\sigma_l^{-\beta_l}N^{d/r_l}\leq C_5L_l.
\end{eqnarray}
Thus, choosing
\begin{eqnarray}
\label{eq27-new:proof-th:lower-bound-in-L_p}
\sigma_l=C_6 A^{1/\beta_l}
L_l^{- 1/\beta_l} N^{\frac{d}{\beta_l r_l}},
\end{eqnarray}
we guarantee the fulfillment of
(\ref{eq2700-new:proof-th:lower-bound-in-L_p})
provided that  $C_6\geq\max_{l=1,\ldots,d}C^{-1/\beta_l}_5$.
Moreover, with this choice (\ref{eq24:proof-th:lower-bound-in-L_p}) is reduced to
\begin{eqnarray}
\label{eq27:proof-th:lower-bound-in-L_p}
&& A^{2+1/\beta}N^{d(1+1/s)}\leq  C_7 L_\beta n^{-1},
\end{eqnarray}
where, as before,  $L_\beta=\prod_{l=1}^d L_l^{1/\beta_l}$.
%$C_5=(8C_4)^{-1}(2C_1)^{-1/\beta}$.
Moreover, we have from (\ref{eq20:proof-th:lower-bound-in-L_p})
\begin{eqnarray}
\label{eq26:proof-th:lower-bound-in-L_p}
&& \rho_n=C_8 AN^{d/p},\;\;\;C_8=C_3(160\kappa)^{-1/p}.
\end{eqnarray}
\par
Let $N^d=C_9A^{-1}$, where  constant $C_9\leq\kappa^{d}$ will be specified below;
then
(\ref{eq17:proof-th:lower-bound-in-L_p}) holds.
Next,  in view of (\ref{eq27:proof-th:lower-bound-in-L_p}) and (\ref{eq26:proof-th:lower-bound-in-L_p})
$$
A=C_{10}(L_\beta/n)^{\frac{1}{1-1/s+1/\beta}},\qquad \rho_n=C_{11}(L_\beta/n)^{\frac{1-1/p}{1-1/s+1/\beta}}=C_{11}\big(L_\beta\alpha_n n^{-1}\big)^{\nu}.
$$
%with sufficiently small  constant $c$.
We remark that $N\to\infty$ as $n\to\infty$.
It remains to check that $\sigma_l$, $l=1,\ldots, d$ are small enough.
It follows from (\ref{eq27-new:proof-th:lower-bound-in-L_p}) that if
$r_l>1$, then $\sigma_l\to 0$ as $n\to\infty$ since $A\to 0$.
If $r_l=1$, then
$$
\sigma_l=C_{12} C_9^{1/\beta_l }
L_l^{- 1/\beta_l}\leq C_{12}\big(C_9/L_{0}\big)^{1/\beta_l}.
$$
Choosing $C_9$ small enough we guarantee that $\sigma_l\leq (20\kappa)^{-1}$, for all $l=1,\ldots,d$.
This condition is required
in the construction of the family  $G_{m}$, $m\in\cM$.
Thus, Lemma~\ref{lem:tsyb_book-result} can be applied with
$\rho_n=C_{11}\big(L_\beta\alpha_n n^{-1}\big)^{\nu}$, and the result follows.
\subsubsection{Dense zone: $\frac{2+1/\beta}{1+1/s} \leq p\leq s(2+\frac{1}{\beta})$}
Here, as in the previous case, we let $m=|\cM|$.
The relationships (\ref{eq27-new:proof-th:lower-bound-in-L_p})
(\ref{eq27:proof-th:lower-bound-in-L_p}) and (\ref{eq26:proof-th:lower-bound-in-L_p})
remain to be true, but our choice of $N$ will be different.
\par
Let $N=C_{12}$ from some constant $C_{12}$. This
yields in view of (\ref{eq27:proof-th:lower-bound-in-L_p}) and
(\ref{eq26:proof-th:lower-bound-in-L_p})
$$
A=C_{13}(L_\beta/n)^{\frac{\beta}{2\beta+1}},\qquad \rho_n=C_{14}(L_\beta/n)^{\frac{\beta}{2\beta+1}}=C_{14}\big(L_\beta\alpha_n n^{-1}\big)^{\nu}.
$$
The requirement (\ref{eq17:proof-th:lower-bound-in-L_p}) is
obviously fulfilled since $A\to 0,\;n\to\infty$. Moreover, we obtain from
(\ref{eq27-new:proof-th:lower-bound-in-L_p}) that $\sigma_l\to 0$
 as $n\to\infty$ and, therefore,
$\sigma_l\leq (20\kappa)^{-1}$, $l=1,\ldots,d$ for $n$ large enough.
Thus, Lemma~\ref{lem:tsyb_book-result} can be applied with
$\rho_n=C_{14}\big(L_\beta\alpha_n n^{-1}\big)^{\nu}$ and the result follows.
%%%%%%%%%%%%%%%%%%%%%%%%%%%%%%%%%%%%%%%%%%%%%%%%%%%%%%%%%%%%%%%%%
\iffalse
\par
In fact, we note that
\begin{eqnarray*}
\frac{1}{2+ 1/\beta}&=& \frac{1}{2+ 1/\beta}\;\wedge\;
\frac{1- 1/p}{1- 1/s + 1/\beta},
\qquad p\geq \frac{2+ 1/\beta}{1+ 1/s};
\\*[2mm]
\frac{1- 1/p}{1- 1/s+ 1/\beta}&=&
\frac{1}{2+ 1/\beta}\;\wedge\;
\frac{1- 1/p}{1- 1/s + 1/\beta},
\qquad p \leq \frac{2+1/\beta}{1+ 1/s},
\end{eqnarray*}
which means that
Lemma~\ref{lem:tsyb_book-result} can be applied with
\begin{eqnarray}
\label{eq28:proof-th:lower-bound-in-L_p}
&&\rho_n=
C_{13}\max\Big\{
(L_\beta n^{-1})^{\frac{1-1/p}{1-1/s+1/\beta}},\; (L_\beta n^{-1})^{\frac{1}{2+1/\beta}}\Big\}.
\end{eqnarray}
\fi
%%%%%%%%%%%%%%%%%%%%%%%%%%%%%%%%%%%%%%%%%%%%%%%%%%%%%%%%%%%%%%%%%%%%%%%%%%%
\subsubsection{Sparse zone: $s(2+\frac{1}{\beta})<p<\infty$, $s<1$}

Let
 $A=\tilde{C}$ and $N=C_{17}$ and suppose that
$\tilde{C}\leq C^{-1}_{17} \kappa^{d}$;
then
(\ref{eq17:proof-th:lower-bound-in-L_p})
is satisfied.
Moreover (\ref{eq24:proof-th:lower-bound-in-L_p}) and (\ref{eq26-new:proof-th:lower-bound-in-L_p})
are reduced to
\begin{eqnarray}
\label{eq240-new:proof-th:lower-bound-in-L_p}
&&n  \prod_{j=1}^d \sigma_j \leq \tilde{C}^{-2}C_{18},\qquad
\sigma_l^{-\beta_l} \bigg(m\prod_{j=1}^d \sigma_j\bigg)^{1/r_l}\leq \widetilde{C}^{-1}C_{19}L_l,
\quad\forall l=1,\ldots,d.
\end{eqnarray}
Let $\tilde{c_1}$, $\tilde{c}_2$ be constants satisfying $\tilde{c_1}\leq \tilde{C}^{-1}C_{18}$,
and $\tilde{c}_2\leq \tilde{C}^{-1}C_{19}$. It is straightforward to check that if we
choose
%
%According to (\ref{eq240-new:proof-th:lower-bound-in-L_p}) we choose
\begin{eqnarray}
\label{eq2400-new:proof-th:lower-bound-in-L_p}
m=\tilde{c}_1^{-1+s}\tilde{c}_2^{s/\beta}L_\beta^{s}n^{1-s},\qquad
\sigma_l=
(\tilde{c}_2L_l)^{-1/\beta_l}\left(\tilde{c}_1\tilde{c}_2^{1/\beta}L_\beta n^{-1}\right)^{s/(\beta_lr_l)},
\quad l=1,\ldots,d,
\end{eqnarray}
then  inequalities (\ref{eq240-new:proof-th:lower-bound-in-L_p})
are fulfilled.
With this choice
(\ref{eq20:proof-th:lower-bound-in-L_p}) is reduced to
\begin{eqnarray}
\label{eq241-new:proof-th:lower-bound-in-L_p}
\rho_n= \widetilde{C}C_{17}\bigg(m\prod_{j=1}^d\sigma_j\bigg)^{1/p}=\widetilde{C}C_{20}(L_\beta n^{-1})^{s/p}=\widetilde{C}C_{20}\big(L_\beta\alpha_n n^{-1}\big)^{\nu}.
\end{eqnarray}
It remains to verify that $\sigma_l$ are small enough, and that $m\geq 8$, $|\cM|/\mm \geq 1$.
Note that
$m\to\infty$ as $n\to\infty$ because of $s<1$.
Remind also that
$$
|\cM|=\prod_{l=1}^d M_l= (20\kappa)^{-d}N^d\prod_{l=1}^d\sigma_l^{-1}
= (20\kappa)^{-d} C_{17}^d \tilde{c}_1 n\,;
$$
hence
$|\cM|/m\geq (20\kappa)^{-d}C_{17}^d (\tilde{c}_1\tilde{c}_2^{1/\beta})^{-s}L_0^{-s/\beta}n^s$.
Thus $|\cM|/m\geq 1$ for large enough $n$.
%
%and, therefore, $|\cM|/\mm\geq \hat{c}L_0^{s/\beta}$  since $s\geq 0$. Therefore, choosing $\hat{c}$ large enough we guarantee that
%$|\cM|/\mm \geq 1$.
We note also that
$\sigma_l\leq (\tilde{c}_2L_0)^{-1/\beta_l}$
for all $n$ large enough. Therefore, if we choose $\tilde{C}$ large enough and put
$\tilde{c}_2=\tilde{C}^{-1} C_{19}$ we can ensure that
$\sigma_l\leq (20\kappa)^{-1}$ for all $l=1,\ldots,d$.
%The latter requirement, in its turn, is possible if one chooses $\widetilde{C}$ small enough
%in order to guarantee $\tilde{c}^{-\beta_l}\hat{c}^{1/r_l}\leq \widetilde{C}^{-1}C_{16}$ for all $l=\overline{1,d}$.
Thus, Lemma~\ref{lem:tsyb_book-result} can be applied with $\rho_n=\widetilde{C}C_{20}\big(L_\beta\alpha_n n^{-1}\big)^{\nu}$ and the result follows.
%Finally remark that condition
%$p\geq s(2+ 1/\beta)$ is equivalent to $s/p\leq \beta/(2\beta+1)$ which means that in view of  (\ref{eq28:proof-th:lower-bound-in-L_p}) and (\ref{eq241-new:proof-th:lower-bound-in-L_p})
%\begin{eqnarray}
%\label{eq280:proof-th:lower-bound-in-L_p}
%&&\rho_n=
%C_{22}\max\left[
%(L_\beta n^{-1})^{\frac{1-1/p}{1-1/s+1/\beta}}, (L_\beta n^{-1})^{\frac{1}{2+1/\beta}},\mathbf{1}_{s<1}(L_\beta n^{-1})^{\frac{s}{p}}\right].
%\end{eqnarray}

%
\subsubsection{Sparse zone: $s(2+\frac{1}{\beta})<p<\infty$, $s\geq 1$}
Here we consider another choice of the set $W$. Let
$W=\{e_1,e_2,\ldots,e_{|\cM|}\}$, where $e_j$,  $j=1,\ldots,|\cM|$ is the canonical
basis in $\bR^{|\cM|}$. With this choice
$$
S_W=1,\qquad |W|= N^d\prod_{j=1}^d\sigma_j^{-1},
$$
and (\ref{eq18:proof-th:lower-bound-in-L_p})
holds with $B=1$. Let  $N=C_{14}$; then
(\ref{eq14:proof-th:lower-bound-in-L_p}) and (\ref{eq22:proof-th:lower-bound-in-L_p})
take the form
\begin{eqnarray}
\label{eq29:proof-th:lower-bound-in-L_p}
A\sigma_l^{-\beta_l}
\bigg(\prod_{j=1}^d \sigma_j\bigg)^{1/r_l}
&\leq&(2C_1)^{-1} L_l,\quad\forall l=1,\ldots, d;
\\
\label{eq30:proof-th:lower-bound-in-L_p}
A^{2}n\prod_{j=1}^d \sigma_j &\leq& C_{15}\ln{\Big(\prod_{j=1}^d \sigma^{-1}_j\Big)}.
\end{eqnarray}
Moreover, we get from  (\ref{eq20:proof-th:lower-bound-in-L_p})
\begin{eqnarray}
\label{eq31:proof-th:lower-bound-in-L_p}
&& \rho_n=C_{16} A\Big(\prod_{j=1}^d\sigma_j\Big)^{1/p}.
\end{eqnarray}
Put $\e=\sqrt{\ln n/n}$ and
\begin{eqnarray}
\label{eq31-new:proof-th:lower-bound-in-L_p}
&&\;\;A=c_1L_\beta^{\frac{1}{2-2/s+1/\beta}}\e^{\frac{1-1/s}{1-1/s+
1/(2\beta)}},\qquad
\sigma_{l}=c_2L_\beta^{\frac{1- 2/r_{l}}{\beta_l(2-2/s+1/\beta)}}
\e^{\frac{1-1/s+1/(\beta r_l)}{\beta_l(1-1/s+
1/(2\beta))}} L_l^{-1/\beta_l}.
%\sigma_{l}=c_2c_1^{1/\beta_{l}- 2/\beta_{l}r_{l}}L_\beta^{\frac{1/\beta_{l}- 2/(\beta_{l}r_{l})}{2-2/s+1/\beta}}
%\e^{\frac{1-1/s+1/\beta r_l}{\beta_l(1-1/s+1/2\beta)}}L_l^{-1/\beta_l}.
\end{eqnarray}
We have
\begin{eqnarray*}
%\label{eq32-new:proof-th:lower-bound-in-L_p}
&&\prod_{l=1}^d\sigma_{l}=c_2^{d}L_\beta^{-\frac{2}{2-2/s + 1/\beta}}
\e^{\frac{1/\beta}{1- 1/s +
1/(2\beta)}},
\end{eqnarray*}
and it is evident that
$\prod_{l=1}^d\sigma_{l}\leq \e^{1/(\beta+1/2)}$ for all $n$ large enough;
hence
$\ln(\prod_{l=1}^d\sigma^{-1}_{l})\geq \ln n/(2\beta+1)$.
Then is is easily checked
that our choice (\ref{eq31-new:proof-th:lower-bound-in-L_p}) satisfies
(\ref{eq29:proof-th:lower-bound-in-L_p}) and (\ref{eq30:proof-th:lower-bound-in-L_p})
provided that
\begin{eqnarray}
\label{eq33-new:proof-th:lower-bound-in-L_p}
&& c_1\leq (2C_1)^{-1},\qquad c_2\leq 1\qquad
c_1^2c_2^{d}\leq C_{15}/(2\beta+1).
\end{eqnarray}
Here we have also used that $d-1/s\geq 0$. Note also that if $s>1$ then
$$
A\to 0,\qquad \max_{l=1,\ldots, d}\sigma_l\to 0,\quad n\to\infty,
$$
which ensures
(\ref{eq17:proof-th:lower-bound-in-L_p}) and
$\sigma_l\leq (20\kappa)^{-1}$, $l=1,\ldots, d$ for all n large enough.
\par
On the other hand, if $s=1$ then
we should add to (\ref{eq33-new:proof-th:lower-bound-in-L_p}) the conditions
$$
c_1L_\beta^{\frac{1}{2-2/s+1/\beta}}\leq C_{14}\kappa^{d},\qquad
c_2\max_{l=1,\ldots,d}\bigg[L_\beta^{\frac{1/\beta_{l}- 2/(\beta_{l}r_{l})}{2-2/s+1/\beta}}
 L_0^{-1/\beta_l}\bigg]\leq (20\kappa)^{-1}.
$$
Obviously, both restrictions hold if we choose $c_1$ and $c_2$ small enough,
but now these constants may depend on $\vec{L}$.
Note, however, that if
$\max_{l=1,\ldots,d}L_l\leq L_\infty$ then $c_1$ and $c_2$ can be chosen depending
on $L_0$ and $L_\infty$ only.
\par
Using (\ref{eq31:proof-th:lower-bound-in-L_p}) and (\ref{eq31-new:proof-th:lower-bound-in-L_p})
we conclude that Lemma \ref{lem:tsyb_book-result} is applicable with
\begin{eqnarray}\label{eq:sparse-s>1}
\rho_n=C_{16}L_\beta^{\frac{1/2- 1/p}{1- 1/s + 1/(2\beta)}}
\left(\frac{\ln n}{n}\right)^{\frac{1- 1/s + 1/(p\beta)}{2(1- 1/s +
1/(2\beta))}}=C_{16}L_\beta^{\frac{1/2- 1/p}{1- 1/s + 1/(2\beta)}-\nu}\left(\frac{L_\beta\alpha_n}{n}\right)^{\nu}.
\end{eqnarray}
that completes the proof of statement~(i) of the theorem.
%%%%%%%%%%%%%%%%%%%%%%%%%%%%%%%%%%%%%%%%%%%%%%%%%%%%%%%%%%%%%%%%%%%%%%%%
\iffalse
It remains to note that
\begin{eqnarray*}
\frac{1- 1/s+  1/(p\beta)}{2(1- 1/s+
1/(2\beta))}&=& \frac{\beta}{2\beta+1}\;\wedge\;
\frac{1- 1/s + 1/(p\beta)}{2(1- 1/s+
1/(2\beta))}\,,
\qquad
\end{eqnarray*}
whenever  $p\geq s(2+1/\beta)$ and  $s\geq 1$.
This fact together with
%(\ref{eq241-new:proof-th:lower-bound-in-L_p}),
(\ref{eq280:proof-th:lower-bound-in-L_p}) and (\ref{eq:sparse-s>1})
\fi
%%%%%%%%%%%%%%%%%%%%%%%%%%%%%%%%%%%%%%%%%%%%%%%%%%%%%%%%%%%%%%%%%%%

%
\subsubsection{Proof of statement~(ii): sparse zone, $p=\infty$, $s\leq 1$}
The proof in this case coincides with the one
for the sparse zone with $s<1$.
Thus, we keep (\ref{eq240-new:proof-th:lower-bound-in-L_p}),
(\ref{eq2400-new:proof-th:lower-bound-in-L_p}), and, in view of (\ref{eq190:proof-th:lower-bound-in-L_p}),
(\ref{eq241-new:proof-th:lower-bound-in-L_p}) is replaced by
$
\rho_n=\widetilde{C}C_{17}.
$
Since $\rho_n$ does not tend to $0$ as $n\to\infty$,
a consistent estimator does not exist. All other details of the proof remain unchanged.
%
%It remains to note that $\mm\geq \hat{c}L_0^{s/\beta}$ and choosing $\hat{c}$ large enough we guarantee that $\mm>8$.
%To obey simultaneously $\tilde{c}^{-\beta_l}\hat{c}^{1/r_l}\leq \widetilde{C}^{-1}C_{16}$ for all $l=\overline{1,d}$
%it suffices to choose $\widetilde{C}$
%small enough.
This completes the proof of Theorem~\ref{th:lower-bound-in-L_p}.
\epr

\subsection{Proof of statement~(ii) of Theorem~\ref{th:upper-new}}
\label{subsec:(ii)-th-4}
The proof goes along the lines of the proof of
Theorem~\ref{th:lower-bound-in-L_p} with modifications indicated below.

\par
We start with the following simple observation: for any $M>0$ and $y>0$ one has
\begin{equation}
\label{eq:inclusion-comp-sup-*}
 \|g\|_\infty\leq M,\;\;
{\rm supp}\{g\}\subseteq [-y,y]^d \;\;\Rightarrow\;\; g\in
\bG_\theta\big(M(2y+4)^{d/\theta}\big),\;\;\forall\theta\in (0,1].
\end{equation}
This is an immediate consequence of the fact that
conditions $\|g\|_\infty\leq M$, ${\rm supp}\{g\}\subseteq [-y,y]^d$
imply that $\|g^*\|_\infty\leq M$ and ${\rm supp}\{g\}\subseteq [-y-2, y+2]^d$.
\iffalse
For any $y>0$ denote $\bK_y$ the set of all borel functions compactly supported on $[-y,y]^d$. Then, for any  $M>0$ and $y>0$
 \begin{eqnarray}
 \label{eq:inclusion-comp-sup-*}
&&\bF(M)\cap\bK_y\subset\bG_\theta\left(M(2y+4)^{d/\theta}\right),\quad\forall\theta\in (0,1).
 \end{eqnarray}
The latter inclusion follows from the obvious observation: if $g\in\bF(M)\cap\bK_y$ then $g^*\in\bF(M)\cap\bK_{y+2}$.
\fi
\par
Next, we note that the lower bounds of
Theorem~\ref{th:lower-bound-in-L_p} in the dense and sparse zones
are proved over the set of compactly
supported densities. Hence they are valid also on
$\bG_\theta(R)\cap\bN_{\vec{r},d}(\vec{\beta},\vec{L}, M)$, provided that $R$ is large enough.
Hence, if $p\geq \frac{2+1/\beta}{1/\theta+1/s}$  the assertion of the theorem follows.
\par
Let $p< \frac{2+1/\beta}{1/\theta+1/s}$.
The proof of the lower bound here differs from
the proof of Theorem~\ref{th:lower-bound-in-L_p} only in
construction of the function $f^{(0)}$.
\par
Let $f^{(0)}$ be the function constructed exactly as
in the proof of Theorem~\ref{th:lower-bound-in-L_p}
with $N=N_0$ fixed throughout the asymptotics $n\to \infty$,
and such that
$f^{(0)}\in \bN_{\vec{r},d}(\vec{\beta},4^{-1}\vec{L}_0, 4^{-1}M)$.
Since $N_0$ is fixed, $f^{(0)}$ is compactly supported
and, by (\ref{eq:inclusion-comp-sup-*}) we have that  $f^{(0)}\in\bG_\theta(R_1)$ for some large enough $R_1>0$.
Define
$$
\bar{f}^{(\theta)}(x)=
\prod_{l=1}^d\bigg[ N^{-1/\theta}
\int_{\bR}\Lambda(y-x_l)\mathbf{1}_{[-\frac{N}{2},\frac{N}{2}]}(y)\rd y\bigg],\;\;\;\;\;
x=(x_1,\ldots,x_d)\in \bR^d,
$$
where $N=N(n)\to\infty$ will be specified later. Let
$\tilde{f}^{(\theta)}(x)=\varsigma^{d}\bar{f}^{(\theta)}(\varsigma x)$,
where $\varsigma>0$ is chosen to guarantee
$\tilde{f}^{(\theta)}\in \bN_{\vec{r},d}(\vec{\beta},4^{-1}\vec{L}_0, 4^{-1}M)$.
We note however that, in contrast to the case $\theta=1$,
$\tilde{f}^{(\theta)}$ is not a probability density. In particular,
$\int \tilde{f}^{(\theta)}\to 0$ as $N\to \infty$,
because $\theta<1$.
Define
$$
f^{(\theta)}=(1-p_N)f^{(0)}+\tilde{f}^{(\theta)},
$$
where $p_N:=\int \tilde{f}^{(\theta)}$ ensures $\int f^{(\theta)}=1$.
%Obviously $p_N\to 0,\; N\to \infty$ and $p_N\to 0,\;\theta\to 1$.

Note also that $f^{(1)}=\tilde{f}^{(1)}$ since $\tilde{f}^{(1)}$ is a probability density and, therefore $p_N=1$.
Thus, we can assert that
\begin{eqnarray*}
%\label{eq:f-theta}
f^{(\theta)}\in
\bN_{\vec{r},d}(\vec{\beta},2^{-1}\vec{L}_0, 2^{-1}M),\quad \int f^{(\theta)}=1,\quad f^{(\theta)}\geq 0.
\end{eqnarray*}
Note that, by construction, $\tilde{f}^{(\theta)}$
is supported on the cube $\left[(-N/2-1)/\varsigma,(N/2+1)/\varsigma\right]^d$
and bounded by $N^{-d/\theta}\varsigma^d$. Therefore, in view of (\ref{eq:inclusion-comp-sup-*}),
$\tilde{f}^{(\theta)}\in\bG_\theta(R_2)$
for some large enough $R_2$.
\par
Let $W$ be the parameter set as defined in the proof
of Theorem~\ref{th:lower-bound-in-L_p}.  For any $w\in W$ and any
$\theta\leq 1$
we let
$$
f^{(\theta)}_w(x)=f^{(\theta)}(x)+F_{w}(x),\;\;\;\;\; x\in\bR^d,
$$
where functions $F_w$ are constructed as in the proof of Theorem~\ref{th:lower-bound-in-L_p}.
If instead of (\ref{eq17:proof-th:lower-bound-in-L_p}) we require
\begin{eqnarray}
\label{eq17-new:proof-th:lower-bound-in-L_p}
&& A\leq \big[\kappa^d+ \varsigma^{d}\big]N^{-d/\theta},
\end{eqnarray}
then we obtain in view of
(\ref{eq8-new:proof-th:lower-bound-in-L_p}) and
%(\ref{eq16:proof-th:lower-bound-in-L_p}), (\ref{eq:inclusion-comp-sup-*}) and
(\ref{eq17-new:proof-th:lower-bound-in-L_p})
that $\{F_w,\;w\in W\}\subset\bG_\theta(R_3)$ for some large enough~$R_3$.
All said above one allows to conclude that
$\{f^{(\theta)},\;f^{(\theta)}_w,\;w\in W\}$ is a finite set of probability densities
from
$\bG_\theta(R)\cap\bN_{\vec{r},d}(\vec{\beta},\vec{L}, M)$ for some large enough $R>0$,
 and Lemma~\ref{lem:tsyb_book-result} is applicable with $\cJ_n=W$ and $\bF=\bG_\theta(R)\cap \bN_{\vec{r},d}(\vec{\beta},\vec{L}, M)$.

Note also that if $\theta=1$ we come to the construction used in the proof of Theorem \ref{th:lower-bound-in-L_p} and, therefore, the statement of the theorem 
in the case $\theta=1$ follows.
 
\par
Suppose now that $\theta<1$. We will follow construction of the set $W$ for the tail zone which is given in
Subsection~\ref{subsec:tail}.  Choose $m=|\cM|$ and
 note that (\ref{eq2700-new:proof-th:lower-bound-in-L_p}),
(\ref{eq27-new:proof-th:lower-bound-in-L_p}),
(\ref{eq26:proof-th:lower-bound-in-L_p}) remain unchanged, while
(\ref{eq27:proof-th:lower-bound-in-L_p})  should be replaced by
\begin{eqnarray}
\label{eq27-new-new:proof-th:lower-bound-in-L_p}
&& A^{2+ 1/\beta}N^{d(1/\theta+1/s)}\leq  C_7 L_\beta n^{-1}.
\end{eqnarray}
Now we choose
$N^d=cA^{-\theta}$ with  $c\leq\kappa^{d}+\varsigma^{d}$; then
(\ref{eq17-new:proof-th:lower-bound-in-L_p}) is valid.
We obtain from (\ref{eq27-new-new:proof-th:lower-bound-in-L_p})
that
$$
A=C_{8}(L_\beta/n)^{\frac{1}{1-\theta/s+1/\beta}},\qquad
\rho_n=C_{9}(L_\beta/n)^{\frac{1-\theta/p}{1-\theta/s+1/\beta}}.
$$
Finally, because
(\ref{eq27-new:proof-th:lower-bound-in-L_p}) remains intact,
$\sigma_l\to 0$ as $n\to \infty$ for any $l=1,\ldots, d$;
this follows from  $A\to 0$ and $\theta<1$.
This  completes the proof.
\epr

\subsection{Proof of the lower bound in (\ref{eq:lower-bound-new})}
\label{sec:proof-of-eq:lower-bound-new}
The required result will follow
from the lower bound
of Theorem~\ref{th:lower-bound-in-L_p}
in the tail zone (see Section~\ref{subsec:tail}) if
we will show that
for any given $R>0$ and $\theta\in (0,1)$
\begin{eqnarray}
\label{eq300-new-new:}
f^{(0)}\notin\bG_\theta(R),\;\;f_w \notin \bG_\theta(R),\;\;\;\forall w\in W.
%\{f^{(0)},\;f_w,\;w\in W\}\nsubseteq\bG_\theta(R).
\end{eqnarray}
First we note that
$f^{(0)}=N^{-d}$ for $x\in [-(N-2)/(2\kappa), (N-2)/(2\kappa)\big]^d$;
therefore,
$\|[f^{(0)}]^*\|_\theta\to\infty$ as $N\to\infty$,
because $\theta<1$.
\par
Next, in view of (\ref{eq1600:proof-th:lower-bound-in-L_p}),
$f_w(x)=f^{(0)}(x)$ for any  $x\notin [-(N-4)/(4\kappa),
\;(N-4)/(4\kappa)]^d$, which also implies
$$
\inf_{w\in W}\big\|f_w^*\big\|_\theta\to\infty,\quad N\to\infty.
$$
It remains to note that in the tail
zone the parameter $N$ is chosen so that $N=N(n)\to\infty$ as $n\to\infty$.
This completes the proof
of (\ref{eq300-new-new:}).
\epr

\appendix\label{app:A}
\section{Proofs of auxiliary results of Section~\ref{sec:Proof-th-3}}

\subsection{Measurability}\label{meas}
Write $f(x,X^{(n)}):=\hat{f}_{\hat{h}(x)}(x)$, and note that
the map $f:\bR^d\times\bR^{dn}\to\bR$
is  completely determined by the kernel $K$ and the set $\cH$.
We need to show that $f$ is a Borel function.
\par
Let $R_h(x,X^{(n)}):=\hat{R}_h(x)$, and note that for every $h\in \cH$,
the map
$R_h:\bR^d\times\bR^{dn}\to\bR$ is a continuous function.
This follows from the continuity of  the kernel $K$ and from the fact that
$\cH$ is a finite set. The continuity of $K$ also implies that the map
$f_h:\bR^d\times\bR^{dn}\to\bR$ is a continuous function for any $h\in\cH$,
where $f_h(x,X^{(n)}):=\hat{f}_h(x)$.
Next, denote by $\mB$  the Borel $\sigma$-algebra
on $\bR^d\times\bR^{dn}$, and let
$b:\bR^d\times\bR^{dn}\to\cH$ be the function
$b(x,X^{(n)}):=\hat{h}(x)$. We obviously have for any given $h\in\cH$
$$
\Big\{(x,y)\in\bR^d\times\bR^{dn}: b(x,y)=h\Big\}=
\bigcup_{\eta\in\cH}\Big\{(x,y)\in\bR^d\times\bR^{dn}:
R_h(x,y)-R_\eta(x,y)\leq 0\Big\}\in\mB,
$$
where the last inclusion
follows from the continuity of $R_\eta$, $\eta\in\cH$.
Here we have also used that $\cH$ is a finite set.
It remains to note that
$$
\hat{f}_{\hat{h}(x)}(x)=\sum_{h\in\cH}f_h(x,X^{(n)})
\mathbf{1}\big\{ b\big(x,X^{(n)}\big)=h\big\},
$$
and the required statement follows.
\subsection{Proof of Lemma \ref{lem:proportional}}
\label{proofs-lem:proportional-and-lem:1}
%\paragraph{Proof of Lemma \ref{lem:proportional}}
1$^0$. Note that $\check{M}_h(g,x)=4^{-1} \hat{M}_h(g,x)$ and
let $\cH_0=\{h\in \cH: A_h(g,x)\geq 4\kappa\ln n/(nV_h)\}$.

For any $h\in \cH_0$ we have
\[
\sqrt{\frac{\kappa A_h(g,x)\ln n}{nV_h}}\geq \frac{2\kappa\ln n}{nV_h}\;\;\Rightarrow\;\;
M_h(g,x)\leq \frac{3}{4} A_h(g,x).
\]
Therefore,
\begin{eqnarray*}
 |\hat{A}_h(g,x)-A_h(g,x)| \leq \chi_h(g, x) + M_h(g,x)\leq \chi_h(g,x) + (3/4) A_h(g,x).
\end{eqnarray*}
We have for any $h\in \cH_0$
\begin{eqnarray*}
 \left|\check{M}_h(g,x) - M_h(g,x)\right|& =&\left|\sqrt{\frac{\kappa \ln n}{nV_h}}
\frac{\hat{A}_h(g,x)-A_h(g,x)}{\hat{A}_h^{1/2}(g,x) + A_h^{1/2}(g,x)}\right|
\\*[2mm]
&\leq&  \sqrt{\frac{\kappa \ln n}{nV_h}}\left(
\frac{\chi_h(g,x) + (3/4) A_h(g,x)}{ A_h^{1/2}(g,x)}\right)
\leq \frac{1}{2}\chi_h(g,x)+\frac{3}{4}M_h(g,x).
\end{eqnarray*}
It yields for any $h\in \cH_0$
\begin{eqnarray}\label{eq:M-bounds}
&& \big[\check{M}_h(g,x) - \frac{7}{4}M_h(g,x)\big]_+\leq \frac{1}{2}\chi_h(g,x),\;\;\;\;
\big[M_h(g,x) - 4\check{M}_h(g,x)\big]_+\leq 2\chi_h(g,x).
\end{eqnarray}
\par\medskip
(b). Now consider the set $\cH_1:=\cH\setminus \cH_0$. Here  $A_h(g,x)\leq 4\kappa\ln n/(nV_h)$,
and, by definition of~$M_h$ we have
\begin{eqnarray}\label{eq:A-00}
 \frac{1}{4} A_h(g,x) \leq \frac{\kappa\ln n}{nV_h} \leq M_h(g,x) \leq \frac{3\kappa\ln n}{nV_h},\;\;\;
\forall h\in \cH_1.
\end{eqnarray}
Note that we have $\hat{M}_h(g,x)\geq \kappa\ln n /(nV_h)$ for all $h$. This together with
(\ref{eq:A-00}) shows that
\begin{equation}\label{eq:000}
 [M_h(g,x)-3 \tilde{M}_h(g,x)]_+ =0,\;\;\;\forall h\in \cH_1.
\end{equation}
Furthermore, for any $h\in \cH_1$
\begin{eqnarray*}
 \hat{A}_h(g,x) \leq A_h(g,x)+\chi_h(g,x) + M_h(g,x) \leq \chi_h(g,x) + \frac{7\kappa \ln n}{nV_h}.
\end{eqnarray*}
Therefore
\begin{eqnarray*}
\check{M}_h(g,x) &=& \sqrt{\frac{\kappa \hat{A}_h(g,x)\ln n}{nV_h}} + \frac{\kappa\ln n}{nV_h}
\;\leq\; \sqrt{\frac{\kappa \chi_h(g,x) \ln n}{nV_h}} + (\sqrt{7}+1)\frac{\kappa\ln n}{nV_h}
\\*[2mm]
&\leq& \frac{1}{2}\chi_h(g,x) + \Big(\sqrt{7}+\frac{3}{2}\Big) \frac{\kappa\ln n}{nV_h}
\;\leq\; \frac{1}{2}\chi_h(g,x) + \Big(\sqrt{7}+\frac{3}{2}\Big) M_h(g,x).
\end{eqnarray*}
To get the penultimate inequality we have used that $\sqrt{|ab|}\leq 2^{-1}(|a|+|b|)$.
Thus, it is shown that
\begin{equation}\label{eq:M-00}
 \Big[\check{M}_h(g,x)- \Big(\sqrt{7}+\frac{3}{2}\Big) M_h(g,x)\Big]_+ \leq \frac{1}{2}\chi_h(g,x),\;\;\;
\forall h\in \cH_1.
\end{equation}
Relations (\ref{eq:M-00}), (\ref{eq:000}) and (\ref{eq:M-bounds}) imply statement of the lemma.
\epr
\subsection{Proof of Lemma \ref{lem:1}}
1$^0$.
Let $g:\bR^d\to\bR^1$ be a fixed bounded function, and let
\[
 \xi_h(g, x) := \frac{1}{n}\sum_{i=1}^n g_h(X_i-x) - \int g_h(t-x)f(t)\rd t,\;\;\;h\in \cH.
\]
With this notation
$\xi_h(x) = \xi_h(K,x)$ and $\hat{A}_h(g,x)-A_h(g, x) = \xi_h(|g|, x)$.
Therefore moment bounds on $\zeta_1(x)$, $\zeta_3(x)$ and $\zeta_4(x)$  will follow from those on
$\xi_h(g,x)$ with substitution $g \in \{ K,Q, |K|, |Q|\}$.
Since $M_h(g,x)$ depends on $g$ only via $|g|$ and $\|g\|_\infty$
[see (\ref{eq:A-A})--(\ref{eq:A-M})],
$M_h(g, x)=M_h(|g|,x)$,
and moment bounds
on  $\zeta_1(x)$ and $\zeta_3(x)$ are identical. The bound on $\zeta_4(x)$
will follow from bounds on $\zeta_1(x)$ and $\zeta_3(x)$ with only one modification:
kernel $K$ should be replaced by $Q$.
As for $\zeta_2(x)$,  $\xi_{h,\eta}(x)$ cannot be represented in terms of $\xi_h(g,x)$ with
function $g$ independent of $h$ and $\eta$; see (\ref{eq:Q}). However,
the bounds on $\zeta_2(x)$ will be obtained similarly  with minor modifications.
Thus it suffices to bound $\bE_f[\zeta_1(x)]^q$ and $\bE_f[\zeta_2(x)]^q$.
\par
2$^0$.
We start with bounding $\bE_f[\zeta_1(x)]^q$.
For any $z>0$, $h\in \cH$ and $q\geq 1$ one has
\begin{equation}\label{eq:ineq-0}
 \bE_f \bigg[|\xi_h(x)| - \sqrt{\frac{2\rk_\infty A_h(K,x)z}{nV_h}}-
\frac{2\rk_\infty z}{3nV_h}\bigg]_+^q \leq
2\Gamma(q+1)
\bigg[\sqrt{\frac{2\rk_\infty A_h(K,x)}{nV_h}} + \frac{2\rk_\infty }{3nV_h}\bigg]^q e^{-z}.
\end{equation}
This inequality follows by integration of the Bernstein inequality and the following bound
on the second moment of $\xi_h(x)$:
\begin{eqnarray*}
 \bE_f|\xi_h(x)|^2 &\leq& \frac{\rk_\infty}{nV_h}\int |K_h(t-x)| f(t) \rd t =
\frac{\rk_\infty A_h(K,x)}{nV_h}.
\end{eqnarray*}
\par
We will show that $\bE_f[\zeta_1(x)]^q$ is bounded by the expression appearing
on the right hand side of (\ref{eq:ineq-1}). In fact,
we will prove a stronger inequality.
Let for some $l>0$
\begin{equation*}
 \lambda_h:= (1+q)\ln \big(1/V_h\big) +
 \ln\big(F^{-1}(x) \wedge n^{l}\big).
\end{equation*}
%here $\ln_+(\cdot)=\max\{\ln(\cdot), 0\}$.
In suffices to show that
(\ref{eq:ineq-1}) holds when in the definition of $\zeta_1(x)$ the quantity
$M_h(K,x)$ replaced by $\tilde{M}_h(K,x)$, where
\[
 \tilde{M}_h(K,x)=\sqrt{\frac{2\rk_\infty A_h(K,x)\lambda_h}{nV_h}} + \frac{2\rk_\infty\lambda_h}{3nV_h}.
\]
Indeed, since  and $n^{-d}\leq V_h\leq 1$ for any $h\in\cH$, we have that
\[
\lambda_h \leq (q+1)d \ln n +  l\ln n.
\]
Therefore $\tilde{M}_h(K,x) \leq M_h(K,x)$ for all $x\in \bR^d$ and $h\in \cH$ provided that
%Since $V_{\max}\leq 1$ and $V_h\geq V_{\min}\geq n^{-1}$, we have that
%$\lambda_h \leq 1+ (\frac{q}{2}+1) \ln n + sl\ln n$. Therefore
%$M_h(x)\leq \tilde{M}_h(x)$ for all $x\in \bR^d$ and $h\in \cH$ provided that
\[
\kappa\geq \rk_\infty\big[d(2q+4)+2l\big].
\]
Thus if we establish
(\ref{eq:ineq-1}) with $M_h(K,x)$ replaced by $\tilde{M}_h(K,x)$, the
required bound for $\bE_f[\zeta_1(x)]^q$ will be proved.
\par
We have for any $h\in \cH$
\begin{eqnarray*}
 \exp\{-\lambda_h\} = \big(V_h\big)^{q+1}
\big\{F(x)\vee n^{-l}\big\}.
\end{eqnarray*}
Furthermore, taking into account that $A_h(g,x)\leq V_h^{-1}\|g\|_\infty$ for any $g$, we obtain
\[
\sqrt{\frac{2\rk_\infty A_h(K,x)}{nV_h}} + \frac{2\rk_\infty }{3nV_h}
\leq \frac{2\rk_\infty}{\sqrt{n}V_h}.
\]
Here we have used that $n\geq 3$.
If we set $z=\lambda_h$ then
(\ref{eq:ineq-0}) together with two previous display formulas yields
\begin{eqnarray}
\bE_f[\zeta_1(x)]^q &=&\bE_f \sup_{h\in \cH} \Big[|\xi_h(x)|- M_h(K,x)\Big]^q_+  \;\leq\;
\sum_{h\in \cH} \bE_f \Big[|\xi_h(x)|- \tilde{M}_h(K,x)\Big]^q_+
\nonumber
\\
&\leq&2\Gamma(q+1)  \big(2\rk_\infty\big)^q  n^{-q/2}
\big\{F(x)\vee n^{-l}\big\}
\sum_{h\in \cH} V_h
\nonumber
\\
&\leq&
2^{d+1}\Gamma(q+1)\big(2\rk_\infty\big)^q  n^{-q/2}
\big\{F(x)\vee n^{-l}\big\}.
\label{eq:zeta-2}
\end{eqnarray}
As it was mentioned above, under the same conditions
inequality (\ref{eq:zeta-2}) holds for $\bE_f[\zeta_3(x)]^q$.
\par
As for the moment bound for  $\zeta_4(x)$, in all formulas
above $K$ should be replaced by $Q$ and $\rk_\infty$ by $\rk^2_\infty$ since
\mbox{$\|Q\|_\infty\leq \rk_\infty^2$}. Specifically, if $\kappa\geq \rk^2_\infty[d(2q+4)+2l]$
then
\begin {eqnarray*}
 \bE_f[\zeta_4(x)]^q
&\leq&
2^{d+1}\Gamma(q+1)\big(2\rk^2_\infty\big)^q n^{-q/2}
\big\{F(x)\vee n^{-l}\big\}.
\label{eq:zeta-1}
\end{eqnarray*}
\par\medskip
3$^0$.  Now we turn to bounding $\bE_f[\zeta_2(x)]^q$. We have similarly to
(\ref{eq:ineq-0})
\begin{eqnarray*}
&& \bE_f \bigg[|\xi_{h,\eta}(x)| - \sqrt{\frac{2\rk^2_\infty A_{h\vee\eta}(Q,x)z}{nV_{h\vee \eta}}}-
\frac{2\rk^2_\infty z}{3nV_{h\vee \eta}}\bigg]_+^q
\nonumber
\\*[2mm]
&&\;\;\;\;\hspace{15mm} \leq\;
2\Gamma(q+1)
\bigg[\sqrt{\frac{2\rk^2_\infty A_{h\vee\eta}(Q,x)}{nV_{h\vee \eta}}} +
\frac{2\rk^2_\infty }{3nV_{h\vee \eta}}\bigg]^q e^{-z}.
%\label{eq:ineq-2}
\end{eqnarray*}
Here we have used
the following bound on the second moment of $\xi_{h,\eta}(x)$:
\begin{eqnarray*}
\bE_f|\xi_{h,\eta}(x)|^2 &\leq& \frac{\|Q_{h,\eta}\|_\infty}{nV_{h\vee \eta}} \int
\Big|\frac{1}{V_{h\vee\eta}}Q_{h,\eta}\Big(\frac{t-x}{h\vee\eta}\Big)\Big| f(t)
\rd t
\\
&\leq& \frac{\|Q\|_\infty}{nV_{h\vee \eta}} \int
|Q_{h\vee\eta}(t-x)| f(t) \rd t = \frac{\rk_\infty^2 A_{h\vee \eta}(Q,x)}{nV_{h\vee\eta}}~.
\end{eqnarray*}
The further proof goes along the same lines as the above proof
with the following minor modifications:
in all formulas $\rk_\infty$ should be replaced with
$\rk_\infty^2$, $V_{h\vee \eta}$ should be written instead of $V_h$, and $\kappa$ should satisfy
$\kappa\geq \rk^2_\infty[d(2q+4)+2l]$.  The statement of the lemma holds with constant
$C_0=2^{d^2+1}\Gamma(q+1)(2\rk^2_\infty)^q$.
Combining the above bounds we complete the proof.
\epr

\section{Proofs of auxiliary results of
Section~\ref{sec:proofs}}

\subsection{Proof of Lemma \ref{lem:bias-norm-bound}}
\label{sec:bias-lemma}
% and
%Propositions \ref{prop:first-bound} and \ref{prop:new}}
%\label{sec:proof_tech-results-th-upper}
%
%\subsubsection{Proof of Lemma \ref{lem:bias-norm-bound}}
%
We have
\begin{eqnarray*}
&& B_h(f,x) = \int K(u) \big[f(x+uh)-f(x)\big] \rd u
=  \int \prod_{j=1}^d w_\ell(u_j) \big[f(x+uh)-f(x)\big] \rd u.
\end{eqnarray*}
First, we note that $f(x+uh)-f(x)$ can be represented by the telescopic sum
\begin{equation}\label{eq:telescope}
 f(x+uh)-f(x) = \sum_{j=1}^{d} \Delta_{u_jh_j, j} f(x_1,\ldots,x_j, x_{j+1}+u_{j+1}h_{j+1},\ldots, x_d+
u_dh_d),
\end{equation}
where we put formally $h_{d+1}u_{d+1}=0$.
\par
Next, for any function $g:\bR^d\to\bR^1$ and $j=1,\ldots, d$ we have
\begin{eqnarray}
&&\int w_\ell(u_j) \Delta_{u_jh_j, j}g(x)\rd u_j =
\int \sum_{i=1}^\ell \binom{\ell}{i} (-1)^{i+1}\frac{1}{i} w\Big(\frac{u_j}{i}\Big)\Delta_{u_jh_j, j}g(x)
\rd u_j
\nonumber
\\
&&= (-1)^{\ell-1}\int w(z) \sum_{i=1}^\ell \binom{\ell}{i} (-1)^{i+\ell}\Delta_{izh_j, j}g(x)
\rd z
%\nonumber\\
= (-1)^{\ell-1} \int w(z) \Delta^\ell_{zh_j, j}\, g(x) \rd z.
\label{eq:int-representation}
\end{eqnarray}
The last equality follows from the definition of $\ell$-th order difference operator
%in the last line we have used that for any $u\in \bR^1$ and $j=1,\ldots, d$
%\[
%\Delta_{u, j}^\ell g(x)=(-1)^{\ell-1}\Big[ \sum_{i=1}^\ell \binom{\ell}{i}(-1)^{i+1}g(x+ue_j)-g(x)\Big];
%\]
(\ref{eq:Delta}).
Thus (\ref{eq:int-representation}) and (\ref{eq:telescope}) imply that
$B_h(f,x) =\sum_{j=1}^d B_{h,j}(f,x)$,
where
\begin{equation*}\label{eq:B-h-j}
\begin{array}{lll}
&& (-1)^{\ell-1}  B_{h,j}(f,x) :=
\\*[2mm]
&&
{\displaystyle
\int
  \int w(z) \Delta^\ell_{zh_j, j}\, f(x_1,\ldots,x_j,x_{j+1}+u_{j+1}h_{j+1},
\ldots,x_d+u_dh_d) \rd z \prod_{m=j+1}^d w_\ell(u_m)\rd u_m.
}
\end{array}
\end{equation*}
Therefore, by the Minkowski inequality for integrals [see, e.g., \cite[Section~6.3]{folland}]
\begin{eqnarray*}
&& \big\|B_{h,j}(f, \cdot)\big\|_{r_j}
\\
&&\;
\leq\; \int
  \int |w(z)|\Big\| \Delta^\ell_{zh_j, j}\, f(\cdot,\ldots,\cdot,\cdot+u_{j+1}h_{j+1},
\ldots,\cdot+u_dh_d)\Big\|_{r_j} \rd z \prod_{m=j+1}^d |w_\ell(u_m)| \rd u_m
\\
&&\;=\int
  \int |w(z)|\big\| \Delta^\ell_{zh_j, j}\, f\big\|_{r_j} \rd z \prod_{m=j+1}^d |w_\ell(u_m)| \rd u_m.
\end{eqnarray*}
Since $f\in \bN_{\vec{r},d}\big(\vec{\beta},\vec{L}\big)$ one has
\begin{eqnarray*}
\big\|B_{h,j}(f, \cdot)\big\|_{r_j}\leq\;
L_j h_j^{\beta_j} \int\int |w(z)|\,|z|^{\beta_j} \rd z \prod_{m=j+1}^d |w_\ell(u_m)| \rd u_m
\;\leq\; C_1 L_jh_j^{\beta_j}.
\end{eqnarray*}
This proves (\ref{eq:bias-norms-1}). To get (\ref{eq:bias-norms-2}) we first note that the condition $s\geq 1$ implies
$\tau(p)>0$ and $\tau_j>0$, $j=1,\ldots,d$. Then
the inequality in (\ref{eq:bias-norms-2}) follows by the same reasoning with $r_j$ replaced by
$q_j$, $\beta_j$ replaced by $\gamma_j$
and with the use of embedding (\ref{eq:embedd-nik}).
\epr
\subsection{Proof of Proposition \ref{prop:first-bound}}
\label{sec:proof-prop-1}
By definition of  $J_m$ and $\cX_m$,
\begin{equation}\label{eq:J-m}
J_m \leq 2^{p(m+1)}\varphi^p|\cX_m|,
\end{equation}
and
now we bound from above $|\cX_m|$.
By definition of $\cX_m$ we have for any $h\in \cH$
\begin{eqnarray}
 |\cX_m|\leq \Big|\{x: \sup_{\eta\geq h} M_\eta (x) >2^{m-1} \varphi\}\Big|
+ \sum_{j=1}^d \Big|\{x: B_{h, j}^*(f,x)>2^{m-1} \varphi\}\Big|
&&
\nonumber
\\
=: J_{m,1}(h) + J_{m,2}(h). &&
\label{eq:X-m}
\end{eqnarray}
Recall that with the introduced notation
\[
M_\eta(x)=\sqrt{\kappa A_\eta(x) \delta V^{-1}_\eta} + \kappa\delta V_\eta^{-1}, \;\;\;\eta\in \cH.
\]
\par
For any $h\in\cH$ we have
\begin{eqnarray}
 J_{m,2}(h) &=&
\sum_{j\in I\setminus I_\infty} \Big|\{x: B_{h, j}^*(f,x)>2^{m-1} \varphi\}\Big|
\;+\;
\sum_{j\in I_\infty} \Big|\{x: B_{h, j}^*(f,x)>2^{m-1} \varphi\}\Big|
\nonumber
\\
&=:& J_{m,2}^{(1)}(h) + J_{m,2}^{(2)}(h).
\label{eq:J-m-2-2}
\end{eqnarray}
By the Chebyshev inequality and (\ref{eq:B*-1}) for any $h$
\begin{eqnarray}\label{eq:J-m-1}
&&J_{m,2}^{(1)}(h)
 \leq \sum_{j\in I\setminus I_\infty} \big[2^{(m-1)}\varphi\big]^{-r_j}
\|B^*_{h, j}(f,\cdot)\|_{r_j}^{r_j}
\;\leq \; c_1  \sum_{j\in I\setminus I_\infty} 2^{-r_jm}\varphi^{-r_j}  L_j^{r_j}h_j^{\beta_jr_j}.
\end{eqnarray}
In addition, if  $s\geq 1$
then the Chebyshev inequality and (\ref{eq:B*-2}) yield
\begin{eqnarray}\label{eq-new:J-m-1}
&&J_{m,2}^{(1)}(h)
\leq \sum_{j\in  j\in I\setminus I_\infty} \big[2^{(m-1)}\varphi\big]^{-q_j}
\|B^*_{h, j}(f,\cdot)\|_{q_j}^{q_j}
% \nonumber\\&=&
\leq \tilde{c}_1 \sum_{j\in I\setminus I_\infty} 2^{-q_jm}\varphi^{-q_j}  L_j^{q_j}h_j^{\gamma_jq_j}.
\end{eqnarray}
%where, recall, $\vec{\gamma}$ and $\vec{q}$ are defined in (\ref{eq:gamma-and-q}).
%Here we have also used (\ref{eq:B*-2}).
\par
In order to prove
statements (i)--(iv) of the proposition we bound quantities
$J_{m,1}(h)$ and $J_{m,2}(h)$ with bandwidth
$h=h[m]$ specified in an appropriate way.
\subsubsection{Proof of statements {\rm (i)} and {\rm (ii)}}
Note that the bound in the statement  
{\rm (i)} coincides formally with that 
of the statement {\rm (ii)} when $\theta=1$. 
This implies that the statement {\rm (ii)}
in the case $\theta=1$ follows from the statement {\rm (i)}. 
So, with slight abuse of notation, 
we will identify the case $\theta=1$ with 
the assumption that $f$ is a probability density.

1$^0$. We start with bounding the term $J_{m,1}(h)$ on the right hand side of (\ref{eq:X-m}).
Assume that $h\in \cH$ is such that
\begin{equation}\label{eq:cond-0}
\kappa \delta V_h^{-1} <2^{m-2}\varphi;
\end{equation}
then by  the Chebyshev inequality
\begin{eqnarray*}
 J_{m, 1}(h) &\leq & \Big|
\big\{x: \sup_{\eta\geq h}\sqrt{\kappa A_\eta(x)\delta V^{-1}_\eta} > 2^{m-2}\varphi\big\}\Big|
\nonumber
\\
&\leq& \sum_{\substack{\eta\geq h\\ \eta\in \cH}}
\Big| \{x: A_\eta (x) > 2^{2m-4}\varphi^2 \kappa^{-1}\delta^{-1}V_\eta\}\Big|
\nonumber
\\
&\leq& \left(2^{-2m+4}\varphi^{-2}\delta\kappa\right)^{\theta} \sum_{\substack{\eta\geq h\\ \eta\in \cH}} \|A_\eta\|^{\theta}_\theta V_\eta^{-\theta}
\;\leq\; c_2 \left(2^{-2m}\varphi^{-2}\delta\right)^{\theta} \sum_{\substack{\eta\geq h\\ \eta\in \cH}} V_\eta^{-\theta},
% V_h^{-1},
\end{eqnarray*}
where we have taken into account that, for any $\eta$,
$\|A_\eta\|_\theta\leq R$ if $f\in\bG_\theta(R)$ with $\theta<1$, and
$\|A_\eta\|_1\leq \rk_\infty^2$.
By definition of $\cH$,
for any $\eta\geq h$, $\eta\in \cH,$ we have
$V_\eta = V_h 2^{k_1+\cdots+k_d}$ for some $k_1,\ldots,k_d\geq 1$, which implies that
$\sum_{\eta\geq h} V_\eta^{-\theta} \leq (1-2^{-\theta})^{-d} V_h^{-\theta}$.
Thus, we conclude that for any $h$ satisfying (\ref{eq:cond-0}) one has
\begin{eqnarray}
\label{eq:J-m-2}
&& J_{m, 1}(h) \leq c_3 \left(2^{-2m}\varphi^{-2}\delta V_h^{-1}\right)^{\theta}.
\end{eqnarray}
\par
2$^0$. Let $\tilde{h}=(\tilde{h}_1,\ldots,\tilde{h}_d)\in (0,\infty]^d$ be given by
\begin{eqnarray}
\label{eq01:proof-prop}
&&\tilde{h}_j=\big(c_4L_j^{-1}\varphi)^{1/\beta_j}
2^{\frac{m}{\beta_j}\left(1-\frac{\theta(2+1/\beta)}{r_j(1+\theta/s)}\right)},
\;\;j=1,\ldots,d,
\end{eqnarray}
where constant $c_4$ will be specified later.
Let us prove that
$\tilde{h}\in \big[n^{-1},1\big]^{d}$ for large enough $n$.
\par
Denote
\[
a=\frac{1-\theta/s+1/\beta}{1+\theta/s}, \;\;\;\;b_j=1-\frac{\theta(2+1/\beta)}{r_j(1+1/s)},
\]
and remark that $a>0$. We note also that
$$
a^{-1} b_j=\frac{(2+1/\beta)(1-\theta/r_j)}{1-\theta/s+1/\beta}-1.
$$
\par
If
$b_j\leq 0$
then, because $m\leq 0$,
\[
 \tilde{h}_j \geq \big(c_4L_j^{-1}\varphi)^{1/\beta_j} = (c_4L_j^{-1})^{1/\beta_j}
(L_\beta \delta)^{\frac{1}{\beta_j(2+1/\beta)}} > n^{-1}
\]
for all large enough $n$. On the other hand, since $0\geq m\geq m_0(\theta)$ and
$2^{m_0(\theta)a} \leq 2^{a}\hat{c}_1\kappa \varphi$ by definition of $m_0(\theta)$,
\begin{eqnarray*}
 \tilde{h}_j \leq (c_4L_j^{-1}\varphi)^{1/\beta_j}
2^{\frac{m_0(\theta)b_j}{\beta_j}}=(c_4L_j^{-1}\varphi)^{1/\beta_j}
\big(2^{m_0(\theta)a}\big)^{\frac{b_j}{a\beta_j}}
\\
\leq (c_4L_j^{-1})^{1/\beta_j} (2^{-a}\hat{c}_1\kappa)^{\frac{b_j}{a\beta_j}}
\varphi^{\frac{1+b_j/a}{\beta_j}} \leq c_5 (c_4L_0^{-1})^{1/\beta_j},
\end{eqnarray*}
where we took into account that $1+b_j a^{-1}>0$ and $\min_{j=1,\ldots,d} L_j\geq L_0>0$.
Then choosing constant $c_4$ small enough we have $\tilde{h}_j\leq 1$.
Thus we showed that $\tilde{h}_j\in [n^{-1}, 1]$ for $j$ such that  $b_j\leq 0$.
\par
Now consider the case $b_j>0$. Here
\begin{eqnarray*}
\tilde{h}_j \geq
(c_4L_j^{-1}\varphi)^{1/\beta_j} \big(2^{m_0(\theta)a}\big)^{b_ja^{-1}/\beta_j} \geq
(c_4L_j^{-1}\varphi)^{1/\beta_j}\big(\hat{c}_1\kappa\varphi \big)^{b_ja^{-1}/\beta_j}
\nonumber
\\
\geq
c_6\big(c_4L_j^{-1})^{1/\beta_j}
\varphi^{\frac{(2+1/\beta)(1-\theta/r_j)}{\beta_j(1-\theta/s+1/\beta)}}~.
%\label{eq02:proof-prop}
%\\
%&&
%\big(c_4L_j^{-1})^{1/\beta_j}
%2^{m\ma\left(\frac{(2+1/\beta)(1-\theta/r_j)}{\beta_j(1-\theta/s+1/\beta)}-\frac{1}{\beta_j}\right)}
%\varphi^{1/\beta_j}\geq c_5\big(c_4L_j^{-1})^{1/\beta_j}
%\varphi^{\frac{(2+1/\beta)(1-\theta/r_j)}{\beta_j(1-\theta/s+1/\beta)}}
\end{eqnarray*}
It remains to note that
$$
\frac{1-\theta/r_j}{\beta_j(1-\theta/s+1/\beta)}
=\frac{1/\beta_j-\theta/(\beta_jr_j)}{1-\theta/s+1/\beta}<1,\quad\forall j=1,\ldots,d,
$$
in view of the obvious inequality $1/\beta-\theta/s\geq 1/\beta_j-\theta/(\beta_jr_j)$,
which, in its turn, follows from
the fact that $\theta\in (0,1]$.
Thus, we have that $\tilde{h}_j>n^{-1}$ for all large enough $n$.
Furthermore, if $b_j>0$ then  since $m\leq 0$
\[
 \tilde{h}_j \leq (c_4L_j^{-1})^{1/\beta_j} \varphi^{1/\beta_j} \leq 1
\]
for all large enough $n$. Thus we have shown that $\tilde{h}\in [n^{-1}, 1]^d$.
\par
3$^0$.
Now we proceed with bounding $J_{m,2}(h)$ for a specific choice of $h=h[m]$, which is defined as follows.
Let
$h[m]\in \cH$ such that $h[m]< \tilde{h}\leq 2h[m]$.
Let constant $c_4$
in (\ref{eq01:proof-prop}) be chosen so that  $c_4<(2\bar{c}_1)^{-1}$, where $\bar{c}_1$ appears
on the right hand side of (\ref{eq:B*-1}).
With this choice of $c_4$ by (\ref{eq:B*-1})
$$\|B^*_{h[m],j}(f,\cdot)\|_\infty\leq
\bar{c}_1 L_j(h_j[m])^{\beta_j}\leq \bar{c}_1 L_j\tilde{h}_j^{\beta_j}\leq 2^{m-1}\varphi.
$$
Therefore, $J_{m,2}^{(2)}\big(h[m]\big)=0$, where $J_{m,2}^{(2)}(\cdot)$ is defined in (\ref{eq:J-m-2-2}).
Moreover, we obtain from (\ref{eq:J-m-1}) and $h[m]\leq \tilde{h}_j$ that
\begin{equation}
\label{eq03:proof-prop}
J_{m,2}^{(1)}\big(h[m]\big)\leq
c_1\sum_{j\in I\setminus I_\infty} 2^{-r_jm}\varphi^{-r_j} L_j^{r_j}\tilde{h}_j^{\beta_j r_j} \leq
c_1\sum_{j\in I\setminus I_\infty} c^{r_j}_4
2^{-m \left(\frac{2+1/\beta}{1/\theta+1/s}\right)} \leq  c_72^{-m\left(\frac{2+1/\beta}{1/\theta+1/s}\right)}~.
\end{equation}
Note that
\begin{eqnarray}\label{eq:V-h-tilde}
V_{h[m]}\geq 2^{-d} V_{\tilde{h}}= 2^{-d}c_4^{1/\beta}
L_\beta^{-1}\varphi^{1/\beta}2^{m\left(\frac{1}{\beta}-\frac{2+1/\beta}{1+s/\theta}\right)}.
\end{eqnarray}
This together with  (\ref{eq:J-m-2}) yields
\begin{eqnarray}
\label{eq04:proof-prop}
J_{m,1}\big(h[m]\big)\leq c_{8}2^{-m\left(\frac{2+1/\beta}{1/\theta+1/s}\right)}.
\end{eqnarray}
Then it follows from (\ref{eq03:proof-prop}) and (\ref{eq04:proof-prop})
that
\begin{eqnarray*}\label{eq:j1+j2}
J_{m,1}\big(h[m]\big)+J_{m,2}\big(h[m]\big)\leq
c_{9}2^{-m\big(\frac{2+1/\beta}{1/\theta+1/s}\big)},
\end{eqnarray*}
which combined with (\ref{eq:J-m}) results in
 \begin{equation*}%\label{eq:J-m-bound-0}
 J_m \leq c_{10} 2^{m\big(p-\frac{2+1/\beta}{1/\theta+1/s}\big)} \varphi^p.
%,\;\;\;b_1=\frac{2+1/\beta}{1+1/s}.
\end{equation*}
\par
Inequality (\ref{eq:J-m}) is valid only if
(\ref{eq:cond-0}) is fulfilled for $h[m]$, i.e.,
$\kappa \delta V_{h[m]}^{-1} <2^{m-2}\varphi$; now we verify this condition.
It is sufficient to check that
$\kappa\delta 2^{d}V^{-1}_{\tilde{h}}< 2^{m-2}\varphi$. In view of
(\ref{eq:V-h-tilde})  this inequality  will follow if
\[
 c_4^{1/\beta}\varphi^{1+ 1/\beta}2^{m\left(1+\frac{1}{\beta}-
\frac{2+1/\beta}{1+s/\theta}\right)}> 2^{d+2}\kappa (L_\beta\delta).
\]
Taking into account that $L_\beta\delta=\varphi^{2+1/\beta}$
we conclude that (\ref{eq:cond-0}) is fulfilled for $h[m]$ if
\[
2^{m\left(\frac{1-\theta/s+1/\beta}{1+\theta/s}\right)} >
\hat{c}_1\kappa \varphi,
\]
which is ensured by the condition $m\geq m_0(\theta)$.
This completes the proof of (\ref{eq01:m-large}).

\iffalse
\smallskip
\noindent\textsf{Proof of} (\ref{eq01-new:m})
 $1^{0}.\;$ Let $\hat{c}_4$ be the constant whose choice will be done later and
let $C_j=c_4$ if $j\in I_\infty$ and $C_j=\hat{c}_4$ if $j\in I\setminus I_\infty$. Here the constant $c_4$ is the same as before.
Let $\tilde{h}=(\tilde{h}_1,\ldots,\tilde{h}_d)\in (0,\infty]^d$ be given by
\begin{eqnarray}
\label{eq1-new:proof-prop}
&&\tilde{h}_j= \big(C_jL_j^{-1})^{1/\beta_j}\varphi^{\frac{1}{\beta_j}+\frac{s}{\beta_jr_j}}
2^{m\left(\frac{1}{\beta_j}-\frac{s(1+1/\beta)}{\beta_jr_j}\right)},\;\;j=\overline{1,d}.
\end{eqnarray}
Note that if $j\in I_\infty$ the corresponding coordinates of $\tilde{h}$ given by (\ref{eq01:proof-prop}) and (\ref{eq1-new:proof-prop}) are the same.

\smallskip

$2^{0}.\;$ Let us show that $\tilde{h}\in \big[n^{-1},1]^d$. Indeed, since $m\leq 0$ then if $\ma_j:=1-\frac{s(1+1/\beta)}{r_j}\geq 0$ we have
taking into account that $2^m\varphi\geq \hat{c}^{-1}\delta\;\Leftrightarrow\; 2^m\geq \hat{c}^{-1}L_\beta^{-1}\varphi^{1+1/\beta}$
$$
\tilde{h}_j\geq C_j(\vec{L})\varphi^{(2+1/\beta)\left((\beta_j-\right))}
> \delta> n^{-1},
$$
\fi

\subsubsection{Proof of statement~{\rm (iii)}}
1$^{0}$. Let $\hat{c}_4$ be a constant to be specified later, and
let $c_4$ be the constant given in
(\ref{eq01:proof-prop}). Let $C_j=c_4$ if $j\in I_\infty$ and
$C_j=\hat{c}_4$ if $j\in I\setminus I_\infty$.
Define  $\tilde{h}=(\tilde{h}_1,\ldots,\tilde{h}_d)\in (0,\infty]^d$ by the formula
\begin{eqnarray}
\label{eq1:proof-prop}
&&\tilde{h}_j=\big(C_jL_j^{-1}\varphi)^{1/\beta_j}
2^{m\big(\frac{1}{\beta_j}-\frac{s(2+1/\beta)}{\beta_jr_j}\big)},\;\;j=1,\ldots,d.
\end{eqnarray}
Note that if $j\in I_\infty$ the corresponding
coordinates of $\tilde{h}$ given by (\ref{eq01:proof-prop}) and (\ref{eq1:proof-prop}) are the same.
\par
Let us show that $\tilde{h}\in \big[n^{-1},1]^d$ for large enough $n$.
First consider the coordinates $\tilde{h}_j$ such that
$1-\frac{s}{r_j}(2+1/\beta)\geq 0$.
Because $m\geq 0$ we have for all $n$ large enough
$$
\tilde{h}_j\geq
\big(C_jL_j^{-1}\varphi)^{1/\beta_j}\geq
\big(C_jL_j^{-1})^{1/\beta_j}(L_\beta\delta)^{\frac{1}{2\beta_j+\beta_j/\beta}}
> \delta> n^{-1},
$$
where we have used the obvious inequality  $\beta_j/\beta>1$ for any $j=1,\ldots, d$.
On the other hand,
because $2^m\leq \hat{c}_2\varphi^{-1}$ we obtain
\begin{eqnarray*}
\tilde{h}_j\leq
c_{11}\big(\hat{c}_4L_j^{-1})^{1/\beta_j}\varphi^{\frac{s(2+1/\beta)}{\beta_jr_j}}\to 0,\;\;n\to\infty, \quad\forall j\in I\setminus I_\infty;
\\*[2mm]
\tilde{h}_j\leq c_{11}\big(C_jL_j^{-1})^{1/\beta_j}\varphi^{\frac{s(2+1/\beta)}{\beta_jr_j}}\leq
c_{11}\big(c_4L_0^{-1})^{1/\beta_j},
\quad\forall j\in I_\infty.
\end{eqnarray*}
Thus $\tilde{h}_j \leq 1$ for large enough $n$ if $j\in I\setminus I_\infty$, and $\tilde{h}_j\leq 1$
by choice of constant $c_4$ if $j\in I_\infty$.
\par
Now consider the case $1-\frac{s}{r_j}(2+1/\beta)< 0$.
Since $2^m\leq \hat{c}_2\varphi^{-1}$
$$
\tilde{h}_j\geq c_{12}\big(C_jL_j^{-1})^{1/\beta_j}\varphi^{\frac{s(2+1/\beta)}{\beta_jr_j}}=
c_{12}\big(C_jL_j^{-1})^{1/\beta_j}
(L_\beta\delta)^{\frac{s}{\beta_jr_j}}> \delta> n^{-1},
$$
for all $n$ large enough. Here we have used the obvious inequality $1/s>1/\beta_jr_j$
$\forall j=1,\ldots,d$.
On the other hand, since $m\geq 0$,
$\tilde{h}_j\leq \big(C_jL_j^{-1}\varphi)^{1/\beta_j}\leq 1$ for large enough $n$.
Thus we have proved that $\tilde{h}\in [n^{-1}, 1]^d$ for all large enough $n$.
\par
2$^{0}$. Let
$h[m]\in \cH$ such that $h[m]< \tilde{h}\leq 2h[m]$ and choose
constant $c_4$ satisfies $c_4<(2\bar{c}_1)^{-1}$ [see (\ref{eq:B*-1})].
Recall that formulas (\ref{eq01:proof-prop}) and (\ref{eq1:proof-prop})
coincide for $j\in I_\infty$. Therefore,
as before, with the indicated choice of $c_4$ we have
\begin{eqnarray}
\label{eq2:proof-prop}
&&J_{m,2}^{(2)}\big(h[m]\big)=0.
\end{eqnarray}
\par
Let $\beta_{\pm}$ and $\beta_\infty$ be defined by expressions
$
1/\beta_\pm:=\sum_{j\in I_+\cup I_-} 1/\beta_j$ and
$1/\beta_\infty:=\sum_{j\in I_\infty} 1/\beta_j$.
We  have
\begin{eqnarray}
\label{eq3:proof-prop}
&& V_{h[m]}\geq 2^{-d} V_{\tilde{h}}=2^{-d}
c_4^{1/\beta_\infty}\hat{c}_4^{1/\beta_\pm}L_\beta^{-1}\varphi^{1/\beta}2^{-2m}.
\end{eqnarray}
This together with $2^m\leq \hat{c}_2\varphi^{-1}$ shows that
$V_{h[m]}\geq c_{13}\hat{c}_4^{1/\beta_\pm} L^{-1}_\beta\varphi^{2+1/\beta}=c_{13}
\hat{c}_4^{1/\beta_\pm}\delta$ and, therefore,
\begin{equation}\label{eq:cond-1}
 \kappa\delta V_{h[m]}^{-1} \leq c_{13}^{-1}\hat{c}_4^{-1/\beta_\pm}\kappa.
\end{equation}
Remark that $A_\eta(x) \leq 2^{d}M\rk_\infty^2 $ for all $x\in \bR^d$ and
$\eta\in \cH$.
% and $f\in\bF(M)$.
Hence,  in view of (\ref{eq:cond-1})
\begin{eqnarray*}
\sup_{\eta\geq h[m]} M_\eta (x)
&\leq&
\rk_\infty\sqrt{2^{d}M} \sqrt{\kappa \delta V_{h[m]}^{-1}} + \kappa \delta V_{h[m]}^{-1}
\\
&\leq& \bigg(\rk_\infty
\sqrt{2^{d}M}+\sqrt{c_{13}^{-1}\hat{c}_4^{-1/\beta_\pm}}\bigg)
\sqrt{\delta V_{h[m]}^{-1}}
=c_{14}\sqrt{\delta V_{h[m]}^{-1}}.
%\big(c_{16}+c_{17}\big(\hat{c}_4)^{1/2\beta_\pm}\big)\sqrt{\delta V_h^{-1}}.
\end{eqnarray*}
It yields together with (\ref{eq3:proof-prop})
$$
\sup_{\eta\geq \tilde{h}} M_\eta (x)\leq c_{15}\,\hat{c}_4^{-1/(2\beta_\pm)}2^{m}
\sqrt{(L_\beta\delta)\varphi^{-1/\beta}}=c_{15}\,\hat{c}_4^{-1/(2\beta_\pm)} 2^{m}\varphi.
$$
Setting $\hat{c}_4$ so that $c_{15}\hat{c}_4^{-1/(2\beta_\pm)}<2^{-1}$, we
obtain $\sup_{\eta\geq \tilde{h}} M_\eta (x) \leq  2^{m-1}\varphi$. This implies that
\begin{eqnarray}
\label{eq4:proof-prop}
J_{m,1}\big(h[m]\big)=0.
\end{eqnarray}
Moreover, it follows from (\ref{eq:J-m-1}) and from inequality $h[m]\leq \tilde{h}$ that
\begin{eqnarray}
\label{eq5:proof-prop}
&&J_{m,2}^{(1)}\big(h[m]\big)\leq J_{m,2}^{(1)}\big(\tilde{h}\big)\leq \bigg[c_1\sum_{j\in I\setminus I_\infty} \hat{c}^{r_j}_4\bigg]
2^{-ms\left(2+1/\beta\right)}.
\end{eqnarray}
Then
(\ref{eq02:m}) is a consequence of
(\ref{eq:J-m}), (\ref{eq2:proof-prop}), (\ref{eq4:proof-prop}) and (\ref{eq5:proof-prop}).
The statement~(ii) is proved.
\subsubsection{Proof of statement~{\rm (iv)}}
1$^{0}$. Let $C_j$, $j=1,\ldots, d$ be the same constants in the proof of statement~(ii)
in the previous section.  Define
$\tilde{h}=(\tilde{h}_1,\ldots,\tilde{h}_d)\in (0,\infty]^d$ by the following formula
\begin{eqnarray}
\label{eq100:proof-prop}
\tilde{h}_j&=&
\big(C_jL_j^{-1}\varphi)^{1/\gamma_j}
2^{m\left(\frac{1}{\gamma_j}-\frac{\upsilon(2+1/\gamma)}{\gamma_jq_j}\right)}
\left[\frac{L_\gamma\varphi^{1/\beta}}{L_\beta\varphi^{1/\gamma}}\right]^{\frac{\upsilon}{\gamma_jq_j}},
\;\;\;\;j=1,\ldots,d,
\end{eqnarray}
where $\gamma_j$, $q_j$ are defined in (\ref{eq:gamma-and-q}) and
$\gamma$, $\upsilon$ and $L_\gamma$ are given in (\ref{eq:gamma-upsilon}).
\par
Let us show that $\tilde{h}\in \big[n^{-1},1]^d$ for large $n$.
Let $b_j=1-\frac{\upsilon}{q_j}(2+1/\gamma)$.
\par
First, assume that $b_j<0$.
Since $m> 0$ and $2^m\leq\hat{c}_2\varphi^{-1}$,
\begin{eqnarray*}
%&&
\tilde{h}_j \geq c_{16}\big(C_jL_j^{-1})^{1/\gamma_j}\varphi^{\frac{\upsilon(2+1/\gamma)}{\gamma_jq_j}}
\left[\frac{L_\gamma\varphi^{1/\beta}}{L_\beta\varphi^{1/\gamma}}\right]^{\frac{\upsilon}{\gamma_jq_j}}
= c_{16}\big(C_jL_j^{-1})^{1/\gamma_j}
\left[L_\gamma/L_\beta\right]^{\frac{\upsilon}{\gamma_jq_j}}\varphi^{\frac{\upsilon(2+1/\beta)}{\gamma_jq_j}}
\\
%&&
= c_{16}\big(C_jL_j^{-1})^{1/\gamma_j}
\left[L_\gamma\right]^{\frac{\upsilon}{\gamma_jq_j}}\delta^{\frac{\upsilon}{\gamma_jq_j}}
>\delta>n^{-1},
\end{eqnarray*}
where we have used the obvious inequality $1/\upsilon>1/(\gamma_jq_j)$ for any $j=1,\ldots,d$.
On the other hand,
in view of $m\geq m_1$ and
by (\ref{eq:m-1-def})
\begin{eqnarray*}
\tilde{h}_j &\leq&
\big(C_jL_j^{-1}\varphi)^{1/\gamma_j}2^{\frac{m_1b_j}{\gamma_j}}
\big[(L_\gamma/L_\beta) \varphi^{1/\beta-1/\gamma}\big]^{\frac{\upsilon}{\gamma_jq_j}}
\\
&\leq&
(C_jL_j^{-1}\varphi)^{1/\gamma_j}2^{\frac{m_1b_j}{\gamma_j}}
\left[ 2^{m_1[\upsilon(2+1/\gamma)-s(2+1/\beta)]}\right]^{\frac{1}{\gamma_jq_j}}
= \big(C_jL_j^{-1}\varphi)^{1/\gamma_j}
2^{m_1\big(\frac{1}{\gamma_j}-\frac{s(2+1/\beta)}{\gamma_jq_j}\big)}.
\end{eqnarray*}
Then by (\ref{eq:m-1-bounds})
\begin{equation}\label{eq:h-tilde<}
\tilde{h}_j\leq c_{17}C_1(\vec{L})C_{j}^{1/\gamma_j}
\varphi^{\frac{1}{\gamma_j}\big(1+[1-\frac{s}{q_j}{(2+1/\beta)]
\frac{\upsilon(1/\beta-1/\gamma)}{\upsilon[2+1/\gamma]-s[2+1/\beta]}\big)}}~~~,
\end{equation}
where the expression for constant $C_1(\vec{L})$ is easily found.
It remains to note that
\begin{eqnarray*}
%\label{eq10011:proof-prop}
\upsilon(2+1/\gamma)-s(2+1/\beta)=s\upsilon\left[(2+1/\beta)(1/s-1/\upsilon)+
(1/\gamma-1/\beta)s^{-1}\right],
\end{eqnarray*}
and
in view of (\ref{eq:gamma<beta}) and (\ref{eq:s<upsilon})
%\begin{eqnarray}
%\label{eq10011:proof-prop}
%\upsilon[2+1/\gamma]-s[2+1/\beta]=s\upsilon\left([2+1/\beta][1/s-1/\upsilon]+[1/\gamma-1/\beta]s^{-1}\right).
%\end{eqnarray}
%which leads to
\[
1+
\frac{[q_j-s(2+1/\beta)]\upsilon(1/\beta-1/\gamma)}
{q_j[\upsilon(2+1/\gamma)-s(2+1/\beta)]}=\frac{2+1/\beta}{q_j}
\left[\frac{(1/s-1/\upsilon)q_j+(1/\gamma-1/\beta)}{(1/s-1/\upsilon)(2+1/\gamma)+(1/\gamma-1/\beta)\upsilon^{-1}}
\right] >0.
\]
This shows that $\tilde{h}_j \leq 1$ for large $n$.
\par
Now assume that $b_j\geq 0$.
Then, similarly to the reasoning that resulted in (\ref{eq:h-tilde<})
we have
\begin{eqnarray}
\label{eq101:proof-prop}
\tilde{h}_j&\geq&\big(C_jL_j^{-1}\varphi)^{1/\gamma_j}2^{\frac{m_1b_j}{\gamma_j}}
\big[(L_\gamma/L_\beta) \varphi^{1/\beta-1/\gamma}\big]^{\frac{\upsilon}{\gamma_jq_j}}
\nonumber\\
%&=&\big(C_jL_j^{-1}\varphi)^{1/\gamma_j}2^{m_1\left(\frac{1}{\gamma_j}-\frac{s(2+1/\beta)}{\gamma_jq_j}\right)}
%2^{-m_1\left(\frac{\upsilon(2+1/\gamma)-s(2+1/\beta)}{\gamma_jq_j}\right)}
%\left[\frac{L_\gamma\varphi^{1/\beta}}{L_\beta\varphi^{1/\gamma}}\right]^{\frac{\upsilon}{\gamma_jq_j}}
%\nonumber\\
%&\geq& 2^{-\frac{1}{\gamma_jq_j}}\big(C_jL_j^{-1}\varphi)^{1/\gamma_j}
%2^{m_1\left(\frac{1}{\gamma_j}-\frac{s(2+1/\beta)}{\gamma_jq_j}\right)},\;\;j=\overline{1,d}.
%\\*[2mm]
%\tilde{h}_j&=&\big(C_jL_j^{-1}\varphi)^{1/\beta_j}2^{m\left(\frac{1}{\beta_j}-\frac{s(2+1/\beta)}{\beta_jr_j}\right)},\;\;j\in I_-.
&\geq& C_1(\vec{L})C_{j}^{1/\gamma_j}
\varphi^{\frac{1}{\gamma_j}\big(1+\big[1-\frac{s}{q_j}(2+1/\beta)\big]
\frac{\upsilon(1/\beta-1/\gamma)}{\upsilon[2+1/\gamma]-s[2+1/\beta]}\big)}.
\nonumber
\end{eqnarray}
Since $\varphi^{2+1/\beta}=L_\beta\delta$,
$$
\tilde{h}_j\geq c_{18}C_1(\vec{L})C_{j}^{1/\gamma_j}
\delta^{\frac{(1/s-1/\upsilon)(1/\gamma_j)+(1/\gamma-1/\beta)(1/\gamma_j r_j)}{(1/s-1/\upsilon)(2+1/\gamma)+(1/\gamma-1/\beta)\upsilon^{-1}}}>\delta>n^{-1}
$$
for all $n$ large enough. Here we have used
(\ref{eq:gamma<beta}), (\ref{eq:s<upsilon}),
and   obvious inequalities: $2+1/\gamma>1/\gamma_j$
and $1/\upsilon>1/\gamma_j r_j$  for all
$j=1,\ldots, d$.
On the other hand,
%Now we show that $\tilde{h}_j\leq 1$ for large enough $n$ whenever $b_j\geq 0$.
since $2^{m}\leq\hat{c}_2\varphi^{-1}$
\begin{eqnarray*}
\tilde{h}_j&\leq&c_{19}\big(C_jL_j^{-1})^{1/\gamma_j}\varphi^{\frac{\upsilon(2+1/\gamma)}{\gamma_jq_j}}
\left[\frac{L_\gamma\varphi^{1/\beta}}{L_\beta\varphi^{1/\gamma}}\right]^{\frac{\upsilon}{\gamma_jq_j}}
= c_{19}\big(C_jL_j^{-1})^{1/\gamma_j}\left[L_\gamma/L_\beta\right]^{\frac{\upsilon}{\gamma_jq_j}}\varphi^{\frac{\upsilon(2+1/\beta)}{\gamma_jq_j}}
\\
&=&c_{19}\big(C_jL_j^{-1})^{1/\gamma_j}\left[L_\gamma\right]^{\frac{\upsilon}{\gamma_jq_j}}\delta^{\frac{\upsilon}{\gamma_jq_j}}.
\end{eqnarray*}
 Therefore,
$\tilde{h}_j\to 0,\; n\to\infty$, $\forall j\in I\setminus I_\infty$,
and $\tilde{h}_j\leq c_{19}\big(c_4L_0^{-1})^{1/\gamma_j}$,
$\forall j\in I_\infty$.
Choosing $c_4$ small enough we come to required assertion.
\par
3$^{0}$.  Let
$h[m]\in \cH$ be such that $h[m]< \tilde{h}\leq 2h[m]$, and let
constant $c_4$ satisfy $c_4<(2\bar{c}_1)^{-1}$, where $\bar{c}_1$
is given in (\ref{eq:B*-1}).
With this choice of $c_4$, if $j\in I_\infty$ then
the corresponding coordinates of $\tilde{h}$ given
by  (\ref{eq1:proof-prop}) and (\ref{eq100:proof-prop}) coincide.
Hence we have as before
\begin{eqnarray}
\label{eq20:proof-prop}
&&J_{m,2}^{(2)}\big(h[m]\big)=0.
\end{eqnarray}
Let
$\frac{1}{\gamma_\pm}:=\sum_{j\in I_+\cup I_-} \frac{1}{\gamma_j}$;
then
\begin{eqnarray}
\label{eq30:proof-prop}
&& V_{h[m]}\geq 2^{-d}
V_{\tilde{h}}=2^{-d}\big(c_4)^{1/\beta_\infty}\big(\hat{c}_4)^{1/\gamma_\pm}
L_\beta^{-1}\varphi^{1/\beta}2^{-2m}.
\end{eqnarray}
We remark that (\ref{eq30:proof-prop}) and
(\ref{eq3:proof-prop}) coincide up
to the change in notation $\beta_\pm\leftrightarrow\gamma_\pm$.
Hence all the computations preceding
(\ref{eq4:proof-prop}) remain valid, and we have as before
\begin{eqnarray}
\label{eq40:proof-prop}
J_{m,1}\big(h[m]\big)=0.
\end{eqnarray}
Moreover, we obtain from (\ref{eq-new:J-m-1})
\begin{eqnarray}
\label{eq50:proof-prop}
&&J_{m,2}^{(1)}\big(h[m]\big)\leq
\bigg[\tilde{c}_1\sum_{j\in I\setminus I_\infty} \hat{c}^{r_j}_4\bigg]
\bigg[\frac{L_\gamma\varphi^{1/\beta}}{L_\beta\varphi^{1/\gamma}}
\bigg]^{\upsilon}2^{-m\upsilon\left(2+1/\gamma\right)}.
\end{eqnarray}
The bound given in (\ref{eq03:m})
follows now from (\ref{eq:J-m}), (\ref{eq20:proof-prop}),
(\ref{eq40:proof-prop}) and (\ref{eq50:proof-prop}).
\epr
\subsection{Proof of Proposition \ref{prop:new}}\label{sec:proof-prop-2}
In view of (\ref{eq:pointwise-oracle})
and (\ref{eq01:oracle-inequality})
\begin{equation}\label{eq:oracle-inequality-1}
 |\hat{f}(x)-f(x)| \leq c_0 [\bar{U}_f(x)
 %\inf_{h\in\cH}\big\{\bar{B}_h(f,x) + \sup_{\eta\geq h} M_\eta (K\vee Q, x)\big\}
+ \omega(x)] \leq c_1 [U_f(x)
 %\inf_{h\in\cH}\big\{\bar{B}_h(f,x) + \sup_{\eta\geq h} M_\eta (K\vee Q, x)\big\}
+ \omega(x)] ,
\end{equation}
where
$c_0$ and $c_1$ are appropriate constants, $\bar{U}_f(x)$ and $U_f(x)$ are given by
(\ref{eq:U-bar}) and (\ref{eq:U}) respectively, and
$\omega(x):=\zeta(x)+\chi(x)$ with
$\zeta(x)$ and $\chi(x)$ defined in  (\ref{eq:zeta}) and (\ref{eq:chi}).
\subsubsection{Proof of statement~{\rm (i)}}
%1$^0.\;$ For any $a>0$ let $\cX_-(a)=\left\{x\in\bR^d:\;\; \bar{U}_f(x)\leq a\right\}$.
Here for brevity we will write $m_0=m_0(1)$.
By (\ref{eq:oracle-inequality-1})
%for any $a>0$
\begin{eqnarray}
&&\int_{\cX_{m_0}^-} |\hat{f}(x)-f(x)|^p \rd x
\;\leq\;
c_1^{p-1} \int_{\cX_{m_0}^-} \big[U_f(x)+\omega(x)\big]^{p-1}|\hat{f}(x)-f(x)|\rd x
\nonumber
\\
&&\;\;\;\leq\;  c_2\bigg[ (2^{m_0}\varphi)^{p-1}\int_{\bR^d} |\hat{f}(x)-f(x)|\rd x
+
\int_{\bR^d} \omega^{p-1}(x) [2^{m_0}\varphi + \omega(x)] \rd x\bigg].
\nonumber
\end{eqnarray}
Noting that $\|\hat{f}\|_1\leq \ln^d(n)\|K\|_1\leq \ln^d(n)\rk_\infty$,
we have $\|\hat{f}-f\|_1\leq \ln^d(n)\rk_\infty +1$
 and, therefore,
\[
 (2^{m_0}\varphi)^{p-1}\int_{\bR^d} |\hat{f}(x)-f(x)|\rd x \;\leq\; \big(\ln^d(n)\rk_\infty +1\big)(2^{m_0}\varphi)^{p-1}.
\]
Moreover,
since $\kappa= \rk_\infty^2 [(4d+2)p+4(d+1)]$,
the second statement of Theorem~\ref{th:oracle-inequality}
implies
\[
 \bE_f \int_{\bR^d} \omega^{p-1}(x) [2^{m_0}\varphi + \omega(x)] \rd x
\leq c_3  (2^{m_0}\varphi) n^{-(p-1)/2} + c_4 n^{-p/2}.
\]
Combining these inequalities
and taking into account that $2^{m_0}\varphi \leq 1$ we obtain
\begin{eqnarray}
J_{m_0}^-=\bE_f\int_{\cX^-_{m_0}} |\hat{f}(x)-f(x)|^p \rd x
\leq c_5\big[\ln^d(n) (2^{m_0}\varphi)^{p-1}+2^{m_0}\varphi n^{-(p-1)/2} + n^{-p/2}\big]
\nonumber
\\
\leq 2c_5\big[\ln^d(n)(2^{m_0}\varphi)^{p-1}+ n^{-p/2}\big].
%\label{eq1:prop-new}
\nonumber
\end{eqnarray}
By definition of $m_0=m_0(1)$,
$2^{m_0}\varphi\leq c_6 \big(L_\beta \delta\big)^{1/(1+1/\beta -1/s)}$; therefore
$$
J_{m_0}^- \;\leq\; c_7\ln^d(n)\big(L_\beta \delta\big)^{\frac{p-1}{1-1/s+1/\beta }}+c_7n^{-p/2}.
$$
It remains to note that for large $n$
that
$$
\big(L_\beta \delta\big)^{\frac{p-1}{1+1/\beta -1/s}}\leq
\big(L_\beta \delta\big)^{\nu p},\qquad n^{-p/2}<\big(L_\beta \delta\big)^{\nu p},
$$
and  (\ref{eq:prop-new_1}) follows.
\subsubsection{Proof of statement {\rm (ii)}}
Let $f^*$ be
the maximal operator of $f$ defined in
(\ref{eq:maximal-function}).
It follows from the definition of $M_\eta(x)$ that for any $h\in\cH$
\begin{eqnarray*}
\sup_{\eta\geq h} M_\eta(x) \leq  c_8\sqrt{\frac{\kappa f^*(x)\ln n}{nV_h}} + \frac{\kappa\ln n}{nV_h}.
\end{eqnarray*}
Moreover, by definition of $\bar{B}_h(f,x)$,
$\bar{B}_h(f,x) \leq c_{9} [f^*(x)+f(x)]\leq 2c_{9} f^*(x)$ almost everywhere,
where the last inequality follows from the Lebesgue differentiation theorem.
Using these two inequalities  and setting $h=(1,\ldots,1)$ in (\ref{eq:U-bar})
we come to the following upper bound on $\bar{U}_f(x)$
\begin{equation}
\label{eq2:prop-new}
\bar{U}_f(x)\leq c_{10} \big[f^*(x)+ \sqrt{f^*(x)\delta} + \delta\big].
\end{equation}
In view of   (\ref{eq01:oracle-inequality}) we have that
$\cX_{m_0(\theta)}^-\subseteq \cX^-:=\{x\in \bR^d: \bar{U}_f(x) \leq \rk_\infty 2^{m_0(\theta)}\varphi\}$; therefore
if we put
\[
D_1:= \cX^- \cap \{x\in \bR^d: f^*(x) \leq \delta\},\;\;D_2:=\cX^- \cap \{x\in \bR^d: f^*(x) > \delta\}
\]
then
\begin{equation}\label{eq:S1+S2}
 J_{m_0(\theta)}^- \leq \bE_f \int_{D_1} |\hat{f}(x)-f(x)|^p \rd x + \bE_f
\int_{D_2} |\hat{f}(x)-f(x)|^p \rd x =: \bE_f S_1 + \bE_f S_2.
\end{equation}
We bound from above the two terms on the right hand side of the above inequality.
\par
First consider $\bE_fS_1$. By (\ref{eq:oracle-inequality-1})
for any $\theta\in (0,1]$ we have
\begin{eqnarray}
S_1 &=& \int_{D_1} |\hat{f}(x)-f(x)|^p \rd x
\;\leq\;
c_0^{p-\theta} \int_{D_1} \big[\bar{U}_f(x)+\omega(x)\big]^{p-\theta}|\hat{f}(x)-f(x)|^{\theta}\rd x
\nonumber
\\
&\leq&  c_{11}\bigg\{\delta^{p-\theta}\int_{\bR^d} |\hat{f}(x)-f(x)|^{\theta}\rd x
+
\int_{\bR^d} \omega^{p-\theta}(x) [\delta + \omega(x)]^{\theta} \rd x\bigg\}.
\nonumber
\end{eqnarray}
Here we have used that, by (\ref{eq2:prop-new}),
$\bar{U}_f(x)\leq 2 c_{10}\delta$ for all $x\in D_1$.
Remind that $\hat{f}(x)=\hat{f}_{\hat{h}(x)}(x)$;
therefore, for any $\theta\in (0,1]$
\begin{eqnarray}
&&\bE_f|\hat{f}(x)|^{\theta}\leq \big(\bE_f|\hat{f}(x)|\big)^{\theta}
\;\leq\;
\Big(\sum_{h\in\cH}\bE_f\big|\hat{f}_h(x)\big|\Big)^{\theta}\leq c_{12}\big[(\ln{n})^{d}f^*(x)\big]^\theta.
\nonumber
\end{eqnarray}
Thus, for any $f\in \bG_\theta(R)$,
\[
\delta^{p-\theta} \bE_f \int_{\bR^d} |\hat{f}(x)-f(x)|^{\theta}\rd x
 \;\leq\; \delta^{p-\theta}\Big\{\|f\|_\theta^\theta +
c_{12}(\ln{n})^{d\theta} \|f^*\|_\theta^\theta\Big\} \;\leq\; c_{13} \delta^{p-\theta}R^\theta
(\ln{n})^{d\theta}.
\]
Furthermore, because $\kappa= \rk_\infty^2 [(4d+2)p+4(d+1)]$, by
the second statement of Theorem~\ref{th:oracle-inequality}
\begin{equation*}
%\label{eq6:prop-new}
 \bE_f \int_{\bR^d} \omega^{p-\theta}(x) [\delta + \omega(x)]^{\theta} \rd x
 \leq c_{14}  \delta^{\theta} n^{-(p-\theta)/2} + c_{15} n^{-p/2}\leq c_{15} n^{-p/2}.
\end{equation*}
Combining the last two inequalities we obtain
\begin{equation}\label{eq:S-1}
\bE_f S_1 =
 \bE_f \int_{D_1} |\hat{f}(x)-f(x)|^p \rd x \leq  c_{16} \big[\delta^{p-\theta}R^\theta
(\ln{n})^{d\theta} + n^{-p/2}\big].
\end{equation}
\par
Now we proceed with bounding $\bE_f S_2$.
We have
\begin{eqnarray}
\bE_f S_2&=& \bE_f \int_{D_2} |\hat{f}(x)-f(x)|^p \rd x
\;\leq\;
c_0^{p} \;\bE_f \int_{D_2} \big[\bar{U}_f(x)+\omega(x)\big]^{p}\rd x
\nonumber
\\
& \stackrel{{\rm (a)}}{\leq}&
c_{17} \Big[ (2^{m_0(\theta)}\varphi)^{p-\theta} \int_{D_2} |\bar{U}_f(x)|^{\theta}\rd x
+ n^{-p/2}\Big]
\nonumber
\\
&\stackrel{{\rm (b)}}{\leq}&
c_{17} \big[ (2^{m_0(\theta)}\varphi)^{p-\theta}R^\theta
+ n^{-p/2}\big]\;\leq\;
c_{18} \big[ R^\theta (L_\beta\delta)^{\frac{p-\theta}{1-\theta/s+1/\beta}}
+ n^{-p/2}\big].
\label{eq:S-2}
\end{eqnarray}
Here (a)
follows from the second statement
of Theorem~\ref{th:oracle-inequality} and $\bar{U}_f(x) \leq 2^{m_0(\theta)}\varphi$ for $x\in D_2$, and
(b)
is valid because $\bar{U}_f(x)\leq 3c_{10}f^*(x)$ for all $x\in D_2$,
see (\ref{eq2:prop-new}).
\par
Note that
$\delta^{1-\frac{1}{1-\theta/s+1/\beta }}(\ln{n})^{\frac{d\theta}{p-\theta}}\to 0$ as
$n\to \infty$ since $\theta\leq 1$ and $1/\beta>1/s$, where the latter inequality follows from $r\in (1,\infty]^d$.

Thus, combining (\ref{eq:S-1}) and (\ref{eq:S-2}) with (\ref{eq:S1+S2}),
we obtain
\[
 J_{m_0(\theta)}^-\;\leq\; c_{19}\Big[\big(L_\beta \delta\big)^{\frac{p-\theta}{1-\theta/s+1/\beta }}+
n^{-p/2}\Big]\leq c_{20}\big(L_\beta \delta\big)^{p\nu(\theta)},
\]
as claimed.
\epr

\par\bigskip

\bibliographystyle{agsm}

\end{document}